\documentclass[12pt,leqno]{article}

\usepackage{amsmath,amssymb,theorem,latexsym}
\usepackage[refpage,text]{glosstex}

\usepackage[all]{xy}

\xyoption{2cell}
\UseAllTwocells

\makeatletter
\let\@internalcite\cite
\def\cite{\def\citeauthoryear##1##2{##1{##2}}\@internalcite}
\def\@biblabel#1{\def\citeauthoryear##1##2{##1{##2}}[#1]\hfill}
\makeatother

\newcommand{\Appendix}[1]{%
  \refstepcounter{section}%
  \addcontentsline{toc}{section}%
    {\bfseries\appendixname~\thesection\ #1}%
    {\medskip\noindent \Large\bfseries\appendixname\ \thesection\ #1}%
\sectionmark{#1}\smallskip\noindent
\renewcommand{\theequation}{{\bf 
{{\thesection}.\arabic{subsection}}.{\arabic{equation}}}}
}

\theoremstyle{plain} 

\newtheorem{thm}{Theorem}
\newtheorem{theo}{Theorem}[section]
\newtheorem{lem}[theo]{Lemma}
\newtheorem{cor}[theo]{Corollary}
\newtheorem{defi}[theo]{Definition}
\newtheorem{prop}[theo]{Proposition}
\newtheorem{prop-defi}[theo]{Proposition-Definition}
\newtheorem{lemma-defi}[theo]{Lemma-Definition}
\newtheorem{que}[theo]{Question}
{\theorembodyfont{\rmfamily} \newtheorem{rem}[theo]{Remark}}
{\theorembodyfont{\rmfamily} \newtheorem{ex}[theo]{Example}}
{\theorembodyfont{\rmfamily}
\newtheorem{ex-speculation}[theo]{Example-Speculation}}

\newcounter{inner}
\newcounter{rom}

\renewcommand{\theequation}{{\bf 
{\arabic{section}.\arabic{subsection}}.{\arabic{equation}}}}

\newcommand{\thick}[1]{{\mathfrak #1}}
\newcommand{\mycal}[1]{{\mathcal #1}}
\newcommand{\op}[1]{\operatorname{#1}}
\newcommand{\lm}[1]{\mycal{M_{\! D}^{\;\, l}}(#1)}
\newcommand{\hb}[3]{H^{#1}_{\op{B}}(#2,#3)}
\newcommand{\hdr}[3]{H^{#1}_{\op{DR}}(#2,#3)}
\newcommand{\hdol}[3]{H^{#1}_{\op{Dol}}(#2,#3)}
\newcommand{\hhod}[3]{H^{#1}_{\op{Hod}}(#2,#3)}
\newcommand{\hdel}[3]{H^{#1}_{\op{Del}}(#2,#3)}

\newcommand{\stackb}[2]{\mycal{M}_{\op{B}}(#1,#2)}
\newcommand{\stackdr}[2]{\mycal{M}_{\op{DR}}(#1,#2)}
\newcommand{\stackdol}[2]{\mycal{M}_{\op{Dol}}(#1,#2)}
\newcommand{\stackhod}[2]{\mycal{M}_{\op{Hod}}(#1,#2)}

\newcommand{\stackdrreg}[2]{\mycal{M}^{o}_{\op{DR}}(#1,#2)}
\newcommand{\stackdolreg}[2]{\mycal{M}^{o}_{\op{Dol}}(#1,#2)}
\newcommand{\stackhodreg}[2]{\mycal{M}^{o}_{\op{Hod}}(#1,#2)}

\newcommand{\rb}[3]{R_{\op{B}}(#1,#2,#3)}
\newcommand{\rdr}[3]{R_{\op{DR}}(#1,#2,#3)}
\newcommand{\rdol}[3]{R_{\op{Dol}}(#1,#2,#3)}

\newcommand{\pdr}{\pi_{\op{DR}}}
\newcommand{\pdol}{\pi_{\op{Dol}}}
\newcommand{\phod}{\pi_{\op{Hod}}}

\newcommand{\fdr}{f_{\op{DR}}}
\newcommand{\fdol}{f_{\op{Dol}}}
\newcommand{\GL}[2]{\op{GL}_{#1}(#2)}
\newcommand{\bun}[2]{\mycal{B}un(#1,#2)}
\newcommand{\bunreg}[2]{\mycal{B}un^{o}(#1,#2)}
\newcommand{\vint}{\mycal{V}_{{\mathbb Z}}}
\newcommand{\vc}{\mycal{V}_{\bbc}}
\newcommand{\lieg}[2]{{\mathfrak g}_{#1}(#2)}
\newcommand{\liea}[2]{\mycal{E}_{#1}(#2)}
\newcommand{\lieat}[2]{\mycal{E}^{\sim}_{#1}(#2)}
\newcommand{\ad}[1]{\op{ad}_{#1}}
\newcommand{\del}{\partial}
\newcommand{\delbar}{\bar{\partial}}
\newcommand{\ellt}{\tilde{\ell}}
\newcommand{\nablaf}{\nabla_{\!\! f}}
\newcommand{\hh}{h^{1}/h^{0}}
\newcommand{\pr}{\op{pr}}
\newcommand{\ver}{\op{ver}}

\newcommand{\bbc}{{\mathbb C}}
\newcommand{\bbz}{{\mathbb Z}}
\newcommand{\bbq}{{\mathbb Q}}
\newcommand{\bbh}{{\mathbb H}}
\newcommand{\bbp}{{\mathbb P}}
\newcommand{\bbr}{{\mathbb R}}
\newcommand{\bbv}{{\mathbb V}}
\newcommand{\cstar}{{\mathbb C}^{\times}}

\newlength{\dgrmone}
\newlength{\dgrmtwo}

\begin{document}

\title{Nonabelian $(p,p)$ classes}
\author{L. Katzarkov\thanks{Partially supported by NSF Grant DMS-
9700605 and A.P. Sloan research fellowship} \and
T. Pantev\thanks{Partially supported by NSF Grant DMS-9800790}} 
\date{}
\maketitle

\begin{abstract}
In this paper we generalize to the non-abelian context a 
classical theorem of Griffiths which studies the behavior of the
$(p,q)$-components of a horizontal section in a variation of Hodge
structures over a smooth projective variety. 
\end{abstract}

\tableofcontents


\section{Introduction} \label{sec-intro}
In this paper we generalize to the non-abelian context a 
classical theorem of Griffiths \cite[Theorem~7.1]{griffiths-periods3} 
which studies the behavior of the
$(p,q)$-components of a horizontal section in a variation of Hodge
structures over some smooth projective variety. We begin by recalling
the abelian story.

\subsection{Abelian $(p,p)$ classes} \label{ssec-abelian-pp}

One of the geometric motivations of  the question addressed by
Griffiths  theorem of the
$(p,q)$ classes comes from the Hodge conjecture. Recall that for
any complex smooth projective variety $X$  the Hodge theorem provides a
canonical isomorphism (see
Section~\ref{ssec-abelian-non-abelian-hodge} for details):
\[
\oplus_{p+q = w} H^{p}(X,\Omega^{q}_{X}) \cong H^{w}(X,\bbc).
\]
The set $\mycal{W}^{p}$ 
consisting of classes of
algebraic cycles of codimension $p$ is contained in the locus of
classes of type $(p,p)$, i.e we have an inclusion
\[
\mycal{W}^{p}\subset H^{p}(X,\Omega^{p}_{X})\cap H^{2p}(X,\bbz)
\subset H^{2p}(X,\bbc).
\]
The rational Hodge conjecture then says that the locus
$\mycal{W}^{p}\otimes \bbq$ of motivic cohomology classes coincides
with the locus $H^{p}(X,\Omega^{p}_{X})\cap H^{2p}(X,\bbq)$ of Hodge
classes. 

It is therefore clear that understanding the Hodge conjecture requires
having a good grasp on the locus of Hodge classes
or more modestly on the the locus $H^{p}(X,\Omega^{p}_{X})\cap
H^{2p}(X,\bbc)$ of $(p,p)$-classes.

One way to gauge the plausibility of the Hodge conjecture is by
deforming $X$. A beautiful geometric method for producing potential
counter-examples to the Hodge conjecture was proposed by Andr\'{e} Weil
\cite{weil}. It is based on the observation that subvarieties tend to
be much more rigid than the cohomology classes they represent. Weil
suggests that one looks for a smooth
variety $X$ and a subvariety $Y \subset X$ with
a cohomology class $\alpha = \op{PD}[Y] \in
H^{p}(X,\Omega^{p}_{X})\cap H^{2p}(X,\bbq)$ such that:
\begin{itemize}
\item the pair $(X,
\alpha)$ deforms to a pair $(f : {\mathfrak X} \to S, a : S \to
R^{2p}f_{*}\bbq)$ with $a \in \Gamma(S, R^{p}f_{*}\Omega^{p}_{\mathcal{X}/S})$
and $a(o) = \alpha$ for some $o \in S$.
\item the deformations of $Y$ in ${\mathfrak X}$ are obstructed.
\end{itemize}
In this setup the classes $\{a(s)\}_{s \neq o}$ are candidates for
counter-examples to the Hodge conjecture.

\begin{rem} \label{rem-cschoen} In fact 
Weil even proposes an
explicit construction of pairs $X \subset Y$ coming from
correspondences on abelian varieties of CM-type. However C. Schoen had
shown \cite{chad-schoen} that the deformed classes in Weil's examples
are also represented by algebraic cycles.
\end{rem}

\bigskip

\noindent
Weil's observation shows that it is important to find ways of
deforming cohomology classes so that properties like Hodge type and
rationality are preserved.

On the other hand, there is an intrinsic way to deform cohomology
classes - a parallel transport via the Gauss-Manin connection. Indeed
given a smooth projective morphism $f : X \to S$ between smooth
varieties 
\glosstex(general)[p]{basic-morphism}
the vector bundles $(R^{m}f_{*}\bbc_{X})\otimes_{\bbc}
{\mathcal O}_{S}$ have a natural integrable connection - {\em the
Gauss-Manin connection} - coming from the inclusion
\[
R^{m}f_{*}\bbz_{X}/\op{torsion} \subset R^{m}f_{*}\bbc_{X}
\]
(see
Section~\ref{ssec-abelian-non-abelian-hodge} for details). Thus, given
a cohomology class $\alpha \in H^{2p}(X_{o},\bbc)$ we have a 
unique horizontal extension $a \in \Gamma(U, R^{m}f_{*}\bbc_{X})$ of
$\alpha$ over some simply-connected analytic neighborhood $o \in U \subset S$.
If the $\alpha$ we started with was a rational cohomology class, then
by the definition of the Gauss-Manin connection we will have $a \in
\Gamma(U, R^{m}f_{*}\bbq_{X})$. So, as far as deformations are
concerned, the rationality of a cohomology class does not pose any
major obstructions. The next question to ask is if the Gauss-Manin
connection preserves the property of being of Hodge type $(p,p)$. This
is a non-trivial question since the images of
$R^{p}f_{*}\Omega^{q}_{X}$ in $(R^{p+q}f_{*}\bbc_{X})\otimes_{\bbc}
{\mathcal O}_{S}$ are not algebraic subsheaves and are not preserved
by the Gauss-Manin connection. 

The answer to this question is given by the following theorem (see
Definition~\ref{def-zvhs} for the definition of a variation of Hodge 
structures):

\bigskip

\noindent
{\bf Theorem(\cite[Theorem~7.1]{griffiths-periods3})} {\em Let $S$ be a
smooth complex projective variety and let $(V_{\bbz}\subset
V,\nabla,F^{\bullet})$ be a variation of Hodge structures on $S$. Then
for any horizontal global section $a$ of $V$, the Hodge
$(p,q)$-components of $a$ are also horizontal.

In particular if $a$ is of pure Hodge type $(p,q)$ at some point, then
$a$ is of pure Hodge type $(p,q)$ everywhere.}

\bigskip

\noindent
As a consequence one has the following immediate

\bigskip

\noindent
{\bf Corollary} {\em Let $f : X \to S$ be a smooth morphism of smooth
projective varieties. Let $a \in \Gamma(S,
R^{2p}f_{*}\bbc_{X}\otimes {\mathcal O}_{S})$ an algebraic  global
section which is horizontal with respect to the Gauss-Manin connection
and such that $a(o)$ is of Hodge type $(p,p)$ for some point $o \in S$.
Then $a(s)$ is of Hodge type $(p,p)$ for every $s$.
}

\bigskip

\noindent
It is exactly this corollary that we wish to generalize to the nonabelian
situation. Before we explain what this means it
is instructive to examine the existing generalizations of Griffiths
theorem in the abelian situation.

The first is Deligne's theorem of the fixed part. Given a morphism $f
: X \to S$ as in the previous corollary consider the sub-local system
${\mathbb V} \subset R^{w}f_{*}\bbc_{X}$ spanned by the global 
horizontal sections. Then Griffiths theorem on the $(p,q)$-classes
is equivalent to saying that the $C^{\infty}$ decomposition of
$(R^{w}f_{*}\bbc_{X})\otimes {\mathcal O}_{S}$ into $(p,q)$-pieces
induces a horizontal $(p,q)$-decomposition of ${\mathbb V}$, i.e. that
${\mathbb V}$ is a sub variation of Hodge structures. One can ask if
this is always the case and this what Deligne's theorem
answers\footnote{This is not exactly the way Deligne states his
theorem but the method of proof gives this slightly stronger result.}:

\bigskip

\noindent
{\bf Theorem(\cite{deligne-hodge2})}  {\em Let $f : X \to S$ be a smooth
projective morphism to a quasi-projective $S$. Let $G$ be the Zariski
closure of the monodromy group of $R^{w}f_{*}\bbc_{X}$ and let $V$ be
a representation of $G$ which is defined over ${\mathbb Q}$. Then
\begin{itemize}
\item[(1)] (see \cite[Section~4.2]{deligne-hodge2}) $G$ is a complex
reductive group.
\item[(2)] (see \cite[Section~4.1]{deligne-hodge2}) The isotypic
component ${\mathbb V} \subset R^{w}f_{*}\bbc_{X}$ of $V$ is a
$\bbq$-sub-variation of Hodge structures.
\end{itemize}
}

\bigskip

\noindent
(Recall that the isotypic component of $V$ is by definition the maximal
sub-local system that corresponds to a direct sum of copies of $V$.)

\bigskip

\noindent
Deligne's theorem of the fixed part was generalized further to the
case of general 
variations of Hodge structures (i.e. variations not necessarily of
geometric origin) by W. Schmid:

\bigskip

\noindent
{\bf Theorem (\cite[Theorem~7.22]{schmid})} {\em Let $S$ be a smooth
complex manifold that can be embedded as a Zariski open subset in a
compact analytic space. Let $G$ be the Zariski
closure of the monodromy group of $R^{w}f_{*}\bbc_{X}$ and let $V$ be
a representation of $G$ which is defined over ${\mathbb Q}$. Then
the isotypic
component ${\mathbb V} \subset R^{w}f_{*}\bbc_{X}$ of $V$ is a
$\bbq$-sub-variation of Hodge structures.
}

\bigskip

\noindent
Clearly the two theorems of Deligne and Schmid quoted 
above specialize to Griffiths theorem of the
$(p,q)$-classes after taking $S$-projective and $V$ - the trivial one
dimensional $G$-module.

\bigskip

\noindent
In a slightly different direction one can forget about the integral
structure of $R^{w}f_{*}\bbc_{X}$ and view it just as a complex
variation of Hodge structures (see Definition~\ref{def-cvhs}). 
In other words we look at the data
$(V,\nabla,F^{\bullet},\overline{F}^{\bullet})$ where $V :=
(R^{w}f_{*}\bbc_{X})\otimes {\mathcal O}_{S}$, $\nabla$ is the
holomorphic part of the Gauss-Manin connection, $F^{p} := \oplus_{r
\geq q} R^{q}f_{*}\Omega^{r}_{X/S}$ and $\overline{F}^{\bullet}$ is
the complex conjugate of $F^{\bullet}$. Since $\nabla$ integrable and
satisfies the
Griffiths transversality condition $\nabla(F^{p}) \subset
F^{p-1}\otimes \Omega^{1}_{S}$ we can  form the associated graded
$(E,\theta) =
(\op{gr}_{F^{\bullet}}(V),\op{gr}_{F^{\bullet}}(\nabla))$. That is 
$E = \oplus_{p+q = w} R^{q}f_{*}\Omega^{p}_{X/S}$ and $\theta : E \to
E\otimes \Omega^{1}_{S}$ is a morphism of coherent sheaves satisfying
$\theta\wedge \theta = 0$. Such pairs $(E,\theta)$ are called {\em
Higgs bundles} and their appearance in the context of variations of
Hodge structures is the starting point of the non-abelian Hodge theory
of Simpson that we will be concerned with. 

Note that if an algebraic section $a : S \to V$ is horizontal with
respect to the non-abelian Gauss-Manin connection, then it will
necessarily lie in the kernel of $\theta$. Actually something more is
true. If we put $\mycal{H}$ for the $C^{\infty}$ bundle underlying
both $V$ and $E$, then for any $C^{\infty}$ section $a$ of $\mycal{H}$
we have (see e.g. the proof of \cite[Theorem~7.22]{schmid}):
\begin{equation} \label{eq-formality}
\left( \begin{minipage}[c]{2in} $a$ is a holomorphic section of $V$
satisfying $\nabla a = 0$ \end{minipage}
\right) 
\Longleftrightarrow \left( \begin{minipage}[c]{2in} $a$ 
is a holomorphic section of $E$
satisfying $\theta a = 0$ \end{minipage}
\right)
\end{equation}
In fact (\ref{eq-formality}) implies Griffiths theorem of the
$(p,q)$-classes. The main point is to observe that the connection 
$D = \delbar_{E} + \nabla + \theta$ is a Hermitian connection on
$\mycal{H}$ corresponding to the Hermitian metric on $\mycal{H}$ that
one obtains by flipping the signs of the polarization on $V$ (see
Definition~\ref{def-cvhs}) on the appropriate $(p,q)$-pieces. In
particular the decomposition of $\mycal{H}$ into $(p,q)$-pieces is
$D$-horizontal and hence if $a \in C^{\infty}(S,\mycal{H})$ is a
solution to the initial value problem $Da = 0$, $a(o) = a_{0}  \in
(R^{q}f_{*}\Omega^{q}_{X/S})$, then $a \in
C^{\infty}(S,R^{q}f_{*}\Omega^{q}_{X/S})$. 

Now $\nabla a = 0$ combined with (\ref{eq-formality}) implies that $Da
= (\delbar_{E} + \nabla + \theta)a = \delbar_{E}a + \nabla a + \theta
a = 0$ and hence one gets the theorem of the $(p,q)$ classes.

\bigskip

\noindent
The statement (\ref{eq-formality}) admits a far reaching
generalization which follows from the higher order K\"{a}hler
identities in non-abelian Hodge theory. More precisely recall that
for a smooth projective variety $S$ the non-abelian Hodge theorem of
Corlette-Simpson provides a correspondence between 
complex reductive local systems on $S$ of rank $n$ and poli-stable
rank $n$ Higgs bundles on $S$ with vanishing rational Chern classes
(see Section~\ref{ssec-abelian-non-abelian-hodge} for details). Now
one has the following formality theorem:

\bigskip

\noindent
{\bf Theorem (\cite[Lemma~2.2]{simpson-higgs})} {\em Let $S$ be a
smooth projective variety and let $(V,\nabla)$ be a complex local
system with reductive monodromy. Let $(E,\theta)$ be the corresponding
poli-stable Higgs bundle. Then there is a canonical isomorphism
\[
H^{\bullet}(S,(V,\nabla)) \cong H^{\bullet}(S,(E,\theta)).
\]
}

\bigskip

\noindent
This general theorem specializes to
(\ref{eq-formality}) when we take $(V,\nabla)$ to be a complex
variation of Hodge structures and use cohomology of degree zero. 

In fact as we will show in Proposition~\ref{prop-smooth-curves} the
proof of the non-abelian version of the theorem of the $(p,p)$ classes
follows formally from the Simpson formality theorem.

\subsection{Statement of the main theorem} \label{ssec-statement-main}

In order to generalize the previous discussion to the non-abelian
setting we have to explain what the non-abelian versions of the
$(p,p)$-classes and the Gauss-Manin connection are. Most of the
relevant concepts were discovered and extensively studied by Simpson
(see Section~\ref{ssec-abelian-non-abelian-hodge}
for a short description and the precise references to the original
works). 

We will discuss only the case of the first non-abelian cohomology since
it is the most geometric one. All essential features of the problem
carry over to the case of the higher degree cohomology 
\cite{simpson-algebraic-nah}, \cite{hirschowitz-simpson-descent} but
this is beyond the scope of the present paper. 

The right way (see Section~\ref{ssec-abelian-non-abelian-hodge} and
references therein) to view the first de Rham cohomology of a smooth 
projective
variety $Y$ is as the moduli (stack) $\stackdr{Y}{n}$ 
of rank $n$ local systems on $Y$. The non-abelian $(p,p)$ classes in
$\stackdr{Y}{n}$ can be defined as the set of fixed points for the
action of certain Weil operators (generalizing the abelian ones) and 
were identified by Simpson \cite[Lemma~4.1]{simpson-higgs} as the
local systems underlying complex variations of Hodge structures.

Next, for a smooth projective morphism $f : X \to S$ with connected fibers
one looks at the the family of de Rham cohomology along the fibers
$\pdr : \stackdr{X/S}{n} \to S$. Similarly to the case of variations
of Hodge structures of geometric origin one can show that $\pdr$ 
carries a natural structure of a local system
of stacks (or a crystal) over $S$. Roughly speaking this means that
locally over $S$ the family $\pdr$ is equipped with a canonical
trivialization over any infinitesimal thickening of $S$ or
equivalently $\stackdr{X/S}{n}$ is equipped with an action of the
sheaf of differential operators $\mycal{D}_{S}$ on $S$.

The algebraic construction of this non-abelian connection was again 
invented
by Simpson \cite{simpson-moduli2} who dubbed it the {\em non-abelian
Gauss-Manin connection}. We discuss this construction at length in
Section~\ref{ssec-geometric-variations}. To make things more explicit
though here we will use the analytic definition of the
non-abelian Gauss-Manin connection which is more down to earth. 

Assume for simplicity that $f : X \to S$ has a section $\xi$.
The geometric Riemann-Hilbert correspondence
\cite[Proposition~7.8]{simpson-moduli2} identifies the underlying
analytic stack of $\stackdr{X/S}{n}$ with the moduli stack
$\stackb{X/S}{n}$ of
$n$-dimensional complex representations of the fundamental groups of
the fibers of $f$. Furthermore since $f$ has a section one can use the
covering homotopy theorem in order to define a monodromy action of
$\pi_{1}(S,o)$ on the fundamental group of the fiber $X_{o}$. By
construction of the monodromy action the elements of $\pi_{1}(S,o)$
act by group automorphisms of $\pi_{1}(X_{o},\xi(o))$ and so one can compose
them with representations $\pi_{1}(X_{o},\xi(o)) \to
\op{GL}_{n}(\bbc)$ to obtain an action of $\pi_{1}(S,o)$ on
$\stackdr{X_{o}}{n}$. This action is precisely the analytic version of
the non-abelian Gauss-Manin connection. In particular an algebraic
section $a : S \to \stackdr{X/S}{n}$ will be horizontal with respect
to the non-abelian Gauss-Manin connection if and only if $a(o)$ is
fixed under $\pi_{1}(S,o)$ (in the sense of group actions on stacks).

With all of this said we are now ready to state our main theorem

\begin{thm} \label{thm-main} Let $f : X \to S$ be a smooth projective
morphism with connected fibers. Assume that $S$ is projective and let
$a : S \to \stackdr{X/S}{n}$ be an algebraic section of $\pdr$ which
is horizontal with respect to the non-abelian Gauss-Manin connection.
If there exists a point $o \in S$ so that $a(o)$ underlies a complex
variation of Hodge structures, then $a(s)$ underlies a complex
variation of Hodge structures for all $s \in S$.
\end{thm}

\

\bigskip
\begin{rem} \label{rem-non-compact-base} It is very natural to ask if
the result of Theorem~\ref{thm-main} will hold for a quasi-projective
base scheme $S$. Due to the lack of the necessary analysis this seems
to be out of reach at the moment. See however
Corollary~\ref{cor-quasi-projective} for a slightly weaker statement.
\end{rem}

\bigskip

\noindent
Geometrically the proof of this theorem boils down to the statement
that if $(F,\nabla)$ is a global local system on $X$ whose restriction
on $X_{o}$ underlies a complex variation of Hodge structures, then the
restriction of $(F,\nabla)$ to any $X_{s}$ underlies a complex
variation of Hodge structures. This fact is the content of
Proposition~\ref{prop-smooth-curves}.  As mentioned before it
follows ultimately  from Simpson's formality theorem
\cite[Lemma~2.2]{simpson-higgs}. To illustrate what is going one we
present a direct argument in the simple (abelian) 
case of rank one local systems.

\begin{ex} \label{ex-main} Let $X$ and $S$ be smooth and projective
varieties and let $f : X \to S$ be a smooth morphism with connected
fibers. Let $(F,\nabla)$ be a rank one complex local system on $X$
such that $(F_{o},\nabla_{o}) := (F,\nabla)_{|X_{o}}$ underlies a
complex variation of Hodge structures. Recall that a rank one local
system underlies a
complex variation of Hodge structures if and only if it is unitary
i.e. if and only if the connection preserves a hermitian metric on the
bundle. In terms of representations of the fundamental group this just
means that the character $\chi$ corresponding to $(F_{o},\nabla_{o})$
is a character $\chi : \pi_{1}(X_{o}) \to U(1) \subset \cstar$. Since
the monodromy action of $\pi_{1}(S,o)$ on $\stackdr{X_{o}}{n}$ is just
the composition of $\chi$ with automorphisms of $\pi_{1}(X_{o})$ it is
clear that for all $\gamma \in \pi_{1}(S,o)$ the characters $\chi$ and
$\gamma^{*}\chi$ will have the same image, i.e. the property of $\chi$
being a $\bbc$VHS is preserved. 

This argument however is completely misleading since in this simple
case the property of $(F_{s},\nabla_{s})$ underlying a $\bbc$VHS is
entirely topological (i.e. the monodromy of $(F_{s},\nabla_{s})$
should be compact) and does not depend on the complex structure on
$X_{s}$. This of course is not true in general so the previous
argument cannot generalize to a proof of
Theorem~\ref{thm-main}. Therefore it will be more helpful to have  
an algebraic argument that
uses the geometric structure on $X$ and $S$. Such an argument is not
hard to find. Indeed - since $F$ has a holomorphic integrable connection it
follows that $c_{1}(F) = 0 \in H^{1}(X,\Omega^{1}_{X})$ and so we can
find \cite[Section~1.2]{griffiths-harris} a flat unitary connection
$D$ on $F$. Let $\theta = D^{1,0} - \nabla$. Then $\alpha$ is a
holomorphic one form on $X$ and the condition that
$(F_{o},\nabla_{o})$ is unitary is equivalent to saying that the image
of $\alpha$ under the restriction map $r_{o} : \Omega^{1}_{X} \to
\Omega^{1}_{X_{o}}$ is zero. From this we would like to deduce that
the image of $\alpha$ under any restriction map $r_{s} : \Omega^{1}_{X} \to
\Omega^{1}_{X_{s}}$ is zero.

More invariantly put $\Omega^{1}_{f}$ for the sheaf of holomorphic one
forms along the fibers of $f$ and let $r : \Omega^{1}_{X} \to
\Omega^{1}_{f}$ be the obvious surjection. Then the global one form
$\alpha$ maps to a section $r(\alpha) = H^{0}(X,\Omega^{1}_{f}) = 
H^{0}(S,f_{*}\Omega^{1}_{f})$ which vanishes at $o$. We want to show
that $r(\alpha) = 0 \in H^{0}(S,f_{*}\Omega^{1}_{f})$. For this
consider the weight one variation of Hodge structures
$((R^{1}f_{*}\bbc_{X})\otimes {\mathcal O}_{S},\op{GM})$ with
$\op{GM}$ being the Gauss-Manin connection. By Griffiths infinitesimal
period relations \cite[Section~9]{griffiths-periods12} the Higgs bundle
corresponding to this variation is 
\[
(E,\theta) := (\op{gr}_{F^{\bullet}}(R^{1}f_{*}\bbc_{X})\otimes
{\mathcal O}_{S},\op{gr}_{F^{\bullet}}\op{GM}) = 
(f_{*}\Omega^{1}_{f}\oplus R^{1}f_{*}{\mathcal O}_{X}, \begin{pmatrix} 0 &
0 \\ c & 0 \end{pmatrix})
\]
where $c : f_{*}\Omega^{1}_{f} \to R^{1}f_{*}{\mathcal O}_{X}\otimes
\Omega^{1}_{S}$ is the cup product with the Kodaira-Spencer class
$\kappa_{X/S} \in R^{1}f_{*}T_{f}\otimes \Omega^{1}_{S}$ of $f$.

Recall next that $\kappa_{X/S}$ is defined as the first edge homomorphism
for the push-forward long exact sequence of the following 
short exact sequence of
sheaves on $X$
\[
0 \to T_{f} \to T_{X} \to f^{*}T_{S} \to 0.
\]
Therefore $c$ is just the first edge homomorphism of the push-forward
long exact sequence for
\[
0 \to f^{*}\Omega^{1}_{S} \to \Omega^{1}_{X} \stackrel{r}{\to}
\Omega^{1}_{f} \to 0
\]
and so $\ker(c) = \op{im}(r)$. This shows that $r(\alpha)$ is a
holomorphic section of $E$ which is annihilated by $\theta$ and so by
(\ref{eq-formality}) $r(\alpha)$ can be interpreted as an algebraic
section of $(R^{1}f_{*}\bbc_{X})\otimes {\mathcal O}_{S}$ which is
$\op{GM}$ horizontal. Hence if $r(\alpha)$ vanishes at one point it
must vanish everywhere.
\end{ex}

\

\bigskip

\noindent
Here is a brief outline of the content of this paper. In
Section~\ref{sec-D} we
recall the notions of a local system and a $\mycal{D}$-module and a
local systems of schemes and a $\mycal{D}$-scheme. In
Section~\ref{sec-nahodge} with the hope of making the paper more
readable we review the necessary background from
non-abelian Hodge theory. The main part of the paper begins in
Section~\ref{sec-GM} where we interpret the non-abelian Gauss-Manin
connection in terms of deformation theory and find an explicit
condition for a section to be
horizontal. Whenever possible we have adopted the \v{C}ech point of
view hoping to make the arguments more transparent.
 Section~\ref{ssec-general-variations} discusses the case
of variations of non-geometric origin and possible directions of
generalizing Theorem~\ref{thm-main}. In Section~\ref{sec-main} we
prove the main theorem first in the case of curves and then in
general. For the convenience of the reader we have collected in an
Appendix some well known fact about algebraic stacks that are used
throughout the paper.

\bigskip

\noindent
The proof of the theorem of the non-abelian $(p,p)$-classes presented
here was finished in late 1996. Since then we have reported
on this work on several occasions explaining details of the
proof. In particular we have given lectures on the subject during the Warwick
Symposium in the Summer of 1997, the Oberwolfach Complex Geometry
Workshop in September 1997, and the ICM Satellite Conference in Essen in
August 1998. We apologize for the delay in writing this result up.

During the preparation of this paper 
C.~Simpson informed us that J. Jost and K. Zuo (who were both present at the
Oberwolfach lecture of L.K.) have announced similar results.

\bigskip

\noindent
{\bf Acknowledgments}  We are very grateful to Carlos Simpson for the
innumerous hours he spent explaining non-abelian Hodge theory to us
and for sharing with us his ideas and insights. 

We would also like to thank Sasha Beilinson, Alexis Kouvidakis 
and Andrew Kresch for many
illuminating discussions on the subject of this paper.

\subsection{Notation and terminology} \label{ssec-notation}

In this section we list the basic notions used throughout the paper and
give page references for the place in the text where they are explained.

\printglosstex(general)[p]
\printglosstex(categories)[p]
\printglosstex(hodge)[p]
\printglosstex(complexes)[p]
\printglosstex(cohomology)[p]
\printglosstex(moduli)[p]

\section{Preliminaries on $\mycal{D}$-varieties} \label{sec-D}
Let $S$ be a smooth scheme over the complex numbers. 
\glosstex(general)[p]{S}
The main subject of our investigation will be certain families 
of schemes and stacks over $S$ that are endowed with an action of the 
tangent sheaf $T_{S}$ of $S$. 
\glosstex(general)[p]{tangent}
We recall some standard facts and
constructions that formalize the notion of a $T_{S}$ action.  
Among the useful references are \cite{bernstein}, \cite{borel} for  
$\mycal{D}$-modules and \cite{beilinson-drinfeld-chiral},
\cite{grothendieck-crystals} \cite{berthelot} \cite{illusie-cotangent}
\cite{simpson-moduli2} and \cite{simpson-santa-cruz} for 
$\mycal{D}$-schemes and crystals.

\subsection{$\mycal{D}$-modules and local systems} 
\label{ssec-D-local}

\begin{defi} \label{def-connection}
Let $F \to S$ be a vector bundle. An algebraic connection 
on $F$ is a $\bbc$-morphism of sheaves $d^{\nabla} : 
F \to F\otimes \Omega_{S}^{1}$ satisfying $d^{\nabla}(f\cdot 
a) = fd^{\nabla}(a) + a\otimes df$ for every $f \in {\mathcal 
O}_{S}$ and $a \in F$. 
\end{defi}

\glosstex(general)[p]{bundle}
\glosstex(general)[p]{connection-as-differential}

The differential $d^{\nabla}$ extends 
to $d^{\nabla} : F\otimes \Omega^{i}_{S} \to F\otimes 
\Omega^{i+1}_{S}$ by the Leibnitz rule:
\[
d^{\nabla}(a\otimes \alpha) = a\otimes d\alpha + (-1)^{i} 
d^{\nabla}a\wedge \alpha.
\]
Equivalently, a connection on
$F$ is a consistent way of lifting infinitesimal symmetries
of $S$ to infinitesimal symmetries of $F$ that are linear along
the fibers. The sheaf of infinitesimal symmetries of $S$ is
the holomorphic tangent bundle $T_{S}$. The sheaf  $\mycal{E}(F) \to
S$
\glosstex(general)[p]{atiyah-algebra}
\glosstex(general)[p]{diff-operators-on-bundle}
of the infinitesimal symmetries of $F$ that are linear along the 
fibers can be described as follows. Put
$\mycal{D}_{S}^{1}(F)$ for the sheaf of differential operators
on $F$ of order $\leq 1$. There is a standard short 
exact sequence of ${\mathcal O}_{S}$-modules
\[
0 \longrightarrow \op{End}(F) \longrightarrow 
\mycal{D}_{S}^{1}(F) \stackrel{\sigma}{\longrightarrow}
\op{End}(F)\otimes T_{S} \longrightarrow 0,
\]
where $\sigma(\del) = [\del, \bullet]$ is the principal
symbol map. Then
\[
\mycal{E}(F) := \left\{ \del \in \mycal{D}^{1}_{S}(F) \left|
\sigma(\del) = \op{id}_{F}\otimes v, \; v \in T_{S} \right. 
\right\}.
\]
The sheaf $\mycal{E}(F)$ has a natural $\bbc$-linear 
Lie bracket $[\del',\del''] = \del'\del'' -
\del''\del'$ and is called the {\em Atiyah algebra}
of $F$. A connection on $F$ is just a ${\mathcal 
O}_{S}$-linear splitting $\nabla$
of the symbol sequence for $\mycal{E}(F)$:
\[
\xymatrix@1{
0 \ar[r] & \op{End}(F) \ar[r] & {\mycal{E}(F)}  
\ar[r]^{\sigma} &  T_{S} \ar[r] \ar@/^1pc/[l]^\nabla  &
 0.}
\]
The existence of an algebraic connection on $F$ is obstructed
by the extension class $e(F) \in H^{1}(S,\op{End}(F)
\otimes \Omega^{1}_{S})$ of this sequence, called the {\em 
Atiyah class} of $F$.

\glosstex(general)[p]{connection-as-splitting}
\glosstex(general)[p]{atiyah-class}
\glosstex(general)[p]{diff-operators}

\begin{defi} \label{def-integrable}
A connection $\nabla$ is called {\em integrable} if it is a 
morphism of sheaves of Lie algebras, i.e. if for any $\xi, 
\eta \in T_{S}$ we have $\nabla_{[\xi,\eta]} = [\nabla_{\xi},
\nabla_{\eta}]$. A vector bundle equipped with an integrable 
connection is called a {\em local system}.
\end{defi}

Equivalently, $\nabla$ is integrable if 
$d^{\nabla}\circ d^{\nabla} = 0$. The algebraic curvature of
$\nabla$ is the endomorphism valued 2-form 
\[
\op{curv}(\nabla) =
d^{\nabla}\circ d^{\nabla} \in H^{0}(S,\op{End}(F)\otimes 
\Omega^{2}_{S}).
\] 
It can be interpreted geometrically as the 
obstruction to lifting the $T_{S}$ action on $F$ to an action 
of the full sheaf $\mycal{D}_{S}$ of differential operators.  
In other words, a connection $\nabla$ is integrable 
if and only if it lifts to an ${\mathcal O}_{S}$-linear 
morphism $\mycal{D}_{S} \to \mycal{D}_{S}(F)$. This implies in 
particular that the coherent sheaf of holomorphic sections
of $F$ is endowed with a left $\mycal{D}_{S}$-action. Sometimes it is 
useful to consider more general objects of this type:

\begin{defi} \label{def-D-module}
A quasi-coherent sheaf on $S$ equipped with a left 
$\mycal{D}_{S}$-action is called a {\em left 
$\mycal{D}_{S}$-module}.
\end{defi}

\noindent
Denote by $\lm{S}$ the category of all left 
$\mycal{D}_{S}$-modules. 

\glosstex(categories)[p]{dmodules}

\begin{rem} \label{rem-D-modules} {\bf (i)} 
One can consider {\em coherent} $\mycal{D}_{S}$-modules.
By definition these are locally finitely
generated $\mycal{D}_{S}$-modules. By an analogue of Oka's
theorem \cite{bernstein}, \cite{borel} the coherent 
$\mycal{D}_{S}$-modules are always locally finitely presented.

There is also a notion of {\em smooth} $\mycal{D}_{S}$-modules.
By definition these are $\mycal{D}_{S}$-modules that are 
coherent as ${\mathcal O}_{S}$-modules. It is a simple 
consequence \cite{bernstein}, \cite{borel} of Kashiwara's lemma that 
the smooth $\mycal{D}_{S}$-modules are precisely the ones 
that are locally free and of finite rank as 
${\mathcal O}_{X}$-modules, i.e. the ones that arise from 
vector bundles with integrable connections.

\medskip

\noindent
{\bf (ii)}
There is a natural forgetful functor 
\[
\lm{S} \stackrel{o}{\longrightarrow} \mycal{M}_{{\mathcal O}}(S)
\] 
to the category $\mycal{M}(S)$ of quasi-coherent ${\mathcal 
O}_{S}$-modules. This functor has a natural left adjoint
functor 
\glosstex(categories)[p]{qcmodules}
\[
\xymatrix@R=1pt{
{\mycal{M}_{{\mathcal O}}(S)} \ar[r]^{\mycal{D}_{S}\otimes} & 
{\lm{S}} \\
F \ar[r] & {\mycal{D}_{S}\otimes_{{\mathcal O}_{S}}F}}
\]
\

\medskip

\noindent
{\bf (iii)} There is another pair of adjoint functors - the first is the 
functor of horizontal sections
\[
\xymatrix@R=1pt{
{\lm{S}} \ar[r]^{\Gamma^{\op{hor}}} & {\op{Vect}_{\bbc}} \\
M \ar[r] & {\op{Hom}_{\lm{S}}({\mathcal O}_{S},M)}}, 
\]
where $\op{Vect}_{\bbc}$ is the category of vector 
spaces over $\bbc$. 
\glosstex(categories)[p]{vector-spaces}
Its left adjoint is the functor assigning to a
vector space the corresponding constant $\mycal{D}$-module, i.e. the
functor
\[
\op{Vect}_{\bbc} \stackrel{\otimes {\mathcal O}_{S}}{\longrightarrow}
\lm{S}.
\] 
For a smooth $\mycal{D}$-module corresponding to a 
local system $(M,\nabla)$ the vector space of horizontal
sections can be identified naturally with the space of all
sections in $M$ that are annihilated by $d^{\nabla}$, i.e. 
$\Gamma^{\op{hor}}((M,\nabla)) = \{ m \in \Gamma(S,M) | 
d^{\nabla}(m) = 0 \}$.

\medskip

\noindent
{\bf (iv)} $\lm{S}$ is an abelian tensor category with a 
tensor product given by 
\[
M_{1}\otimes M_{2} := M_{1}\otimes_{{\mathcal O}_{S}} M_{2}.
\]
Clearly ${\mathcal O}_{S}$ is a unit for $\otimes$ and $o$ is a 
tensor functor. The tensor product $\otimes$ does not preserve
coherency and the category $\lm{S}$ does not have duals. 
Actually, it is easy to check that a $\mycal{D}$-module 
$M$ admits a dual iff $M$ is smooth.
\end{rem}

\subsection{Crystals and local systems of schemes}
\label{ssec-crystals}
The notions of a $\mycal{D}$-module and a smooth 
$\mycal{D}$-module readily generalize to families of schemes 
or stacks. In this section we gather some definitions and 
basic facts about those for further reference.

\begin{defi} \label{def-connection-schemes} Let $f : X \to S$
be a smooth morphism of smooth schemes. An algebraic connection
on $X/S$ is a consistent way of lifting of infinitesimal 
automorphisms of $S$ to infinitesimal automorphisms of $X$, 
i.e. an ${\mathcal O}_{X}$-splitting $\nabla$ of the exact 
sequence of vector bundles on $S$:
\begin{equation}
\xymatrix@1{
0 \ar[r] & T_{f} \ar[r] & T_{X}  
\ar[r]^{df} &  f^{*}T_{S} \ar[r] \ar@/^1pc/[l]^\nabla  &
 0.} \label{eq-tangent-X/S}
\end{equation}
where $T_{X}$ and $T_{S}$ are the tangent sheaves of $X$ and 
$S$ respectively and $T_{f}$ is the sheaf of germs of vector 
fields on $X$ that are tangent to the fibers of $f$.
\end{defi}

The existence of an algebraic connection on $X/S$ is obstructed
by the extension class $e(X/S) \in H^{1}(X,T_{f}\otimes 
f^{*}\Omega^{1}_{S})$ which we will call again the {\em Atiyah 
class} of $X/S$.
\glosstex(general)[p]{atiyah-class-for-schemes}

\begin{rem} \label{rem-atiyah-vs-kodaira-spencer}
If $f$ is proper and with connected fibers that do not have
infinitesimal automorphisms, then the Atiyah class $e(X/S)$
is essentially the Kodaira-Spencer class of $f : X \to S$. 
Recall that the (naive) Kodaira-Spencer class $\kappa_{X/S}$ of a 
family $f : X \to S$ is the first edge homomorphism of the
\glosstex(general)[p]{naive-kodaira-spencer}
direct image sequence of (\ref{eq-tangent-X/S}), i.e.
\[
\xymatrix{
0 \ar[r] & f_{*}T_{f} \ar[r] & f_{*}T_{X}  
\ar[r]^{df} & f_{*} f^{*}T_{S} \ar[r]^{\kappa_{X/S}}  & 
R^{1}f_{*}T_{f}.}
\] 
Since by assumption $f_{*}{\mathcal O}_{X} = {\mathcal O}_{S}$ we have that 
$f_{*}f^{*}T_{S} = T_{S}$ and thus $\kappa_{X/S}$ can be viewed
as an element in $H^{0}(S,R^{1}f_{*}T_{f}\otimes 
\Omega^{1}_{S})$. On the other hand, the Lerray spectral 
sequence gives
\[
\xymatrix{
0 \ar[r] & H^{1}(f_{*}T_{f}\otimes \Omega^{1}_{S}) \ar[r] &
H^{1}(X,T_{f}\otimes f^{*}\Omega^{1}_{S}) \ar[r] &
H^{0}(R^{1}f_{*}T_{f}\otimes \Omega^{1}_{S})),
}
\]
and due to the fact that $f_{*}T_{f} = 0$ we get an
inclusion 
\[H^{1}(X,T_{f}\otimes f^{*}\Omega^{1}_{S})\hookrightarrow 
H^{0}(S,R^{1}f_{*}T_{f}\otimes \Omega^{1}_{S}))
\]
under which $e(X/S)$ goes to $\kappa_{X/S}$.
\end{rem}

\bigskip

Denote by the subsheaf $f^{-1}{\mathcal O}_{S} \subset 
{\mathcal O}_{X}$ consisting of germs of functions that are 
constant along the fibers. Let $f^{-1}T_{S}$ be the sheaf theoretic inverse 
image of the sheaf $T_{S}$ and let $T^{\sim}_{X} \subset
T_{X}$ be the centralizer of $T_{f}$ with respect to the Lie 
bracket on $T_{X}$. 
\glosstex(general)[p]{vertical-tangent}
\glosstex(general)[p]{constant-on-fibers-tangent}
There is a natural 
exact sequence of sheaves on $X$:
\begin{equation}
\xymatrix{
0 \ar[r] & T_{f} \ar[r] & T^{\sim}_{X}  
\ar[r]^{df} &  f^{-1}T_{S} \ar[r]  & 0.} 
\label{eq-tilde-tangent-X/S}
\end{equation}
which is just the pull-back of the extension 
(\ref{eq-tangent-X/S}) via the inclusion   
$f^{-1}T_{S} 
\hookrightarrow f^{*}T_{S}$. Moreover, since $T_{f} 
\triangleleft T^{\sim}_{X}$ is an ideal with respect to 
the $\bbc$-linear Lie bracket we get a natural bracket
on $f^{-1}T_{S}$ and the exact sequence 
(\ref{eq-tilde-tangent-X/S}) becomes a sequence of Lie algebra
sheaves.  Under the assumptions
of Remark~\ref{rem-atiyah-vs-kodaira-spencer} it is clear that
the first edge homomorphism in the direct image of (\ref{eq-tilde-tangent-X/S})
is precisely the Kodaira-Spencer class of $f : X \to S$.
Moreover if $\nabla$ is a connection on 
$X/S$, then its restriction $\nabla_{| f^{-1}T_{S}} :  f^{-1}T_{S} \to 
T^{\sim}_{X}$ is $f^{-1}{\mathcal O}_{S}$-linear. Conversely, any 
$f^{-1}{\mathcal O}_{S}$-linear splitting of the sequence 
(\ref{eq-tilde-tangent-X/S}) determines a connection on $X/S$ when extended
by ${\mathcal O}_{X}$-linearity.

\begin{defi} \label{def-integrable-schemes} A connection 
$\nabla$ on $X/S$ is called {\em integrable} if it is a morphism
of Lie algebra sheaves when restricted on $f^{-1}T_{S}$, i.e. 
if for any $\xi, \eta \in f^{-1}T_{S}$ we have 
$\nabla_{[\xi,\eta]} = [\nabla_{\xi},\nabla_{\eta}]$. A smooth
scheme $X/S$ equipped with an integrable connection is called
a local system of schemes.
\end{defi}

As in the vector bundle case the integrability of a connection on $X/S$
ensures that the infinitesimal $f^{-1}T_{S}$ action on $X$ given by $\nabla$
lifts to a morphism of rings of differential operators $f^{-1}\mycal{D}_{S}
\to \mycal{D}^{\sim}_{X}$ where $ \mycal{D}^{\sim}_{X}$ is the
centralizer of the sheaf of vertical differential operators $\mycal{D}_{f}$
in $\mycal{D}_{X}$. It is again important to consider more general objects
with such an action - the $\mycal{D}_{S}$-schemes.
\glosstex(general)[p]{vertical-differential-operators}
\glosstex(general)[p]{constant-on-fibers-differential-operators}

\medskip

\begin{defi} \label{def-D-algebra} A $\mycal{D}_{S}$-{\em algebra} is an 
associative commutative unital algebra $A$ in the category $\lm{S}$. In 
other words it is a $\mycal{D}_{S}$-module equipped with a horizontal 
associative commutative product and a horizontal section which is a unit 
for this product. The affine $X$-scheme $\op{Spec}(A)$ is called an {\em 
affine  $\mycal{D}_{S}$-scheme}.
\end{defi}

Denote by $\op{Comm}_{{\mathcal O}}(S)$ the category of unital associative 
commutative ${\mathcal O}$-algebras (= category of affine  ${\mathcal 
O}$-schemes  and by  $\op{Comm}_{\mycal{D}}(S)$ the category of associative 
commutative $\mycal{D}_{S}$-algebras (= category of affine  
$\mycal{D}_{S}$-schemes).

\glosstex(categories)[p]{commutative-O-algebras}
\glosstex(categories)[p]{commutative-D-algebras}

\begin{rem} \label{rem-D-algebras}
{\bf (i)} Both  $\op{Comm}_{{\mathcal O}}(S)$ and  $\op{Comm}_{\mycal{D}}(S)$
are tensor categories w.r.t. the tensor product $A\otimes B := A
\otimes_{{\mathcal O}_{S}} B$. Moreover $\otimes$ is a coproduct
and the natural forgetful functor 
\[
\op{Comm}_{{\mathcal O}}(S) \stackrel{o}{\longrightarrow} 
\op{Comm}_{\mycal{D}}(S)
\]   
commutes with $\otimes$.

\medskip

\noindent 
{\bf (ii)} A large supply of $\mycal{D}_{S}$-algebras is provided by the
jet construction. To any $R \in \op{Comm}_{{\mathcal O}}(S)$ one can 
associate a canonical $\mycal{D}_{S}$ algebra $JR$ called the {\em jet 
algebra} of $R$ \cite[$4_{\op{IV}}$ 16]{egaiv} by setting
\[
JR := \op{Sym}\left(\mycal{D}_{S}\otimes_{{\mathcal O}_{S}}R\right)\left/
\left(
\begin{minipage}[c]{2.4in}
the $\mycal{D}_{S}$-ideal generated by $1_{R} -1$ and $1\otimes r_{1}r_{2} -
(1\otimes r_{1})(1\otimes r_{2})$ for all $r_{1}, r_{2} \in R$
\end{minipage}
\right).\right. 
\]
The functor 
\[
\op{Comm}_{{\mathcal O}}(S) \stackrel{J}{\longrightarrow} 
\op{Comm}_{\mycal{D}}(S)
\]
is called the {\em jet functor} and also commutes with the tensor structure.
It is left adjoint to the forgetful functor $o$.
Geometrically $JR$ can be interpreted as follows. Denote by $S^{(n)}$ the
$n$-th infinitesimal neighborhood of the diagonal $S
\stackrel{\Delta}{\hookrightarrow} S\times S$. In other words $S^{(n)}$ is a
nilpotent scheme supported on $S$ with ${\mathcal O}_{S^{(n)}} :=  
{\mathcal O}_{S\times S}/I^{n+1}$ where $I$ is the ideal of $\Delta(S)$ in
$S\times S$. The structure sheaf ${\mathcal O}_{S^{(n)}}$ carries two 
${\mathcal O}_{S}$-structures coming from the two embeddings of 
${\mathcal O}_{S}$ in ${\mathcal O}_{S\times S}$. Let $B$ be an 
${\mathcal O}_{S}$-algebra, then $B$-points of the formal scheme $\op{Spec}
(oJB)$ are given by \cite[$4_{\op{IV}}$ 16]{egaiv}, 
\cite{beilinson-drinfeld-chiral}:
\[
\op{Hom}_{{\mathcal O}_{S}}(oJR,B) = \lim_{\longleftarrow} 
\op{Hom}_{{\mathcal O}_{S^{(n)}}}({\mathcal O}_{S^{(n)}}
\otimes_{{\mathcal O}_{S}}R, B\otimes_{{\mathcal O}_{S}} {\mathcal 
O}_{S^{(n)}}).
\]
In particular we have 
\[
\bbc-\text{points of} \; \op{Spec}(oJR) = \left\{ (s,\gamma) \left| 
\begin{minipage}[c]{2.5in}
$s \in S$, $\gamma$ is a section of $\op{Spec}(R/S)$ over the complete formal 
neighborhood of $s$ in $S$.
\end{minipage}
\right.\right\}
\] 
or equivalently $\op{Spec}(oJR)$ is the space of infinite jets of sections
of $\op{Spec}(R)$ over $S$.
\end{rem}

\

\bigskip

\noindent
This notions can be globalized (see e.g. \cite{beilinson-drinfeld-chiral},  
\cite{simpson-moduli2}). Denote by $(S\times S)^{\wedge}$ and
$(S\times S\times S)^{\wedge}$ the formal neighborhood of the (small)
diagonals in the second and third Cartesian product of $S$. 

\glosstex(general)[p]{formal-ngbhd}

\begin{defi}\label{def-D-scheme} A {\em $\mycal{D}_{S}$-scheme} or a 
{\em crystal of schemes on $S$} is a scheme $F \to S$ together with an 
isomorphism 
\[
\varphi : (F\times S)_{|(S\times S)^{\wedge}} \widetilde{\longrightarrow} 
 (S\times F)_{|(S\times S)^{\wedge}}
\] 
satisfying the cocycle condition 
\[
p_{23}^{*}(\varphi)p_{12}^{*}(\varphi) = p_{13}^{*}(\varphi)
\]
for the resulting isomorphisms  between the restrictions of $F\times S\times S$
and $S\times S\times F$ on $(S\times S\times S)^{\wedge}$.
\end{defi}

This notion generalizes in an obvious way to the relative situation  
\cite{simpson-moduli2}, and we can talk about about crystals on $X/S$. If 
in addition $F \to S$ is a vector bundle and the identification $\varphi$ is
a morphism of vector bundles we will call $F$ a {\em crystal of vector 
bundles} on $S$. By reducing to the case when $X$ is affine over $S$ it is 
easy to check that a crystal of vector bundles on $X/S$ is the same as a
vector bundle on $X$ with a relative integrable connection over $S$ 
\cite[Lemma~8.1]{simpson-moduli2}. 

\begin{rem} \label{rem-crystal-vs-stratification}
Strictly speaking, according to the standard terminology 
\cite[Appendix]{grothendieck-crystals} the object $F$ from 
Definition~\ref{def-D-scheme} should be called a {\em stratification of 
schemes} rather than a crystal of schemes. However, due to the infinitesimal 
lifting property, if $S$ is smooth, or in the relative case, if $X$ is smooth 
over $S$, the notions of a crystal and stratification coincide 
\cite[Section~8]{simpson-moduli2}. Since we will be working only with smooth
projective families we will suppress the distinction between a stratification
and a crystal. 
\end{rem}

\begin{rem} \label{rem-deRham-groupoid} For technical reasons it is
sometimes more convenient to restate the conditions in
Definition~\ref{def-D-scheme} in terms of the action of a formal
groupoid. This turns out to be extremely useful for understanding
the Griffiths transversality condition for variations of non-abelian
Hodge structures (see \cite[Section~8]{simpson-santa-cruz} and
\cite[Section~9]{simpson-algebraic-nah}.  
 Recall \cite{berthelot} \cite[Chapter~8]{illusie-cotangent}
\cite[Section~7]{simpson-santa-cruz} that a {\em formal groupoid of smooth
type} is a groupoid 
\[
(e : X \to N, \xymatrix@1{
N \ar@<.5ex>[r]^{s}\ar@<-.5ex>[r]_{t}  & X
}, m : N\times_{X} N \to N)
\]  
in the category of formal schemes such that $X$ is a scheme, $e(X)$ is
the topological space underlying the formal scheme $N$, $N$ is locally
the completion of a scheme of finite type and $s$ and $t$ are
formally smooth. 

To any smooth scheme $S$ of finite type one associates the formal
groupoid $S_{\op{DR}} := \xymatrix@1{
(S\times S)^{\wedge} \ar@<.5ex>[r]\ar@<-.5ex>[r]  & S
}$. By abuse of notation we will write $S_{\op{DR}}$ for the
corresponding quotient stack as well. Now it is clear that $F \to S$
is a local system of schemes if and only if $F$ is actually a sheaf on
$S_{\op{DR}}$ or equivalently if and only if $F$ is equipped with and
action of the formal groupoid $\xymatrix@1{
(S\times S)^{\wedge} \ar@<.5ex>[r]\ar@<-.5ex>[r]  & S
}$.
\glosstex(general)[p]{SdR}

\end{rem}

We conclude this section with a discussion of the analogies between the 
linear and non-linear $\mycal{D}_{S}$-objects. 
Denote by $(\op{Sch}/S)$ and $\op{Sch}_{\mycal {D}}(S)$ the
categories of ${\mathcal O}_{S}$ and $\mycal{D}_{S}$-schemes respectively.

\glosstex(general)[p]{D-schemes}

\begin{rem} \label{rem-D-schemes}
{\bf (i)} Since forgetting of the $\mycal{D}_{S}$ action and passing to jets
commute with Zariski and \'{e}tale localization we get a pair of 
adjoint functors
\[
\xymatrix{
\op{Sch}_{\mycal {D}}(S) \ar@<1ex>[r]^{o} & (\op{Sch}/S).
\ar@<1ex>[l]^{J}  
}
\]
In particular for any scheme $X/S$ we get a $\mycal{D}_{S}$-scheme $JX$.

\medskip

\noindent
{\bf (ii)} By analogy with the linear situation we can speak about 
crystals over $S$ that have good regularity properties - i.e. 
$\mycal{D}_{S}$-schemes that are smooth and of finite type as 
${\mathcal O}_{S}$-schemes. It  is almost 
immediate that these are precisely the local systems of schemes over $S$.
 It will be natural to call such a scheme a
smooth $\mycal{D}_{S}$-scheme but we will refrain from that for two reasons.
The first is that we don't want to introduce redundant terminology and the
second is that this term is reserved for another intrinsic notion of the 
smoothness of a crystal.  

We will not be using this intrinsic notion in the sequel but we just comment 
on it in order to avoid confusion and to stress the difference between
local systems 
of schemes and smooth crystals.
To define the latter   \cite{beilinson-drinfeld-chiral}  one
first looks at the affine case. Let $A$ be a $\mycal{D}_{S}$-algebra. To 
understand smoothness we need to understand the sheaf of K\"{a}hler 
differentials in the context of the category of $A$-modules i.e. the 
category of 
quasi-coherent sheaves over $X$ that are endowed with two compatible 
actions of $A$ (as an algebra) and $\mycal{D}_{S}$. This category admits a 
better description.

Consider $A[\mycal{D}_{S}]$ - the sheaf of algebras on $X$ that is determined
uniquely by the properties: (a)  $A[\mycal{D}_{S}]$ is equipped with two 
algebraic embeddings $A \hookrightarrow A[\mycal{D}_{S}]$ and $T_{S} 
\hookrightarrow A[\mycal{D}_{S}]$; (b) as an algebra $A[\mycal{D}_{S}]$ is
generated by $A$ and $T_{S}$ with the only relations being the ones coming 
from the $T_{S}$ action on $A$. By construction the inclusion $T_{S} 
\hookrightarrow A[\mycal{D}_{S}]$ extends to an inclusion of algebras
$\mycal{D}_{S} \hookrightarrow A[\mycal{D}_{S}]$. As an 
$A$-$\mycal{D}_{S}$-bimodule $A[\mycal{D}_{S}]$ is isomorphic to $A
\otimes_{{\mathcal O}_{S}} \mycal{D}_{S}$. 

In the case of an affine
$\mycal{D}_{S}$-scheme which is a local system of schemes, we have already 
encountered  the  sheaf of algebras $A[\mycal{D}_{S}]$. Namely 
$A[\mycal{D}_{S}] = f_{*}\mycal{D}^{\sim}_{\op{Spec}(A)}$, where $f : 
\op{Spec}(A) \to X$ is the structure morphism. 

Notice that an $A$-module in the above sense is just a 
sheaf of $A[\mycal{D}_{S}]$-modules with respect to the natural algebra 
structure. Furthermore one checks that the finitely generated projective
$A[\mycal{D}_{S}]$-modules localize properly \cite{beilinson-drinfeld-chiral}.
Beilinson and Drinfeld define the notions of formal $\mycal{D}_{S}$-smoothness
and $\mycal{D}_{S}$-smoothness
for the affine $\mycal{D}_{S}$-scheme $\op{Spec}(A)$ by requiring a suitable 
infinitesimal lifting property in   the category of $A$-modules. 
Moreover, they 
show that $\op{Spec}(A)$ is $\mycal{D}_{S}$-smooth if and only if
the sheaf of K\"{a}hler differentials $\Omega_{A}$ is
projective as a $A[\mycal{D}_{S}]$-module and $\op{Spec}(A)$ is smooth over
$S$ \cite{beilinson-drinfeld-chiral}. In particular the intrinsic notion of 
smoothness of a crystal $X$ over $S$ requires that $X$ be infinite dimensional
as a scheme over $S$ and thus a local system of schemes of finite type is 
never smooth in this sense.

\medskip

\noindent
{\bf (iii)} Let $\op{Comm}_{\bbc}$ be the category of all commutative
algebras over $\bbc$ (= category of affine schemes over $\bbc$).
Again we have a pair of adjoint functors
\[
\xymatrix{
\op{Comm}_{\bbc} \ar@<1ex>[r]^{\otimes {\mathcal O}_{X}} & 
\op{Comm}_{\mycal{D}}(S).
\ar@<1ex>[l]^{\Gamma^{\op{hor}}}
}
\]
Globally, if $X/S$ is a crystal of schemes the horizontal sections of
$X$ over $S$ are just the algebraic sections $a : S \to X$ for which the 
subscheme $a(S) \subset X$ is a sub-crystal. In the case when $(X/S,\nabla)$ 
is a local system of schemes a sub-local system is a variety $Y \subset X$,
smooth over $S$ and such that  $\nabla$ lifts 
infinitesimal symmetries of $S$ to infinitesimal symmetries of $X$ which
at the points of $Y$ preserve $Y$. In other words
the composition
\[
\xymatrix@1{
{\left(f^{*}T_{S}\right)_{|Y}} \ar[r]^{\nabla} & {\left(T_{X}\right)_{|Y}}
\ar[r] & N_{Y/X}
}
\]
must be identically zero (here $N_{X/Y}$ denotes the normal 
bundle of $Y$ in $X$).
In particular an algebraic section $a : S \to X$ of $f$ will be 
horizontal if and only if the 
following diagram 
\[
\xymatrix@R=8pt{
a^{*}f^{*}T_{S} \ar[r]^{a^{*}\nabla} \ar@{=}[d] & a^{*}T_{X} \\
T_{S} \ar[ur]_{da} &
}
\]
commutes.
\end{rem}

\section{Nonabelian Hodge structures} \label{sec-nahodge}
The non-abelian analogues of the $(p,p)$-classes live in the first 
non-abelian de Rham cohomology spaces of a variety. In this section we review 
briefly Simpson's theory of non-abelian Hodge structures on such spaces. 

\subsection{Abelian and nonabelian Hodge theory}
\label{ssec-abelian-non-abelian-hodge}
The (abelian) cohomology groups of a smooth projective $X$ are endowed with 
extra linear-algebraic data - their Hodge structure. Heuristically the Hodge 
decomposition of $H^{\bullet}(X,\bbc)$ can be thought of  as a linearization 
of the geometry of $X$. Similarly, the Hodge decomposition on the first 
cohomology set of $X$ with coefficients in a non-abelian reductive group 
linearizes $X$ in a sense by replacing it with a family of abelian varieties.

In order to get a better perspective of the setup for our problem we list 
some analogous concepts and facts both in abelian and nonabelian settings. 

\bigskip

\noindent
{\bf A(i)}  For a complex variety $X$ denote by $\hb{\bullet}{X}{\bbc}$ and 
$\hdr{\bullet}{X}{\bbc}$ Betti and de Rham cohomology rings of $X$
\glosstex(cohomology)[p]{abelian-betti}
\glosstex(cohomology)[p]{abelian-dR}
with complex coefficients respectively. The space $\hb{1}{X}{\bbc}$ can also 
be interpreted as the space of representations of the fundamental group of 
$X$ into the additive group $\bbc$. That is
\[
\hb{1}{X}{\bbc}  =
\op{Hom}(\pi_{1}(X),\bbc)= \op{Hom}(H_{1}(X, {\mathbb Z}),\bbc) = 
H^{1}(\pi_{1}(X),\bbc).
\]
The de Rham theorem identifies the vector spaces  
$\hb{\bullet}{X}{\bbc}$ and $\hdr{\bullet}{X}{\bbc}$. Algebraically
it is best to think of the de Rham cohomology of $X$ as the
hypercohomology of the (holomorphic) de Rham complex
\[
\Omega^{\bullet}_{X} := 
\xymatrix@1{{\mathcal O}_{X} \ar[r]^-{d} & \Omega^{1}_{X} \ar[r]^-{d}
& \ldots \ar[r]^-{d} & \Omega^{\dim X}_{X}, 
}
\]
i.e. $\hdr{\bullet}{X}{\bbc} = \bbh^{\bullet}(X, \Omega^{\bullet}_{X})$.
The de Rham theorem can be interpreted in these terms as follows. 
The holomorphic Poincare lemma implies that $\bbc$ (thought as
a complex concentrated in degree zero) and $\Omega^{\bullet}_{X}$ are
quasi-isomorphic and thus have isomorphic hypercohomology. In
particular we get $\hb{\bullet}{X}{\bbc} = 
\bbh^{\bullet}(X, \Omega^{\bullet}_{X})$.

\medskip 

Similarly we can define the Betti and de Rham cohomology of $X$ with 
coefficients in $\cstar$. Thinking of $\bbc$ as the Lie algebra of the
reductive group $\cstar$ and using the exponential map $\exp : \bbc
\to \cstar$ we can ``exponentiate'' the de Rham complex
\[
\xymatrix@R=8pt{\bbc \ar[d]_-{\exp} \ar[r] & {\mathcal O}_{X} \ar[r]^-{d}
\ar[d]_-{\exp} & \Omega^{1}_{X} \ar[r]^-{d} \ar@{=}[d]
& \ldots \ar[r]^-{d} & \Omega^{\dim X}_{X} \ar@{=}[d] \\
\cstar  \ar[r] & {\mathcal O}_{X}^{\times} \ar[r]^-{d\log}
& \Omega^{1}_{X} \ar[r]^-{d} 
& \ldots \ar[r]^-{d} & \Omega^{\dim X}_{X}
}
\]
and so interpret the de Rham cohomology with coefficients in $\cstar$
as the hypercohomology of the logarithmic de Rham complex\footnote{The
slightly odd notation used here is standard in the theory of 
Deligne cohomology.} (here ${\mathcal O}_{X}$ is placed in degree zero)
\[
\bbc(\dim X + 1)[1] := \xymatrix@1{{\mathcal O}_{X}^{\times} \ar[r]^-{d\log}
& \Omega^{1}_{X} \ar[r]^-{d} 
& \ldots \ar[r]^-{d} & \Omega^{\dim X}_{X}.
}
\]
In particular $\hdr{\bullet}{X}{\cstar} = \bbh^{\bullet}(X,\bbc(\dim X
+ 1)[1])$.  A
simple calculation with \v{C}ech cocycles now shows that in degree one
we have
\[
\hb{1}{X}{\bbc^{\times}}  := \op{Hom}(\pi_{1}(X),\cstar) 
= H^{1}(\pi_{1}(X),\cstar), 
\]
and $\hdr{1}{X}{\cstar}$ is the space of algebraic local systems of rank 
one on 
$X$. The analogue of the de Rham theorem in this case is the  abelian 
Riemann-Hilbert correspondence which establishes an analytic isomorphism 
between $\hdr{1}{X}{\cstar}$ and $\hb{1}{X}{\cstar}$. Notice also that the 
exponential map   $\exp : \bbc \rightarrow \cstar$ 
induces a surjective homomorphism from $\hdr{1}{X}{\bbc}$ to the
identity component $\hdr{1}{X}{\cstar}_{0}$ of $\hdr{1}{X}{\cstar}$.

\bigskip

\noindent
{\bf NA(i)}  To any variety $X$ one can associate the spaces  
$\hb{1}{X}{\GL{n}{\bbc}}$
and $\hdr{1}{X}{\GL{n}{\bbc}}$ - the first Betti and  de Rham cohomology of 
$X$ with coefficients in $\GL{n}{\bbc}$. They are defined by analogy
with {\bf A(i)}.
$\hb{1}{X}{\GL{n}{\bbc}}   := {\rm Hom}(\pi_{1}(X), \GL{n}{\bbc})//
\GL{n}{\bbc}$
is the moduli space of semi-simplifications of representations of 
$\pi_{1}(X)$ in $\GL{n}{\bbc}$ and $\hdr{1}{X}{\GL{n}{\bbc}}$ is defined as 
the moduli space of rank $n$ algebraic local systems on 
$X$. The general Riemann-Hilbert correspondence \cite{deligne-riemann-hilbert} 
\cite[Proposition~7.8]{simpson-moduli2} gives an isomorphism
\[
\psi_{X} : \hb{1}{X}{\GL{n}{\bbc}} \longrightarrow \hdr{1}{X}{\GL{n}{\bbc}}
\]
of complex analytic spaces.
\glosstex(hodge)[p]{riemann-hilbert}
\glosstex(cohomology)[p]{nonabelian-betti-space}
\glosstex(cohomology)[p]{nonabelian-dR-space}

\bigskip

\noindent
{\bf A(ii)}  One can also consider the Dolbeault cohomology 
\[
\hdol{w}{X}{\bbc} := \oplus_{p+q = w} H^{p}(X,\Omega^{q}_{X}) =
\bbh^{w}(X, \xymatrix@1{{\mathcal O}_{X} \ar[r]^-{0} & \Omega^{1}_{X}
\ar[r]^-{0} & \ldots \ar[r]^-{0} & \Omega^{\dim X}_{X}})
\] 
of $X$ and its multiplicative version 
\[
\hdol{\bullet}{X}{\cstar} := \bbh^{\bullet}(X, \xymatrix@1{{\mathcal
O}_{X}^{\times} \ar[r]^-{0} & \Omega^{1}_{X}
\ar[r]^-{0} & \ldots \ar[r]^-{0} & \Omega^{\dim X}_{X}}).
\]
Concretely in degree one we have
\[
\hdol{1}{X}{\bbc} = H^{1}({\mathcal O}_{X})\oplus H^{0}(\Omega^{1}_{X})
\]
and
\[
\hdol{1}{X}{\bbc^{\times}} = H^{1}({\mathcal O}_{X}^{\times})\times H^{0}(
\Omega^{1}_{X}) = T^{\vee}H^{1}({\mathcal O}_{X}^{\times}),
\]
respectively.
\glosstex(cohomology)[p]{abelian-dolbeault}

\bigskip

\noindent
{\bf NA(ii)} To any polarized variety $(X,{\mathcal O}_{X}(1))$ one can 
associate  the   first Dolbeault cohomology set 
$\hdol{1}{X}{\GL{n}{\bbc}}$ 
of $X$ with coefficients in $\GL{n}{\bbc}$. It is defined by analogy with 
{\bf A(ii)} as an appropriate moduli space. 

\begin{defi} \label{def-higgs-bundle} A Higgs bundle on $X$ is a pair
$(E,\theta)$ consisting of a vector bundle $E$ of rank $n$  and a homomorphism
$\theta : E \to E\otimes \Omega^{1}_{X}$ satisfying the symmetry condition
$\theta\wedge \theta = 0$. A Higgs bundle $(E,\theta)$ is called 
${\mathcal O}_{X}(1))$-stable (semistable) if for every $\theta$ invariant 
subsheaf $F \subset E$ one has $p(F) < (\leq) p(E)$. Here  
$p(F) := \chi(X,F\otimes {\mathcal O}_{X}(n))/\op{rk}(F)$ is the
reduced Hilbert polynomial of $F$.
\end{defi}

The first Dolbeault cohomology of $X$ with coefficients in $\GL{n}{\bbc}$ 
is by definition the moduli space of semistable rank $n$ Higgs bundles on 
$X$ with vanishing $c_{1}$ and $c_{2}$.
It has a component birationally equivalent to the cotangent bundle
$T^{\vee}H^{1}(X,\GL{n}{{\mathcal O}_{X}})^{\op{reg}}$ to the regular locus 
of the moduli space 
$H^{1}(\GL{n}{{\mathcal O}_{X}})$ of semistable vector bundles of rank $n$ and
trivial $c_{1}$ and $c_{2}$.

\bigskip

\noindent
{\bf A(iii)} There is an equivalence between the Dolbeault and the 
De Rham cohomology which depends only on the class of the chosen K\"{a}hler
metric (polarization).  More precisely we have the following

\medskip

\noindent
{\bf Theorem (Hodge theorem).} 
{\it For any K\"{a}hler $X$ and any $k : 0 \leq k \leq 
\dim_{{\mathbb R}} X$ there is a natural isomorphism
\[
\tau_{X} : \hdr{k}{X}{\bbc} \stackrel{\cong}{\longrightarrow}
\hdol{k}{X}{\bbc} :=  \oplus_{p+q =k} H^{q}(X,\Omega_{X}^{p}).
\]}

\medskip

\noindent
The isomorphism $\tau_{X}$ is built in two steps.  First one shows that both 
de Rham and Dolbeault cohomology classes are represented  by harmonic forms 
and then one uses the K\"{a}hler identities to identify the harmonic
representatives.
\glosstex(general)[p]{hodge}

Furthermore, the exponential map $\exp : \bbc \rightarrow \cstar$
combined with the Hodge theorem gives an isomorphism between the 
multiplicative De Rham and Dolbeault cohomology. In particular there is an 
isomorphism
\[
\tau_{X} : \hdr{1}{X}{\cstar}_{0} \widetilde{\longrightarrow}  
\hdol{1}{X}{\cstar}_{0}.
\]
Explicitly one has the Cartan decomposition $\bbc^{\times} = 
S^{1}\times {\mathbb R}^{+}$ which induces 
\[
\tau_{X} : \hdr{1}{X}{\cstar}_{0} 
\rightarrow H^{1}(X, S^{1})_{0}\times H^{0}(\Omega^{1}_{X}) \cong 
H^{1}({\mathcal O}_{X}^{\times})_{0}\times H^{0}(\Omega^{1}_{X}).
\]
\

\bigskip

\noindent
{\bf NA(iii)}  There is an equivalence between the de Rham and Dolbeault
moduli spaces.

\medskip

\noindent
{\bf Theorem (\cite{corlette-flat}, 
\cite{simpson-higgs}).} {\it For any smooth
projective variety $X$ there is a natural homeomorphism
\[
\tau_{X} : \hdr{1}{X}{\GL{n}{\bbc}} \stackrel{\cong}{\longrightarrow} 
\hdol{1}{X}{\GL{n}{\bbc}}.
\]}

\medskip

Similarly to A(iii) the isomorphism $\tau_{X}$ is built by means of 
harmonic representatives. The latter are the so called harmonic bundles, i.e.
triples $(F,\nabla,h)$ where $F$ is a holomorphic bundle of rank $n$,
$\nabla$ is a flat holomorphic connection on $F$ and $h$ is a hermitian metric 
on $F$ s.t. the corresponding $\pi_{1}(X)$-equivariant map $\tilde{h} :
\widetilde{X} \rightarrow \GL{n}{\bbc}/U(n)$ has minimal energy.

For future reference the isomorphisms $\tau_{X}$ and $\psi_{X}$ can be 
combined to yield an isomorphism
\[
\phi_{X} : \hb{1}{X}{\GL{n}{\bbc}} \longrightarrow \hdol{1}{X}{\GL{n}{\bbc}}.
\]
\

\glosstex(hodge)[p]{betti-dolbeault}

\bigskip

\noindent
{\bf A(iv)} Due to the isomorphism $\hdol{1}{X}{\bbc^{\times}} \cong
T^{\vee}H^{1}({\mathcal O}_{X}^{\times})$ the first Dolbeault cohomology group 
can be viewed
as an algebraic symplectic manifold. Moreover, the Hodge decomposition 
induces two transversal Lagrangian fibrations
\[
\xymatrix{
& {\hdol{1}{X}{\cstar}} \ar[dl] \ar[dr] & \\
H^{1}({\mathcal O}_{X}^{\times}) & & H^{0}(\Omega^{1}_{X})
}
\]
parameterized by the $(0,1)$ and $(1,0)$ parts of the Hodge structure 
respectively.

\bigskip

\noindent
{\bf NA(iv)} The first non-abelian Dolbeault cohomology group is an
algebraic symplectic manifold \cite{biswas-green-lazarsfeld}, 
\cite{donagi-markman}. This 
symplectic structure restricts to the standard one  on 
$T^{\vee}H^1(\GL{n}{{\mathcal O}_{X}})$ and the projection on 
$H^1(\GL{n}{{\mathcal O}_{X}})$ and the 
Hitchin map  give two generically transversal Lagrangian fibrations
\cite{arapura-coho-support2}.
\[
\xymatrix{
& {\hdol{1}{X}{\GL{n}{\bbc}}} \ar[dl] \ar[dr]^{h} & \\
H^1(\GL{n}{{\mathcal O}_{X}}) & & B_{X} := \oplus_{i=1}^{n}H^{0}(S^{i}
\Omega^1_{X})
}
\]
By analogy with {\bf A(iv)} the spaces $H^1(\GL{n}{{\mathcal O}_{X}}$
and $B_{X}$
should be interpreted as the $(0,1)$ and the $(1,0)$ part of the non-abelian
Hodge structure.

\bigskip

\noindent
{\bf A(v)} The locus of Hodge classes in $\hdr{2p}{X}{\bbc}$ is the locus 
of all integral cohomology classes whose image under $\tau_{X}$ is of type 
$(p,p)$. In other words the locus of Hodge classes is just the set
$\vint^{p} := 
\tau_{X}^{-1}(\phi_{X}(\hb{2p}{X}{{\mathbb Z}})\cap 
H^{p}(X,\Omega^{p}_{X}))$.
The Hodge conjecture asserts that up to tensoring by ${\mathbb Q}$ the 
Hodge classes are exactly the classes of algebraic cycles and hence the 
importance of this locus. The first crude approximation to 
the $\vint^{p}$ is obtained by dropping the integrality restriction.
In this way we arrive at the locus of all $(p,p)$ classes: $\vc^{p} := 
\tau_{X}^{-1}(H^{p}(X,\Omega^{p}_{X}))$.
\glosstex(hodge)[p]{abelian-hodge-classes}
\glosstex(hodge)[p]{abelian-pp-classes}

\medskip

Sometimes it is convenient to describe the (formal) loci of Hodge type as 
invariants of certain group of symmetries of the Hodge structure.
The {\em Hodge group} $\textbf{Hod}_{2p}(X)$ is a subgroup of $\op{GL}(
\hb{2p}{X}{{\mathbb Q}})$
\glosstex(hodge)[p]{hodge-group}
such that
\[
\vint^{p} \otimes {\mathbb Q} = \hb{2p}{X}{{\mathbb 
Q}}^{\textbf{Hod}_{2p}(X)}.
\]
It can be defined as follows. The natural rescaling action of $U(1)$ on the
cotangent bundle of $X$ induces action on the $(p,q)$ forms and trough the
Hodge decomposition and the de Rham theorem an action on 
$\hb{2p}{X}{{\mathbb Q}}\otimes \bbc = \hb{2p}{X}{\bbc}$. This action can be 
extended to a homomorphism
\[
C : \cstar \longrightarrow \op{GL}(\hb{2p}{X}{\bbc}). 
\]
Explicitly an element $t \in \cstar$ acts on
the piece $H^{p}(X,\Omega^{q}_{X})$ as multiplication by $C_t = 
t^{q-p}$ (Weil operators).

The Hodge group $\textbf{Hod}_{2p}(X)$ then is defined as the smallest 
subgroup in $\op{GL}(\hb{2p}{X}{\bbc})$ that is defined over ${\mathbb
Q}$ and contains $C(\cstar)$. It is not hard to see then that

\begin{align*}
\vc^{p} & = \hb{2p}{X}{\bbc}^{\cstar} \\
\vint^{p} \otimes {\mathbb Q}& = 
\hb{2p}{X}{{\mathbb Q}}^{\textbf{Hod}_{2p}(X)}
\end{align*}

\bigskip

\noindent
{\bf NA(v)} Recall first the following

\begin{defi} \label{def-zvhs}
A {\em polarized integral variation of Hodge structures} on $X$ consists of 
the following data: a) local system $(V,\nabla)$ on $X$; b) a
lattice $V_{{\mathbb Z}} \subset V$; c)  a finite decreasing filtration  
$F^{\bullet} : \ldots \subseteq F^{i} \subseteq F^{i-1} \subseteq \ldots 
\subseteq V $ such that $\overline{F}^{\bullet}$ is the complex conjugate 
filtration to $F^{\bullet}$ (w.r.t. the real structure coming from the 
inclusion of $V_{{\mathbb Z}}$ in $V$), then there exists a 
number $w$ for which $V = F^{i} \oplus \overline{F}^{w-i+1}$ for all $i$;
d) a horizontal hermitian form $\psi$ on $V$.  Furthermore, these objects 
should satisfy the following axioms: 
\begin{list}{{\bf(${\mathbb Z}$V\arabic{inner})}}{\usecounter{inner}}
\item (Holomorphicity)  For all $i$ $F^{i} \subset V$  is a holomorphic 
subbundle.
\item (Griffiths transversality)
The filtration $F^{\bullet}$ satisfies
\[
d^{\nabla} : F^{i} \to F^{i-1}\otimes \Omega^{1}_{X}
\] 
for all $i$.
\item (Polarization) Let $V^{p,q} := F^{p}\cap \overline{F}^{q}$. Then the 
natural 
$C^{\infty}$ decomposition $V = \oplus_{p+q = w} V^{p,q}$ is 
$\psi$-orthogonal
and $\psi_{|V^{p,q}}$ is positive definite for $p$-even and negative 
definite
for $p$-odd.
\item (Integrality) The lattice $V_{{\mathbb Z}}$ is horizontal.
\end{list} 
\end{defi}

\begin{rem}\label{rem-opposed}
In general, two exhaustive decreasing filtrations $F^{\bullet}$ and 
$\underline{F}^{\bullet}$ on a complex vector space $V$ are called $w$-opposed
if $V = F^{i} \oplus \underline{F}^{w-i+1}$ for all $i$.
\end{rem}

The non-abelian analogue of the Hodge cohomology classes are the
integral variations of Hodge structures. Thus the nonabelian counterpart
of the group $\vint^{p}$ is the locus
\[
\vint \subset \hdr{1}{X}{\GL{n}{\bbc}},
\]
consisting of all local systems underlying polarizable integral variations of 
Hodge structures.
\glosstex(hodge)[p]{nonabelian-hodge-classes}

It is very hard to describe this locus geometrically. 
The only
information about the properties $\vint$ in general is provided 
by an Arakelov type theorem proven originally by Faltings for variations of 
Hodge structures of weight one \cite{faltings-arakelov} and in the general 
case by Deligne \cite{deligne-arakelov}. According to this theorem 
$\vint$ is a finite set which makes it even harder to characterize
geometrically. The crudest geometric approximation of $\vint$ in
this case is the locus $\vc \subset  \hdr{1}{X}{\GL{n}{\bbc}}$ consisting
of all local systems underlying polarizable complex variations of Hodge 
structures, i.e. variations satisfying all of the above properties with
the exception of the integrality assumption. More precisely we have the
following 
\glosstex(hodge)[p]{nonabelian-pp-classes}

\begin{defi}[\cite{deligne-arakelov},\cite{simpson-higgs}] \label{def-cvhs}  
A {\em complex variation of Hodge 
structures of weight $w$} on a smooth projective $X$ is a complex local system 
$\bbv$ on $X$  together with a flat hermitian form $\psi$ on $\bbv$ so
that the fibers $\bbv_{x}$, $x \in X$ are furnished with a 
decomposition $\bbv_{x} = \oplus_{p \in {\mathbb Z}}\bbv^{p}_{x}$ satisfying 
the following axioms
\begin{list}{{\bf($\bbc$V\arabic{inner})}}{\usecounter{inner}}
\item (Holomorphicity) The subspaces $F^{p} = \oplus_{i \geq p}\bbv_{x}^{p}$ 
and
$\overline{F}^{q} =   \oplus_{i \leq w-q}\bbv_{x}^{q}$ vary with $x$ 
holomorphically and anti-holomorphically respectively.
\item  (Griffiths transversality) If $v$ is a local differentiable section 
of $\bbv$ which is contained
in $F^{p}$ ($\overline{F}^{q}$), then the Lie derivative of $v$ with respect 
to any vector field on $X$ is contained in $F^{p-1}$ ($\overline{F}^{q-1}$)
\item (Polarization) The decomposition $\bbv_{x} = \oplus_{p \in {\mathbb Z}}
\bbv^{p}_{x}$ is $\psi$-orthogonal and $\psi_{|\bbv^{p}_{x}}$ is positive 
definite for $p$ even and negative definite for $p$ odd. 
\end{list}
\end{defi}

The analogue of the Hodge group for the first non-abelian cohomology of 
$X$ would 
be a subgroup $\textbf{Hod}_{\op{na}}(X) \subset \op{Aut}(
\hdol{1}{X}{\GL{n}{\bbc}})$
in the group of algebraic automorphisms of the Dolbeault moduli space with
the property
\[
\phi_{X}(\mycal{V}_{{\mathbb Q}}) = \hdol{1}{X}{\GL{n}{\bbc}}^{
\textbf{Hod}_{\op{na}}(X)},
\]
where $\mycal{V}_{{\mathbb Q}}$ is the locus of ${\mathbb Q}$-variations of 
Hodge structures. As a preliminary step one would like to have good 
analogues of the Weil operators for the non-abelian cohomology. They were 
found by C.~Simpson:

\medskip

\noindent
{\bf Theorem (\cite{simpson-higgs}).} {\it Consider the standard 
$\cstar$
action on the Dolbeault moduli space given by
\[
\xymatrix@R=1pt{
C : & \cstar \ar[r] & {\op{Aut}(\hdol{1}{X}{\GL{n}{\bbc}})} \\
 & t \ar[r]  & ((E,\theta) \mapsto (E,t\theta)).
}
\]
Then the locus of complex variations of Hodge structures coincides with the
fixed-point set of the action $C$, i.e.
\[
\tau_{X}(\vc) = \hdol{1}{X}{\GL{n}{\bbc}}^{\cstar}.
\]}

\medskip

Thus the nonabelian $(p,p)$ classes are just local systems whose Higgs 
bundles are fixed under the $\cstar$-action.

\begin{rem} \label{rem-na-hodge-group} 
The group $\textbf{Hod}_{\op{na}}(X)$ is more elusive. It is clear that
  $\textbf{Hod}_{\op{na}}(X)$ should consist of automorphisms of 
$\hdol{1}{X}{\GL{n}{\bbc}}$ that preserve the non-abelian Hodge decomposition -
that is, automorphism that preserve the rational map to the stack of all
semistable rank $n$ vector bundles: $q: \hdol{1}{X}{\GL{n}{\bbc}} \rightarrow 
H^{1}(X,\GL{n}{{\mathcal O}_{X}})$. More precisely, consider 
the Mumford-type group
\[
{\mathcal G}_{n}(X) := \left\{ (g,s) \left| \; 
\begin{minipage}[c]{3.2in}
where $g \in \op{Aut}(H^{1}(X,\GL{n}{{\mathcal O}_{X}})$ and $s : 
\hdol{1}{X}{\GL{n}{\bbc}} \to \hdol{1}{X}{\GL{n}{\bbc}}$ is an automorphism
that is linear on the fibers of $q$ and induces $g$ on 
$H^{1}(X,\GL{n}{{\mathcal O}_{X}})$
\end{minipage} \;
  \right. \right\}.
\]
On the other hand we have a natural central extension
\[
1\longrightarrow \cstar \longrightarrow {\mathcal G}_{n}^{\op{Weil}}(X)
\longrightarrow \op{Aut}(H^{1}(X,\GL{n}{{\mathcal O}_{X}})) 
\longrightarrow 1,  
\]
consisting of the automorphisms that preserve $q$ and act by Weil operators on
its fibers. 

It is clear that ${\mathcal G}_{n}^{\op{Weil}}(X) \subset {\mathcal 
G}_{n}(X)$ and we 
can define $\textbf{Hod}_{\op{na}}(X)$ to be the minimal subgroup of 
${\mathcal G}_{n}(X)$ that is defined over the $\bbq$ and contains 
${\mathcal G}_{n}^{\op{Weil}}(X)$. This definition is analogous to the 
one we have in the abelian case but there is no  evidence that it is the 
right one.
\end{rem}

\subsection{Twistors and the Hodge filtration} 
\label{ssec-twistors}

From the view-point of the Hodge-de Rham spectral sequence the (abelian)
cohomology groups of a projective manifold come naturally equipped with the
Hodge filtration rather than the Hodge decomposition. This observation becomes
very important when one studies variations of Hodge structures. To understand
the analogue of the Hodge filtration on the first non-abelian cohomology 
one needs a geometric interpretation of the latter due to Deligne 
\cite{deligne-twistor-letter}, Deninger \cite{deninger-g-factors} and Simpson 
\cite{simpson-twistor-letter},
\cite{simpson-kyoto-talk}, \cite{simpson-santa-cruz} which 
we proceed to describe.

\bigskip

\noindent
{\bf A(vi)} Suppose $V$ is a finite dimensional complex vector space 
furnished with an exhaustive decreasing filtration $F^{\bullet}$.
Simpson defines
the {\em Rees module} corresponding to $(V,F^{\bullet})$ to be the
quasi coherent sheaf $\xi(V,F)$ over ${\mathbb A}^{1}$ given by
\[
\xi(V,F) = \sum_{p} \lambda^{-p}F^{p}V \otimes {\mathcal O}_{{\mathbb A}^{1}}
\subset V\otimes {\mathcal O}_{\cstar},
\]
where $\lambda$ is the coordinate on ${\mathbb A}^{1}$. Equivalently $\xi(V,F)$
is the sheafification of the $\bbc[\lambda]$-module $\oplus_{p} 
\lambda^{-p}F^{p}$. The sheaf $\xi(V,F)$ is locally free over ${\mathbb A}^{1}$
and is equipped with a natural $\cstar$-action covering the action on 
${\mathbb A}^{1}$, namely the element $t \in \cstar$ acts by $t^{p}$ on 
the piece $\lambda^{-p}F^{p}V$. Furthermore $\xi(V,F)$ is provided with an
identification between its fiber at $1 \in \cstar$ and $V$. The fiber of 
$\xi(V,F)$ at $0 \in {\mathbb A}^{1}$  is the associated graded
space of $V$, so one way to think of $\xi(V,F)$ is as a canonical deformation 
of $V$ to $\op{Gr}_{F}(V)$.

Conversely
if ${\mathcal V}$ is a quasi-coherent $\cstar$-sheaf on ${\mathbb A}^{1}$, then
${\mathcal V}_{|\cstar} \cong {\mathcal V}_{1}\otimes {\mathcal O}_{\cstar}$ 
and we obtain a filtration $F^{\bullet}$ on ${\mathcal V}_{1}$ by putting
$F^{p} = \{ v \in {\mathcal V}_{1} | v\otimes \lambda^{-p} \in {\mathcal V} 
\}$. If, in addition, ${\mathcal V}$ is locally free, then the natural map
$\xi({\mathcal V}_{1},F) \to {\mathcal V}$ is an isomorphism.

This construction is compatible with $\oplus$, $\otimes$, $\op{Hom}$ and 
passing to duals. A map $f : W \to V$ of filtered vector spaces 
induces a morphism  $\xi(f) : \xi(V,F^{\bullet}V) \to \xi(W,F^{\bullet}W)$ of 
locally free $\cstar$-sheaves which respects kernels, i.e. $\ker(\xi(f)) =
\xi(\ker(f))$. The map $f$ is strictly compatible with the filtrations iff 
$\xi(f)$ is a morphism of vector bundles and in this case $\op{coker}(\xi(f))
= \xi(\op{coker}(f))$. 

If $\underline{F}^{\bullet}$ is another filtration on $V$, then we may apply
the same construction at $\infty$ and glue the resulting sheaves over $\cstar
\subset {\mathbb A}^{1}$  to obtain a locally free $\cstar$ sheaf 
$\xi(V,F,\underline{F})$ with a fiber over $1$ equal to $V$. Suppose 
$F^{\bullet}$ and $\underline{F}^{\bullet}$ are $w$-opposed filtrations (cf. 
remark~\ref{rem-opposed}), then
$\xi(V,F,\underline{F})$ is semistable of slope $w$, i.e. is a direct sum of 
copies of ${\mathcal O}_{\bbp^{1}}(w)$.
If, in addition, $V$ has a real structure given by an anti-holomorphic 
involution $\rho$ and if $\underline{F} = \rho(F)$, then $\xi(V,F,
\underline{F})$ is provided with
an antilinear involution $\sigma$ covering the antipodal involution
$\sigma_{{\mathbb P}^{1}}(\lambda) :=
-\bar{\lambda}^{-1}$ on $\bbp^{1}$ and restricting to $\rho$ on the fiber at 1.
 
To summarize, the Rees module 
construction $\xi$ provides  the following equivalences
of tensor categories \cite[Proposition~1.2 and
Section~2]{simpson-mixed-twistor} (see also
\cite[Apendix]{kaledin-cotangent}):  

\[
\left\{ 
\begin{minipage}[c]{2in}
The category with objects - complex vector spaces with exhaustive 
decreasing filtrations and with morphisms - filtration preserving linear maps
\end{minipage}
\right\} \leftrightarrow  
\left\{
\begin{minipage}[c]{2.3in}
The category of locally free sheaves on ${\mathbb A}^{1}$ equipped with
a $\cstar$ action covering the standard action on ${\mathbb A}^{1}$ 
\end{minipage}
\right\} 
\]
\[
\left\{
\begin{minipage}[c]{2in}
The category with objects - complex vector spaces with exhaustive 
decreasing filtrations and with morphisms - linear maps strictly preserving
the filtrations
\end{minipage}
\right\} \leftrightarrow 
\left\{
\begin{minipage}[c]{2.3in}
The category of $\cstar$-equivariant vector bundles on ${\mathbb A}^{1}$ 
\end{minipage}
\right\}
\]
\[
\left\{
\begin{minipage}[c]{2in}
The category of complex vector spaces equipped with two decreasing
filtrations $F$ and $\overline{F}$ with morphisms - linear maps preserving 
the two filtrations
\end{minipage}
\right\} \leftrightarrow 
\left\{
\begin{minipage}[c]{2.3in}
The category of $\cstar$-equivariant locally free sheaves on $\bbp^{1}$. 
\end{minipage}
\right\}
\]
\[
\left\{
\begin{minipage}[c]{2in}
The category of pure complex Hodge structures
\end{minipage}
\right\} \leftrightarrow 
\left\{
\begin{minipage}[c]{2.3in}
The category of $\cstar$-equivariant vector bundles on $\bbp^{1}$
\end{minipage}
\right\}
\]
Moreover in this case 
pure Hodge structures of weight $w$ correspond to semistable bundles of slope 
$w$. Finally, we have the equivalence
\[
\left\{
\begin{minipage}[c]{2in}
The category of ${\mathbb R}$-Hodge structures 
\end{minipage}
\right\} \leftrightarrow 
\left\{
\begin{minipage}[c]{2.3in}
The category of $\cstar$-equivariant vector bundles on 
$\bbp^{1}$ equipped with
an anti-linear involution covering the antipodal involution on $\bbp^{1}$
\end{minipage}
\right\}
\]
The next crucial observation of Deligne \cite{deligne-twistor-letter} 
is that if $(V_{{\mathbb R}}\subset V, 
F^{\bullet})$ is an ${\mathbb R}$-Hodge structure of weight one, then the
space $V \cong F^{1}V\otimes_{{\mathbb R}} \bbc$ has a canonical action of 
the quaternions ${\mathbb H}$ given by
\begin{align*}
I(f\otimes a) & = f\otimes \sqrt{-1}a \\
J(f\otimes a) & = (-\sqrt{-1}f)\otimes \bar{a} \\
K(f\otimes a) & = IJ(f\otimes a) = - (\sqrt{-1}f)\otimes(\sqrt{-1}\bar{a})
\end{align*}
for a $f\otimes a \in F^{1}V\otimes_{{\mathbb R}} \bbc$.

But for a quaternionic vector space $V$ (or more generally for a 
pseudo-quaternionic\footnote{Recall that a manifold $M$ is called 
pseudo-quaternionic \cite[Chapter~14]{besse} 
if its tangent bundle is equipped with a linear action of
the quaternions ${\mathbb H}\times T_{M} \to T_{M}$} manifold $M$) one 
has the so called twistor 
construction which puts together all the complex structures on $V$ 
(respectively $M$)  coming from the ${\mathbb H}$-action. This construction 
goes as follows. Let ${\mathbb S}$ be the two sphere of pure quaternions of 
norm 1 with the complex structure at a point $q \in {\mathbb S}$ given by the 
left multiplication with $q$. In other words after identifying the tangent 
space $T_{{\mathbb S},q}$ with all pure quaternions $\perp q$ the 
multiplication by $i \in \bbc$ is given by $q\cdot$. Let $M$ be a pseudo
quaternionic manifold. Its twistor space is $Z = M\times {\mathbb S}$ with
the complex structure on the tangent space $T_{M\times {\mathbb S},(m,q)}=
T_{m}\times T_{q}$ given by $(q\cdot,i\cdot)$. If all the almost complex
structures on $M$ induced by the ${\mathbb H}$-action are integrable (such
manifolds are called {\em hypercomplex} \cite{boyer} 
\cite{kaledin-hypercomplex} ), then
$Z$ is a complex manifold and one has: a) $\zeta : Z = M\times {\mathbb S} \to 
{\mathbb S}$ is holomorphic; b) for any $m \in M$ the assignment $q \to 
(m,q)$ is a holomorphic section; c) the map $(m,q) \to (m,-q)$ is an 
antiholomorphic involution.

Identify  ${\mathbb S}$ with $\bbp^{1}$ so that 
$(J,I,-J)$ are identified with $(0,1,\infty)$.
From all we said above it is now clear that for a real Hodge structure 
of weight one the canonical deformation $\xi(V,F)$ of the Hodge filtration 
to the Hodge decomposition is precisely the twistor space of the quaternionic
vector space $V$.

\begin{rem} \label{rem-hyperkahler} If a pseudo-quaternionic manifold $M$
possesses a Riemannian metric $g$ such that the corresponding Hermitian 
two form $g(\bullet,q\bullet)$ is closed for all almost complex structures $q$,
then all the $q$'s are integrable and $g$ is a K\"{a}hler metric in any of 
them. Such manifolds are called \cite{calabi}, \cite[Chapter~14]{besse} 
{\em hyperk\"{a}hler} and can be characterized \cite{hitchinetal}
from the twistor viewpoint as follows. Let $\Omega^{2}_{\zeta} \to Z$ be the 
sheaf of holomorphic two forms along the fibers of $\zeta$. Then a pseudo 
quaternionic manifold $M$ is hyperk\"{a}hler iff its twistor space $Z$ is a 
complex manifold and there exists a section $\Omega \in H^{0}(Z, 
\Omega^{2}_{\zeta}\otimes \zeta^{*}{\mathcal O}_{\bbp^{1}}(2))$ which 
becomes a holomorphic symplectic form when restricted on any fiber of $\zeta$.

To fit this into the above picture notice that a choice of a polarization 
$\psi$ for a ${\mathbb R}$-Hodge structure of weight one corresponds to a 
hyperk\"{a}hler structure on $V$ with the metric being the flat metric given by
$\psi$.
\end{rem}

\bigskip

\noindent
{\bf NA(vi)} The construction $\xi$ of Simpson and its relation to the 
pseudo-quaternionic
and the hyperk\"{a}hler picture give a way of understanding the analogue of the
Hodge filtration for the first nonabelian cohomology  
\cite{deligne-twistor-letter}, \cite{simpson-santa-cruz},
\cite{simpson-kyoto-talk}. 

The starting point is the following beautiful observation of
Hitchin. Denote by \linebreak
$\hdr{1}{X}{\GL{n}{\bbc}}^{\op{reg}}$ and $\hdol{1}{X}{\GL{n}{\bbc}}^{
\op{reg}}$ the smooth loci of the de Rham and Dolbeault cohomology spaces.
Then the homeomorphism $\tau_{X}$ given by the nonabelian Hodge theorem
(see NA(iii)) restricts to  a smooth isomorphism 
\[
\hdr{1}{X}{\GL{n}{\bbc}}^{\op{reg}} = \hdol{1}{X}{\GL{n}{\bbc}}^{\op{reg}} =: 
M^{\op{reg}}.
\]
Furthermore if $I$ and $J$ denote the complex structures coming from 
 
$\hdr{1}{X}{\GL{n}{\bbc}}^{\op{reg}}$ and $\hdol{1}{X}{\GL{n}{\bbc}}^{
\op{reg}}$ respectively, then $K = IJ$ is also a complex structure and this
triple gives rise to a hypercomplex structure on $M^{\op{reg}}$. Finally
the natural polarization on $\hdol{1}{X}{\GL{n}{\bbc}}^{\op{reg}}$ coming 
from the categorical quotient construction is a hyperk\"{a}hler metric.
This was proven by Hitchin \cite{hitchin-self-duality} in the case when
$X$ is a curve and by Deligne \cite{deligne-twistor-letter} and
Fujiki \cite{fujiki} in general.

In view of this and the discussion in A(vi) it is reasonable to try to 
interpret  the Hodge filtration on $\hdr{1}{X}{\GL{n}{\bbc}}$ as a 
special twistor deformation. The main problem is how to deal with the
singularities of the moduli space.

The sheaf $\mycal{D}_{X}$ of differential operators on $X$ is naturally 
filtered: $F^{-p} =$ operators of order $\leq p$. We can build 
$\xi(\mycal{D}_{X},F)$-a quasi coherent sheaf of algebras on $X\times {\mathbb 
A}^{1}$. For a $\lambda \in {\mathbb A}^{1}$ denote by $i_{\lambda}$ the 
inclusion $i_{\lambda} : X \to X\times \{\lambda \} \subset X\times {\mathbb 
A}^{1}$. Define a {\em $\lambda$-connection} on $X$ to be a sheaf $E$ of left
modules for the sheaf of algebras $i_{\lambda}^{*}\xi(\mycal{D}_{X},F)$ that
is locally free as a ${\mathcal O}_{X\times {\mathbb A}^{1}}$-module. Say
that $E$ is {\em semistable} if the Chern classes of $E$ vanish and the 
degree of any subbundle is less than or equal to zero.

It is easy to see that a $\lambda$-connection is a vector bundle $E$ on $X$
together with  a splitting $\nabla_{\lambda}$ of the twisted symbol  sequence
\[
\xymatrix@1{
0 \ar[r] & {\op{End}(E)} \ar[r] & {\mycal{E}(E)}  
\ar[r]^{\lambda\sigma} &  T_{X} \ar[r] 
\ar@/^1pc/[l]^{\nabla_{\lambda}}  &
 0,}
\]
that is a morphism of sheaves of Lie algebras over $\bbc$. Equivalently a 
$\lambda$-connection is an operator $d^{\nabla_{\lambda}} : E \to E\otimes
\Omega^{1}_{X}$ satisfying $d^{\nabla_{\lambda}}(ae) = e\otimes \lambda da + 
ad^{\nabla_{\lambda}}(e)$ and the integrability condition
$d^{\nabla_{\lambda}}\circ d^{\nabla_{\lambda}} = 0$.

Notice that for $\lambda = 1$ a $\lambda$-connection is just a usual 
integrable connection. For $\lambda = 0$ a $\lambda$-connection is a 
Higgs bundle. Therefore this definition provides a deformation from the 
notion of connection to the notion of a Higgs bundle. If $\lambda \neq 0$, 
and $\nabla_{\lambda}$ is a $\nabla$-connection, then $\lambda^{-1}
\nabla_{\lambda}$ is a usual integrable connection. In particular the 
stability of a $\lambda$-connection is automatic for $\lambda \neq 0$
and specializes to the stability for Higgs bundles for $\lambda = 0$.

The importance of the $\lambda$-connections comes from the resulting moduli 
spaces. Simpson had proved \cite[Theorem~4.7]{simpson-moduli1},
\cite[Proposition~4.1]{simpson-santa-cruz} that for any
complex smooth projective variety $X$ there exists a quasi projective 
coarse moduli space  $\hhod{1}{X}{\op{GL}(n,\bbc)} \to {\mathbb
A}^{1}$ of semi-stable 
$\lambda$-connections of rank $n$.
\glosstex(hodge)[p]{nonabelian-hodge-space}
The fibers 
\[
\hhod{1}{X}{\op{GL}(n,\bbc)}_{1} \text{ and }
\hhod{1}{X}{\op{GL}(n,\bbc)}_{0}
\]
over $1, 0 \in {\mathbb A}^{1}$ respectively are 
the moduli spaces
$\hdr{1}{X}{\GL{n}{\bbc}}$ and $\hdol{1}{X}{\GL{n}{\bbc}}$ 
and the natural 
action of $\cstar$ defined by $t : (E,\nabla_{\lambda}) \to 
(E,\nabla_{t\lambda} := t\nabla_{\lambda})$ is an algebraic action on
$\hhod{1}{X}{\op{GL}(n,\bbc)}$ 
covering the standard action on ${\mathbb A}^{1}$. Furthermore
$\hhod{1}{X}{\op{GL}(n,\bbc)}$ 
can be glued to its complex conjugate
\cite[Section~4]{simpson-santa-cruz},  
i.e. there is a
complex analytic space $\hdel{1}{X}{\op{GL}(n,\bbc)}$ over $\bbp^{1}$ 
\glosstex(hodge)[p]{nonabelian-deligne-space}
characterized uniquely by the
following properties: 
\begin{description}
\item[(a)] There is a $\cstar$ action on
$\hdel{1}{X}{\op{GL}(n,\bbc)}$  covering the
standard action on $\bbp^{1}$ and an anti-linear involution $\sigma$ compatible
with the $\cstar$ action and covering the antipodal involution $\lambda \to 
-\bar{\lambda}^{-1}$ on $\bbp^{1}$; 
\item[(b)] There is an algebraic 
$\cstar$-equivariant identification
\[
\hdel{1}{X}{\op{GL}(n,\bbc)}_{|{\mathbb A}^{1}}  \cong
\hhod{1}{X}{\op{GL}(n,\bbc)};
\] 
\item[(c)] On the fiber $\hhod{1}{\GL{n}{\bbc}}_{1} \cong
\hb{1}{X}{\GL{n}{\bbc}}$  
the involution
$\sigma$ takes a representation to the dual of the
complex conjugate representation.
\end{description}

A harmonic bundle yields a family of of holomorphic bundles with 
$\lambda$-connections which extends to a holomorphic section (a twistor line)
$\bbp^{1} \to \hhod{1}{X}{\GL{n}{\bbc}}$. 
Let $M$ be the moduli space of harmonic bundles. 
Deligne had shown \cite{deligne-twistor-letter},
\cite[Section~4]{simpson-santa-cruz}   that the trivialization 
$\hdel{1}{X}{\GL{n}{\bbc}} \cong M \times \bbp^{1}$ 
given by the twistor lines is a homeomorphism 
and that on the set of smooth points it identifies
$\hdel{1}{X}{\GL{n}{\bbc}}^{\op{reg}}$  with
the twistor space for the pseudo-quaternionic structure on $M^{\op{reg}}$.

In view of all this Simpson defines a {\em nonabelian filtration} to be a
$\cstar$-equivariant  scheme over ${\mathbb A}^{1}$. Thus the space
$\hhod{1}{X}{\GL{n}{\bbc}}$ 
is interpreted as the nonabelian Hodge filtration and the space
$\sigma(\hhod{1}{X}{\GL{n}{\bbc}})  :=
\hdel{1}{X}{\GL{n}{\bbc}})_{|\bbp^{1}\setminus \{\infty\}}$ as the
complex conjugate of 
the Hodge filtration. In particular the Dolbeault space 
$\hdol{1}{X}{\GL{n}{\bbc}}$ should be thought of as the associated graded of 
the nonabelian Hodge filtration. The nonabelian counterpart of the usual 
property that a smooth projective $X$ has a Hodge filtration concentrated in 
positive degrees is the statement \cite[Lemma~16]{simpson-kyoto-talk} that for
any point $z \in \hhod{1}{X}{\GL{n}{\bbc}}$ the limit $\lim_{t \to
0}tz$  exists in $\hhod{1}{X}{\GL{n}{\bbc}}$.

\section{The Gauss-Manin connection} \label{sec-GM}

\subsection{Variations of geometric origin} 
\label{ssec-geometric-variations}

\noindent
{\bf A(vii)}
The most important example of ${\mathbb Z}$-variations of Hodge structures 
are the variations of geometric origin which we review next in the simplest 
geometric situation. 

Let $f : X \to S$ be a smooth projective morphism between quasi projective 
varieties.  Put $\hdr{i}{X/S}{\bbc}$ for the sheaf of relative de Rham 
cohomology of degree $i$. Due to our assumptions about $f$ the sheaf 
$\hdr{i}{X/S}{\bbc}$ is locally free and of finite rank, i.e. corresponds to
a vector bundle which we will denote again by $\hdr{i}{X/S}{\bbc}$.

Algebraically the bundle $\hdr{i}{X/S}{\bbc}$ is constructed from the 
relative de Rham complex on $X$. 
To see how this works recall first the following standard notation.

Given an abelian
category $\mycal{A}$, a complex $K^{\bullet}$ of objects in
\glosstex(categories)[p]{abelian}
$\mycal{A}$, and an integer $n$  one defines the {\em stupid} truncation
$\sigma_{\leq n}K^{\bullet}$ of $K^{\bullet}$ 
as the sub complex with terms 
\[
(\sigma_{\geq n}K^{\bullet})^{i} := \left\{ \begin{array}{ll}
0 & \text{ when } i < n \\
K^{i} & \text{ when } i \geq n
\end{array}\right.
\]
Define also $\sigma_{> n}K^{\bullet} := \sigma_{\geq (n-1)}K^{\bullet}$
and $\sigma_{\leq n}K^{\bullet} = K^{\bullet}/\sigma_{< n}K^{\bullet}$.

\glosstex(complexes)[p]{stupid-truncation}

Denote by $\Omega^{\bullet}_{f}$ the full
de Rham complex of relative differential forms on $X/S$ and consider
the truncations $\sigma_{< p}\Omega^{\bullet}_{f}$ and 
$\sigma_{\geq p}\Omega^{\bullet}_{f}$. Explicitly $\sigma_{<
p}\Omega^{\bullet}_{f}$ is the complex 
\[
{\mathcal O}_{X} \to \Omega^{1}_{f} \to \ldots \to \Omega^{p-1}_{f}
\] 
(${\mathcal O}_{X}$ is in degree zero) and  $\sigma_{\geq
p}\Omega^{\bullet}_{f}$ is the complex 
\[
\ldots \to 0 \to \Omega^{p}_{f} \to  \Omega^{p+1}_{f} \to \ldots
\] 
($\Omega^{p}_{f}$ is in 
degree $p$). One has $\hdr{i}{X/S}{\bbc} = 
\bbr^{i}f_{*}\Omega^{\bullet}_{f}$ since the de Rham complex resolves the 
constant sheaf $\bbc$. Furthermore the degeneration of the Hodge-de Rham
spectral sequence gives $F^{p}\hdr{i}{X/S}{\bbc} = \bbr^{i}f_{*}
\sigma_{\geq p}\Omega^{\bullet}_{f}$ and 
$\hdr{i}{X/S}{\bbc}/F^{p}\hdr{i}{X/S}{\bbc} = 
\bbr^{i}f_{*}\sigma_{< p}\Omega^{\bullet}_{f}$, \cite{deligne-lefschetz}, 
\cite{deligne-hodge2}. By the universal coefficients theorem the image of
$R^{i}f_{*}{\mathbb Z}$ in $\hdr{i}{X/S}{\bbc}$ is a full lattice and 
therefore by the covering homotopy property of the cover $R^{i}f_{*}{\mathbb Z}
\to S$ we get a canonical integrable connection on $\hdr{i}{X/S}{\bbc}$ called
the {\em Gauss-Manin connection}. It satisfies the Griffiths transversality
condition and  the first Chern class of any
ample line bundle on $X$ induces a horizontal Hermitian pairing that 
polarizes the variation $\hdr{i}{X/S}{\bbc}$ \cite{griffiths-periods12}.

\medskip

Even though the above topological definition of the Gauss-Manin connection is 
very intuitive, sometimes it is more convenient to have a cohomological
description of the connection. Manin was the first one to realize that this
connection can be defined in purely algebraic terms 
\cite{manin-function-fields}.
Later Grothendieck \cite{grothendieck-crystals} gave a universal algebraic 
construction of the connection
which lead to the notion of a crystal. In general his
construction gives a connection on the de Rham cohomology as an object in the
derived category on $S$ but in our simple case of a smooth projective
$f : X \to S$  it admits the following explicit description \cite{katz-oda,
katz-nilpotent}.

Let $\Omega^{\bullet}_{X}$ be the global de Rham complex on $X$. Let 
$I^{1}$ denote the sub complex which is the image of
$\Omega^{\bullet}_{X}\otimes_{{\mathcal O}_{X}} f^{*}\Omega^{1}_{S}$
and let $I^{2}$ denote the image of $\Omega_{X}^{\bullet - 1}
\otimes_{{\mathcal O}_{X}}  
f^{*}\Omega^{1}_{S}$. Observe that the relative de Rham complex is the quotient
$\Omega^{\bullet}_{f} = \Omega^{\bullet}_{X}/I^{1}$ and that there is an 
isomorphism $I^{1}/I^{2} = \Omega^{\bullet-1}_{f}\otimes_{
{\mathcal O}_{X}} 
f^{*}\Omega^{1}_{S}$. In particular we have an exact sequence of complexes
\[
0 \longrightarrow \Omega^{\bullet-1}_{f}\otimes_{{\mathcal O}_{X}} 
f^{*}\Omega^{1}_{S} \longrightarrow \Omega^{\bullet}_{X}/I^{2} 
\longrightarrow \Omega^{\bullet}_{f} \longrightarrow 0.
\]
The Gauss-Manin connection then is the first edge homomorphism of the 
hyper derived sequence, i.e. the  map
\[
d^{\op{GM}} : \bbr^{i}f_{*}\Omega^{\bullet}_{f} \longrightarrow  
\bbr^{i+1}f_{*}
\Omega^{\bullet}_{f}[-1]\otimes \Omega^{1}_{S} = \bbr^{i}f_{*}
\Omega^{\bullet}_{f}\otimes \Omega^{1}_{S}.
\]
\

\bigskip

\noindent
{\bf NA(vii)} Suppose again that $f : X \to S$ is a smooth projective morphism
of quasi projective varieties and assume for simplicity that the fibers of
$f$ are connected. Denote by $\stackdr{X/S}{n}$ the relative 
moduli stack of rank $n$ local systems (see 
\cite[Theorem~4.7]{simpson-moduli1} for existence). 
\glosstex(moduli)[p]{stackdr}
\glosstex(moduli)[p]{stackdol}
Simpson had shown 
\cite{simpson-moduli2} that $\stackdr{X/S}{n}$ has a natural 
structure of a crystal of stacks over $S$ dubbed by him the 
{\em nonabelian Gauss-Manin connection}. The construction mimics 
Grothendieck's definition for the abelian case in the moduli context. What 
makes the construction work 
is the observation that if $S'$ is an $S$-scheme which contains
a closed subscheme $S_{0}'$ defined by a nilpotent ideal then a crystal on
$X'/S'$ is canonically  equivalent to a crystal on $X'_{0}/S'$. Here $X' =
X\times_{S} S'$ and $X'_{0} = X\times_{S} S'_{0}$. Because of that the 
functor $M^{\natural}_{\op{crys}}(X/S,n)$ that assigns to a pair
$S'_{0}\subset S'$ the set of isomorphism classes of crystals of rank $n$
vector bundles on $X_{0}'/S'$ is a crystal of functors. Moreover Simpson
proves \cite[Lemmas~8.1 and 8.2]{simpson-moduli2} that 
$M^{\natural}_{\op{crys}}(
X/S,n)$ is isomorphic to the relative de Rham moduli functor 
$M^{\natural}_{\op{DR}}(X/S,n)$. Since the latter is 
represented by $\stackdr{X/S}{n}$ one gets a 
structure of a $S$-crystal on $\stackdr{X/S}{n}$, i.e. 
an isomorphism
\[
\varphi^{\op{GM}} : p_{1}^{*}\stackdr{X/S}{n} 
\widetilde{\longrightarrow} p_{2}^{*}\stackdr{X/S}{n} 
\]
on $(S\times S)^{\wedge}$ satisfying the usual cocycle condition. The same 
construction works also for the relative moduli space 
$\hdr{1}{X/S}{\op{GL}_{n}(\bbc)}$ since the moduli functor 
$M^{\natural}_{\op{DR}}
(X/S,n)$ is universally coarsely represented by the scheme 
$\hdr{1}{X/S}{\op{GL}_{n}(\bbc)}$. It turns out that the regular points of 
$\stackdr{X/S}{n}$ (respectively $\hdr{1}{X/S}{\op{GL}_{n}(\bbc)}$) form
a local system of stacks (respectively schemes) on $S$.
In the same way we may define the moduli stack $\stackdol{X/S}{n}$ of
relative Higgs bundles of rank $n$. 

It turns out that 
the stacks $\stackdr{X/S}{n}$ and $\stackdol{X/S}{n}$ are better
approximations of the first non-abelian de Rham and Dolbeault
cohomology than the moduli schemes $\hdr{1}{X/S}{\op{GL}_{n}(\bbc)}$ and 
$\hdol{1}{X/S}{\op{GL}_{n}(\bbc)}$. This can be seen in many ways and
is especially transparent when one considers non-abelian cohomology of
degree greater than one \cite{simpson-algebraic-nah}.
Another indication is the fact that the set of $\bbc$-valued points of
say $\stackdr{X/S}{n}$ is exactly the set of equivalence classes of
relative local systems on $X/S$ whereas the set of $\bbc$-points of
$\hdr{1}{X/S}{\op{GL}_{n}(\bbc)}$ is the set of equivalence classes of
semisimplifications of rank $n$ relative local systems. 

From now on we will freely work with these stacky non-abelian cohomologies.

\medskip

Our first task  describe the corresponding nonabelian
connection explicitly at 
smooth points of the moduli. The natural framework for such a description
is that of algebraic stacks and their tangent stacks. Our main references 
are \cite{laumon-stacks} and \cite[Appendix]{vistoli-stacks}. For the
convenience of the reader we have reproduced the necessary statements
in Appendix \ref{app-tangent-stacks}.

First we will need an intrinsic modular description of the tangent stacks
(see  Section \ref{appss-truncated-cotangent} for a definition) of
the moduli stacks $\stackdr{X/S}{n}$ and 
$\stackdol{X/S}{n}$. Denote by 
\begin{align*}
\pdr : \stackdr{X/S}{n} & \longrightarrow S \\
\pdol :\stackdol{X/S}{n} & \longrightarrow S  
\end{align*}
the structure morphisms of the relative de Rham and Dolbeault stacks 
respectively. Let $(T;F,\nabla) \in \op{Ob}(\stackdr{X/S}{n})$ and 
$(T;E,\theta) \in \op{Ob}(\stackdol{X/S}{n})$. In other words $T \to S$ is
an $S$-scheme and if $X_{T} := X\times_{S} T$, then $(F,\nabla)$ and 
$(E,\theta)$ are a relative local system and a relative Higgs bundle on 
$X_{T}/T$. As usual we will describe the vertical and the total tangent 
stacks of the de Rham and Dolbeault stacks in terms of suitable
deformation-obstruction complexes. 

\begin{lemma-defi} \label{lemma-def-vertical-tangent} In the above notation
put $f_{T} : X_{T} \to T$ for the natural projection. Then  
there are well defined complexes of sheaves on $X_{T}$
\begin{description}
\item[{\em De Rham version}]
$\lieg{f_{T}}{F,\nabla} := (\op{End}(F)\otimes \Omega^{\bullet}_{f_{T}}, 
\ad{\nabla})$, where 
\[
\ad{\nabla}(m) = [d^{\nabla},m] = d^{\nabla}\circ m - 
(m\otimes \op{id})\circ d^{\nabla}.
\]
\item[{\em Dolbeault version}] $\lieg{f_{T}}{E,\theta} := (\op{End}(E)\otimes 
\Omega^{\bullet}_{f_{T}}, \ad{\theta})$,
where 
\[
\ad{\theta}(m) = [\theta,m] = \theta\circ m - 
(m\otimes \op{id})\circ \theta.
\]
\end{description}
\end{lemma-defi}
{\bf Proof.} The only thing that needs checking is the 
${\mathcal O}_{X_{T}}$-linearity of $\ad{\nabla}(m) : F \to F\otimes 
\Omega^{1}_{f_{T}}$ and $\ad{\theta}(m) : E \to E\otimes \Omega^{1}_{f_{T}}$ 
respectively.

Since both $m$ and $\theta$ are ${\mathcal O}_{X_{T}}$-linear it is clear that
 $\ad{\theta}(m)$ will also be ${\mathcal O}_{X_{T}}$-linear. To check the 
de Rham case we need to show that for any $m \in \op{End}(F)$ and any $a\in F$,
$\phi \in {\mathcal O}_{X_{T}}$ we have $\ad{\nabla}(m)(\phi a) = 
\phi\ad{\nabla}(m)(a)$.
Put $d_{f_{T}} : \Omega^{i}_{f_{T}} \to \Omega^{i+1}_{f_{T}}$ for the 
exterior differentiation along the fibers of $f_{T}$. By definition we have 
\[
\begin{split}
\ad{\nabla}(m)(\phi a) & = [d^\nabla, m](\phi a) = d^\nabla(m(\phi a)) - 
(m\otimes\op{id})(d^\nabla(\phi a)) \\
&  = d^\nabla(\phi\cdot m(a)) - (m\otimes\op{id})(a\otimes d_{f_{T}}\phi 
 + \phi\cdot d^\nabla(a)) \\
& = m(a)\otimes d_{f_{T}}\phi + \phi d^\nabla\circ m(a) - m(a)\otimes 
d_{f_{T}}\phi 
 + \phi(m\otimes\op{id})\circ d^\nabla(a) \\
& = \phi\ad{\nabla}(m)(a),
\end{split}
\]
which proves the lemma since the integrability of $\nabla$ and the 
symmetry condition $\theta\wedge \theta = 0$ on $\theta$ guarantee that 
$\ad{\nabla}$ and $\ad{\theta}$ will be differentials. \hfill $\Box$

\bigskip

\noindent
The complexes of the previous lemma carry information 
about the vertical tangent stacks of the structure morphisms $\pdr$ and 
$\pdol$. Before we explain that, recall the following construction. Given any 
algebraic $S$-stack $\mycal{X}$ and any complex of sheaves of abelian groups 
$E^{0} \to E^{1}$ on $\mycal{X}$ one may consider the stack theoretic 
quotient of the translation action of $E^{0}$ on $E^{1}$. In this way we get
an $S$-stack $\hh(E^{\bullet}) := [E^{1}/E^{0}]$ having also a structure 
of a strictly commutative Picard stack over $\mycal{X}$ (see \cite[Section 
1.4 of Expos\'{e} XVIII]{sga4} and \cite[Section 2]{behrend-fantechi} for
details). For an object $V \in 
\op{Ob}(\mycal{X})$ the groupoid of sections of $\hh(E^{\bullet})$
over $V$ is the category of pairs $(R,r)$, where $R$ is an $E^{0}$-torsor on 
$V$ and $r : R \to E^{1}_{|V}$ is an $E^{0}$-equivariant morphism of sheaves 
on $V$.

Denote by $\stackdrreg{X/S}{n}$ and $\stackdolreg{X/S}{n}$ the parts of 
of the stacks $\stackdr{X/S}{n}$ and $\stackdol{X/S}{n}$ over which
the morphisms $\pdr$ and $\pdol$ are smooth. Put $(F_{\op{un}},
\nabla_{\op{un}}) \to
\stackdr{X/S}{n}\times_{S} X$ and $(E_{\op{un}},\theta_{\op{un}}) \to 
\stackdol{X/S}{n}\times_{S} X$ for the  universal families and let 
$\fdr : \stackdr{X/S}{n}\times_{S} X \to \stackdr{X/S}{n}$ and 
$\fdol : \stackdol{X/S}{n}\times_{S} X  \to \stackdol{X/S}{n}$ denote the 
natural projections. We have the following lemma.

\glosstex(moduli)[p]{stackdrreg}
\glosstex(moduli)[p]{stackdolreg}

\begin{lem} \label{lem-vertical-tangent-spaces} \
 
\begin{list}{{\em (\alph{inner})}}{\usecounter{inner}}

\item The complexes $\lieg{\fdr}{F_{\op{un}},\nabla_{\op{un}}}$ and 
$\lieg{\fdol}{E_{\op{un}},\theta_{\op{un}}}$  are the 
deformation-obstruction complexes
for the moduli functors    $M^{\natural}_{\op{DR}}(X/S,n)$ and 
$M^{\natural}_{\op{Dol}}(X/S,n)$ over $S$. That is, the infinitesimal
automorphisms, deformations and obstructions of the pairs $(F,\nabla)$ and
$(E,\theta)$ on the fixed $X_{T}$ are
parameterized respectively by the vector spaces

\begin{center}
\begin{tabular}{|c||c|c|} \hline
 \begin{minipage}[c]{1.35in} \addtolength{\baselineskip}{-4pt}
 {\footnotesize\em  geometric objects }
\end{minipage} &
 \begin{minipage}[c]{1.35in}  \addtolength{\baselineskip}{-4pt}
{\footnotesize\em
 \medskip \noindent $(F,\nabla)$ 
with  $T\to S$  \\ kept fixed
\medskip \ } 
 \end{minipage} &
  \begin{minipage}[c]{1.35in}  \addtolength{\baselineskip}{-4pt}
 {\footnotesize\em \medskip 
\noindent $(E,\theta)$ with  $T\to S$ \\ kept fixed 
\medskip \ }
\end{minipage} \\
\hline
\hline
 \begin{minipage}[c]{1.35in} \addtolength{\baselineskip}{-4pt}
{\footnotesize\em \medskip 
\noindent infinitesimal \\  automorphisms 
\medskip \ } \end{minipage} &
  $\bbh^{0}(X_{T},\lieg{f_{T}}{F,\nabla})$ &
  $\bbh^{0}(X_{T},\lieg{f_{T}}{E,\theta})$ \\
\hline
\begin{minipage}[c]{1.35in} \addtolength{\baselineskip}{-4pt}
{\footnotesize\em \medskip \noindent infinitesimal \\
 deformations 
\medskip \ }\end{minipage} &
 $\bbh^{1}(X_{T},\lieg{f_{T}}{F,\nabla})$ &
  $\bbh^{1}(X_{T},\lieg{f_{T}}{E,\theta})$  \\
\hline
\begin{minipage}[c]{1.35in} \addtolength{\baselineskip}{-4pt}
{\footnotesize\em \medskip \noindent infinitesimal \\
obstructions \medskip \ } 
\end{minipage} &
 $\bbh^{2}(X_{T},\lieg{f_{T}}{F,\nabla})$ & 
$\bbh^{2}(X_{T}, \lieg{f_{T}}{E,\theta})$
\\
\hline
\end{tabular} 
\end{center}

\item There are isomorphisms of Picard stacks over $\stackdrreg{X/S}{n}$
and $\stackdolreg{X/S}{n}$ respectively:
\begin{align*}
T_{\pdr} & = \hh(\bbr^{0}f_{\op{DR}*}\lieg{\fdr}{F_{\op{un}},
\nabla_{\op{un}}}
\stackrel{0}{\to} \bbr^{1}f_{\op{DR}*}\lieg{\fdr}{F_{\op{un}},
\nabla_{\op{un}}}) \\
T_{\pdol} & = \hh(\bbr^{0}f_{\op{Dol}*}\lieg{\fdol}{E_{\op{un}},
\theta_{\op{un}}}
\stackrel{0}{\to} \bbr^{1}f_{\op{Dol}*}\lieg{\fdol}{E_{\op{un}},
\theta_{\op{un}}}) 
\end{align*}
\end{list} 
\end{lem}
{\bf Proof.} (a) The proof is an easy but lengthy cocycle computation. We will
work it out for the infinitesimal deformations of de Rham case only. 
The other cases are completely
analogous. Put $B$ for the spectrum of the dual numbers, i.e. $B$ is the
local nilpotent scheme $\op{Spec}(\bbc[\varepsilon]/(\varepsilon^{2}))$. 
Suppose we are given an infinitesimal deformation $(\widetilde{F},
\widetilde{\nabla})$ of $(F,\nabla)$. Concretely $\widetilde{F} \to X_{T}
\times B$ is a
locally free sheaf of rank $n$ that specializes to $F$ over the closed point 
of $B$ and $\widetilde{\nabla }$ is a relative integrable connection on 
$\widetilde{F}$ 
which specializes to $\nabla$.  Choose an acyclic \v{C}hech cover 
${\mathfrak U}$  of $X_{T}$ by Zariski open sets which trivializes the 
bundle $F$
and denote by $h_{U} : F_{|U} \to {\mathcal O}_{U}^{\oplus n}$ a  
trivialization of $F$ over $U \in {\mathfrak U}$. In terms of the cover 
${\mathfrak U}$ the bundle $F$ is described by the cocycle $\{ g_{UV} \}
\in Z^{1}({\mathfrak U}, \GL{n}{{\mathcal O}_{X_{T}}})$ where $g_{UV} = 
h_{U}\circ h_{V}^{-1}$.  Similarly the connection $\nabla$ is given by the 
cochain $\{ a_{U} \} \in C^{0}({\mathfrak U}, \op{End}({\mathcal O}_{X_{T}})
\otimes
\Omega^{1}_{f_{T}})$ where $a_{U} = (h_{U}\otimes \op{id}_{\Omega^1})\circ 
d^{\nabla}\circ h_{U}^{-1} - d$. Since $X_{T}\times B$ is a nilpotent scheme 
whose reduced support is $X_{T}$ we can describe the pair $(\widetilde{F},
\widetilde{\nabla})$ in terms of the covering ${\mathfrak U}$ by a pair of
cochains $\{ \tilde{g}_{UV} \} \in C^{1}({\mathfrak U}, \GL{n}{{\mathcal 
O}_{X_{T}}[\varepsilon]/(\varepsilon^{2})})$ and $\{ \tilde{a}_{U} \} \in 
C^{0}({\mathfrak U}, \op{End}({\mathcal O}_{X_{T}}^{\oplus n}[\varepsilon]/
(\varepsilon^{2})))$ satisfying
\begin{align} 
\tilde{g}_{UV} = &  \tilde{g}_{UV}\tilde{g}_{VW} \label{eq-cocycle}\\
\tilde{a}_{U} = & \tilde{g}_{UV}d\tilde{g}_{VU} + \tilde{g}_{UV}\tilde{a}_{V}
\tilde{g}_{VU}. \label{eq-connection}
\end{align}
Also since $(\widetilde{F}, \widetilde{\nabla})$ specializes to $(F,\nabla)$
we have 
\begin{align*}
\tilde{g}_{UV} = & g_{UV} + \varepsilon x_{UV}, \\
\tilde{a}_{U} = & a_{U} + \varepsilon m_{U},
\end{align*}
with $\{ x_{UV} \} \in C^{1}({\mathfrak U},\op{End}({\mathcal 
O}_{X_{T}}^{\oplus n}))$ and
$\{ m_{U} \} \in C^{0}({\mathfrak U}, \op{End}({\mathcal O}_{X_{T}}^{\oplus n})
\otimes 
\Omega^{1}_{f_{T}})$.

The integrability condition on the relative connection $\widetilde{\nabla}$
is 
\[
d\tilde{a}_{U} + \frac{1}{2}[\tilde{a}_{U},\tilde{a}_{U}] = 0.
\] 
with $d = d_{f_{T}}$ being the exterior differentiation along the fibers of 
$f_{T}$. Due to the integrability of $\nabla$ this is equivalent to
\begin{equation}
dm_{U} + [a_{U},m_{U}] = 0 \label{eq-integrable}
\end{equation}
By comparing the coefficients in front of $\varepsilon$ in the cocycle condition
(\ref{eq-cocycle}) we obtain the identity $x_{UW} = g_{UV}x_{VW} + 
x_{UV}g_{VW}$ which can be rewritten as the identity $\thick{x}_{UW} = 
\thick{x}_{UV} + \thick{x}_{VW}$ for $\thick{x}_{UV} := h_{U}^{-1}\circ x_{UV}
\circ h_{V} 
\in \Gamma(U\cap V, \op{End}(F))$. Thus $\{ \thick{x}_{UV} \} \in Z^{1}(
{\mathfrak
U}, \op{End}(F))$. Similarly the connection condition (\ref{eq-connection}) 
gives
\[
m_{U} = g_{UV}dx_{VU} + x_{UV}dg_{VU} + x_{UV}a_{V}g_{VU} + g_{UV}m_{V}g_{VU}
+ g_{UV}a_{V}x_{VU}.
\] 
The next step is to rewrite $dx_{VU}$ and $dg_{VU}$ via  Cartan homotopy 
formula: 
\begin{align*}
m_{U} = & g_{UV}\circ d\circ x_{VU} + g_{UV}\circ x_{VU}\circ d  + x_{UV}\circ 
d \circ g_{VU} + x_{UV}\circ g_{VU}\circ d \\
+ & x_{UV}a_{V}g_{VU} + g_{UV}m_{V}g_{VU} + g_{UV}a_{V}x_{VU}.
\end{align*}
Taking into account that $g_{UV} = h_{U}\circ h_{V}^{-1}$, $\thick{x}_{UV}
= - \thick{x}_{VU}$ and $(d^\nabla)_{|U\cap V} = (h_{U}\otimes\op{id})
\circ (d + a_{U}) \circ h_{U}^{-1} = (h_{V}\otimes\op{id})\circ (d + 
a_{V}) \circ h_{V}^{-1}$ the last identity can be rewritten as 
\[
\thick{m}_{U} - \thick{m}_{V} = d^\nabla \circ \thick{x}_{UV} - (\thick{x}_{UV}
\otimes 
\op{id}_{\Omega^{1}}) \circ d^\nabla 
\]   
where $\thick{m}_{U} = (h_{U}\otimes\op{id})\circ m_{U} \circ h_{U}^{-1}$. 
In combination with (\ref{eq-integrable}) this implies that the pair 
$(\{ \thick{x}_{UV} \}, \{ \thick{m}_{U} \}) \in 
C^{1}({\mathfrak U}, \lieg{f_{T}}{F,\nabla})$ is a cocycle and by 
construction 
completely determines the infinitesimal deformation $(\widetilde{F}, 
\widetilde{\nabla})$. It is straightforward to check that cohomologous cocycles
correspond to isomorphic deformations and so $(\widetilde{F}, 
\widetilde{\nabla})$ determines and is determined by an element in 
$\bbh^{1}(X_{T}, \lieg{f_{T}}{F,\nabla})$.

\medskip

\noindent
(b) A cheap way to prove this part is to combine
Simpson's formality result \cite[Lemma~3.5]{simpson-higgs} 
\cite[Proposition~10.5]{simpson-moduli2} with the
isomorphism (\ref{tangent-stack-isomorphism}).

One can also give a direct argument which we proceed to explain in
the de Rham case. 
Fix an object $(T;F,\nabla) 
\in \op{Ob}(\stackdr{X/S}{n})$ and let $\alpha_{T} : T \to \stackdr{X/S}{n}$ be
the corresponding morphism of $S$-stacks.  The fiber stack $T_{\pdr,
(T;F,\nabla)} :=
T_{\pdr}\times_{\alpha_{T}} \stackdr{X/S}{n}$ is naturally a $T$-stack.
Explicitly for a given $T$-scheme $U \to T$ let $X_{U} := X_{T}\times_{T} U$ 
and let $f_{U} : X_{U} \to U$ denote the natural projection. We abuse notation
and write $(F,\nabla)$ for the pull-back of $(F,\nabla)$ to $X_{U}$. The  
groupoid of sections $T_{\pdr, (T;F,\nabla)}$ over $U$ is the category of 
$S$-relative 
connections $(\widetilde{F}, \widetilde{\nabla}) \to X_{U}\times B$ such that
if $j : X_{U} \hookrightarrow X_{U}\times B$ is the inclusion coming from the
inclusion of the closed point in $B$, then there exists an isomorphism
$j^{*}(\widetilde{F}, \widetilde{\nabla}) \cong (F_{U},\nabla_{U})$. The 
choice of such an isomorphism however is not considered part of the data 
$(\widetilde{F}, \widetilde{\nabla})$. Similarly, if we denote by 
$\hh_{\op{DR}, (T;F,\nabla)} $ the fiber-product 
\[
\hh(\bbr^{0}f_{\op{DR}*}\lieg{\fdr}{F_{\op{un}}, 
\nabla_{\op{un}}}\stackrel{0}{\to} \bbr^{1}f_{\op{DR}*}\lieg{\fdr}{F_{\op{un}},
\nabla_{\op{un}}})\times_{\alpha_{T}} \stackdr{X/S}{n}
\]
we can identify the groupoid of sections of $\hh_{\op{DR}, (T;F,\nabla)} $ 
over $U$ as the 
category of all pairs $(R,r)$ where $R$ is an $\bbh^{0}(X_{U}, 
\lieg{f_{U}}{F,\nabla})$ torsor and 
$r \in \bbh^{1}(X_{U}, \lieg{f_{U}}{F,\nabla})$. Next we have a functor 
\[
a : T_{\pdr, (T;F,\nabla)} \longrightarrow \hh_{\op{DR}, (T;F,\nabla)}
\]
given by  $(\widetilde{F}, \widetilde{\nabla})  \mapsto (R_{(\widetilde{F}, 
\widetilde{\nabla})}, r_{(\widetilde{F}, \widetilde{\nabla})})$ where 
$R_{(\widetilde{F}, \widetilde{\nabla})}$ is the set of all surjective
morphisms of crystals 
$(\widetilde{F}, \widetilde{\nabla}) \to (F,\nabla)$ on $X_{T}\times B$ and
$r_{(\widetilde{F}, \widetilde{\nabla})}$ is the (Kodaira-Spencer)
class of $(\widetilde{F}, \widetilde{\nabla})$ constructed in (a). 
Equivalently $R_{(\widetilde{F}, \widetilde{\nabla})}$ can be thought of as 
the set of all realizations of the $\Omega^{\bullet}_{f_{T}}$-dg-module 
$(\widetilde{F}\otimes \Omega^{\bullet}_{f_{T}}, d^{\widetilde{\nabla}})$ 
as an extension of $F$ by the $\Omega^{\bullet}_{f_{T}}$-dg-module 
$(F\otimes \Omega^{\bullet}_{f_{T}}, d^{\nabla})$ and $r_{(\widetilde{F}, 
\widetilde{\nabla})}$ as the extension class of any such realization.

\hfill $\Box$

\bigskip

To describe the full tangent spaces to the relative de Rham and Dolbeault 
stacks we will need to extend the complexes $\lieg{f_{T}}{F,\nabla}$ and 
$\lieg{f_{T}}{E,\theta}$ in a suitable fashion. For a coherent sheaf
$F \to X_{T}$ denote by $\mycal{E}_{f_{T}}(F)$ the $f_{T}$-relative Atiyah
algebra of $F$. Then we have

\begin{lemma-defi} \label{lem-def-full-tangent}
For any $F$ algebraic vector bundle on $X_{T}$ one has:
\begin{list}{{\em (\roman{inner})}}{\usecounter{inner}}
\item The natural morphism 
\[
\xymatrix@R=1pt{
{\ell} : & {\mycal{E}_{f_{T}}(F)} \ar[r] & {\mycal{E}_{f_{T}}(F\otimes
\Omega^{1}_{f_{T}})} \\
 & \del  \ar[r] & {\del\otimes \op{id} + \op{id}\otimes L_{\sigma(
\del)}}}
\]
is well defined and ${\mathcal O}_{X_{T}}$-linear.
\item There are well defined complexes on $Y$
\begin{description}
\item[{\em De Rham version}] 
\[
\liea{f_{T}}{F,\nabla} := \mycal{E}_{f_{T}}(F)  \stackrel{\ad{\nabla}}{\to}
\op{End}(F)\otimes \Omega^{1}_{f_{T}} \stackrel{\ad{\nabla}}{\to}
\op{End}(F)\otimes \Omega^{2}_{f_{T}} \to \ldots
\]
where $\ad{\nabla}(\del) := d^{\nabla}\circ \del - \ell(\del) 
\circ d^{\nabla}$ for $\del \in \mycal{E}_{f_{T}}(F)$ and is as in 
lemma-definition~\ref{lemma-def-vertical-tangent} for $m \in \op{End}(F)
\otimes \Omega^{i}_{Y}$.
\item[{\em Dolbeault version}]
\[
\liea{f_{T}}{E,\theta} := \mycal{E}_{f_{T}}(E)  \stackrel{\ad{\theta}}{\to}
\op{End}(E)\otimes \Omega^{1}_{f_{T}}\stackrel{\ad{\theta}}{\to}
\op{End}(E)\otimes \Omega^{2}_{f_{T}} \to \ldots
\]
where $\ad{\theta}(\del) := \theta\circ \del - \ell(\del)\circ
\theta$ for $\del \in \mycal{E}_{f_{T}}(E)$ and is as in 
lemma-definition~\ref{lemma-def-vertical-tangent} for $m \in \op{End}(E)
\otimes \Omega^{i}_{f_{T}}$.
\end{description} 
\end{list}
\end{lemma-defi}
{\bf Proof.} To prove (i) we first need  to show that for any local section
$\del \in \mycal{E}_{f_{T}}(F)$ the $\bbc$-linear map $\ell(\del) : F\otimes 
\Omega^{1}_{f_{T}} \to F\otimes \Omega^{1}_{f_{T}}$ belongs to the  
Atiyah algebra
of $F$. For any local sections $\varphi \in {\mathcal O}_{X_{T}}$, $a \in F$, 
$\alpha \in \Omega^{1}_{f_{T}}$ we have by definition
\begin{align*}
\ell(\del)((\varphi a)\otimes \alpha) & = \del(\varphi a)\otimes \alpha + 
(\varphi a)\otimes L_{\sigma(\del)}\alpha = \\
& = (L_{\sigma(\del)}\varphi)a\otimes \alpha + \varphi\del a\otimes\alpha +
\varphi a\otimes L_{\sigma(\del)}\alpha = \\
& = \varphi\del a\otimes \alpha + a\left( \varphi L_{\sigma(\del)}\alpha 
+ (L_{\sigma(\del)}\varphi))\alpha \right) = \\
& = \del a\otimes (\varphi \alpha) + a \otimes L_{\sigma(\del)}(\varphi 
\alpha) = \ell(\del)(a\otimes \varphi\alpha).
\end{align*}
In particular $\ell(\del)$ is a well defined $\bbc$-linear endomorphism of
$F\otimes \Omega^{1}_{f_{T}}$ since it respects the 
${\mathcal O}_{X_{T}}$-linearity
of the tensor product of $F$ and $\Omega^{1}_{f_{T}}$. Moreover 
the second line 
in the above identity shows that for any $A \in F\otimes \Omega^{1}_{f_{T}}$
we have $\ell(\del)(\varphi A) = \varphi \ell(\del)(A) + (L_{\sigma(\del)}
\varphi)A$ and thus $\ell(\del)$ is actually an element in $\mycal{E}(
 F\otimes \Omega^{1}_{f_{T}})$ with $\sigma(\ell(\del)) =
\sigma(\del)$. 
Since the 
${\mathcal O}_{X_{T}}$-linearity of $\ell$ is clear from its definition this 
completes the proof of (i). 

\medskip

For part (ii) we will only check that $\ad{\nabla}(\del)$ is 
${\mathcal O}_{X_{T}}$-linear for every $\del \in
\mycal{E}_{f_{T}}(F)$ 
since the 
corresponding statement for $\ad{\theta}$ is checked in exactly the same way.
Let $a$ and $\varphi$ be as before. Then
\begin{align*}
\ad{\nabla}(\del)(\varphi a) & = (d^{\nabla}\circ \del - \ell(\del)\circ 
d^{\nabla})(\varphi a) = \\
& = d^{\nabla}\circ ((L_{\sigma(\del)}\varphi)a + \varphi(\del a )) - 
\ell(\del) (a\otimes d\varphi + \varphi d^{\nabla}a) = \\
& = a\otimes (d\circ L_{\sigma(\del)})\varphi + \varphi (d^{\nabla}\circ 
\del)a \\
& - a\otimes (L_{\sigma(\del)}\circ d)\varphi - \varphi (\ell(\del)\circ 
d^{\nabla})a.
\end{align*}
Also for a vertical vector field $v \in T_{f_{T}}$ the 
Cartan homotopy formula gives 
\[
L_{v} = d\circ i_{v} + i_{v}\circ d  \;\; \text{on} \;\; 
\Omega^{\bullet}_{f_{T}},
\]
where $i_{v}$ is the contraction with $v$.

Therefore $d\circ L_{v} = d\circ i_{v}\circ d = L_{v}\circ d$ and hence
\[
\ad{\nabla}(\del)(\varphi a ) = \varphi \ad{\nabla}(\del)(a)
\]
which completes the proof of the lemma.
\hfill $\Box$

\bigskip

\noindent
Again the complexes $\mycal{E}_{f_{T}}(F,\nabla)$ and 
$\mycal{E}_{f_{T}}(E,\theta)$
can be interpreted as deformation-ob\-struc\-tion complexes but before we
state the result we need to recall some standard notation. 

Given an abelian
category $\mycal{A}$, a complex $K^{\bullet}$ of objects in
$\mycal{A}$, and an integer $n$  one defines the {\em canonical} truncation
$\tau_{\leq n}K^{\bullet}$ of $K^{\bullet}$ 
as the sub complex with terms 
\[
(\tau_{\leq n}K^{\bullet})^{i} := \left\{ \begin{array}{ll}
K^{i} & \text{ when } i < n \\
\ker(d) & \text{ when } i = n \\
0 & \text{ when } i > n
\end{array}\right.
\]
Define also $\tau_{< n}K^{\bullet} := \tau_{\leq (n-1)}K^{\bullet}$
and $\tau_{\geq n}K^{\bullet} = K^{\bullet}/\tau_{< n}K^{\bullet}$. 
By definition we have $H^{i}\tau_{\leq n}K^{\bullet} =
H^{i}K^{\bullet}$ when $i \leq n$ and $H^{i}\tau_{\leq n}K^{\bullet} =
0$ for $i > n$. Similarly $H^{i}\tau_{\geq n}K^{\bullet} =
H^{i}K^{\bullet}$ when  $i \geq n$ and $H^{i}\tau_{\geq n}K^{\bullet} =
0$ for $i < n$.

\glosstex(complexes)[p]{canonical-truncation}

We are now ready to state 

\begin{lem} \label{lem-full-tangent-spaces} \

\begin{list}{{\em (\alph{inner})}}{\usecounter{inner}}

\item The complexes $\liea{\fdr}{F_{\op{un}},\nabla_{\op{un}}}$ and 
$\liea{\fdol}{E_{\op{un}},\theta_{\op{un}}}$  are the 
deformation-obstruction complexes
for the moduli functors    $M^{\natural}_{\op{DR}}(X/S,n)$ and 
$M^{\natural}_{\op{Dol}}(X/S,n)$. That is, the infinitesimal
automorphisms, deformations and obstructions of the triples $(X_{T},
F,\nabla)$ and
$(X_{T},E,\theta)$  over $S$ are
parameterized respectively by the vector spaces

\begin{center}
\begin{tabular}{|c||c|c|} \hline
 \begin{minipage}[c]{1.35in} \addtolength{\baselineskip}{-4pt}
 {\footnotesize\em  geometric objects }
\end{minipage} &
 \begin{minipage}[c]{1.35in}  \addtolength{\baselineskip}{-4pt}
{\footnotesize\em
 \medskip \noindent $(X_{T},F,\nabla)$ over $S$
\medskip \ } 
 \end{minipage} &
  \begin{minipage}[c]{1.35in}  \addtolength{\baselineskip}{-4pt}
 {\footnotesize\em \medskip 
\noindent $(X_{T},E,\theta)$ over $S$ 
\medskip \ }
\end{minipage} \\
\hline
\hline
 \begin{minipage}[c]{1.35in} \addtolength{\baselineskip}{-4pt}
{\footnotesize\em \medskip 
\noindent infinitesimal \\  automorphisms 
\medskip \ } \end{minipage} &
  $\bbh^{0}(X_{T},\liea{f_{T}}{F,\nabla})$ &
  $\bbh^{0}(X_{T},\liea{f_{T}}{E,\theta})$ \\
\hline
\begin{minipage}[c]{1.35in} \addtolength{\baselineskip}{-4pt}
{\footnotesize\em \medskip \noindent infinitesimal \\
 deformations 
\medskip \ }\end{minipage} &
 $\bbh^{1}(X_{T},\liea{f_{T}}{F,\nabla})$ &
  $\bbh^{1}(X_{T},\liea{f_{T}}{E,\theta})$  \\
\hline
\begin{minipage}[c]{1.35in} \addtolength{\baselineskip}{-4pt}
{\footnotesize\em \medskip \noindent infinitesimal \\
obstructions \medskip \ } 
\end{minipage} &
 $\bbh^{2}(X_{T},\liea{f_{T}}{F,\nabla})$ & 
$\bbh^{2}(X_{T}, \liea{f_{T}}{E,\theta})$
\\
\hline
\end{tabular} 
\end{center}

\item There are isomorphisms of Picard stacks over $\stackdrreg{X/S}{n}$
and $\stackdolreg{X/S}{n}$ respectively:
\begin{align*}
T_{\pdr} & = \hh(\tau_{\leq 1}\bbr f_{\op{DR}*}\liea{\fdr}{F_{\op{un}},
\nabla_{\op{un}}}) \\
T_{\pdol} & = \hh(\tau_{\leq 1}\bbr f_{\op{Dol}*}\liea{\fdol}{E_{\op{un}},
\theta_{\op{un}}}) 
\end{align*}
\end{list} 
\end{lem}
{\bf Proof.} Again we will give a proof only for the de Rham
case. Also we will only consider the case when $T = \{s \} \to S$ is a
closed point. The proof carries over verbatim to the case of a general
$T \to S$ and is left to the reader. To simplify notation put $Y :=
X_{T} = X_{s}$.

Let as 
before $B = \op{Spec}(\bbc[\varepsilon]/(\varepsilon^{2}))$ and suppose we are
given a  deformation $(\widetilde{Y},\widetilde{F},\widetilde{\nabla})$ of
$(Y,F,\nabla)$ over $B$. Choose again an acyclic \v{C}ech covering ${\mathfrak
 U}$ of $Y$ consisting of affine open sets over which $F$ trivializes.

The structure sheaf of $\widetilde{Y}$ is an extension of ${\mathcal O}_{Y}$ 
by an ideal of square zero which is isomorphic to ${\mathcal O}_{Y}$ as an
${\mathcal O}_{Y}$-module. Since $Y$ is smooth this extension will
split over any affine open subset due to the infinitesimal lifting
property. Therefore for every $U \in {\mathfrak U}$ we can choose a 
ring isomorphism $c_{U} : {\mathcal O}_{\widetilde{Y}|U} \to {\mathcal O}_{U}
\oplus {\mathcal O}_{U}$ with the ring structure on ${\mathcal O}_{U}
\oplus {\mathcal O}_{U}$ given by $(f,a)\cdot (g,b) := (fg,fb+ga)$. Thus the 
sheaf of rings ${\mathcal O}_{\widetilde{Y}}$ on the topological space $Y$ is
described by the 1-cocycle on the nerve of ${\mathfrak U}$ given by 
$D_{UV} := c_{U}\circ c_{V}^{-1}$. Clearly $D_{UV} \in 
\Gamma(U\cap V, \op{End}_{\bbc}({\mathcal O}_{Y}\oplus {\mathcal O}_{Y}))$
and since $D_{UV}$ is an isomorphism of the corresponding split extensions
we can write it in the form
\[
D_{UV} = \left( \begin{array}{cc} 1 & 0 \\ \del_{UV} & 1 \end{array}\right),
\]
where $\del_{UV}$ is some $\bbc$-linear homomorphism from the first copy of 
${\mathcal O}_{U\cap V}$ to the second one. Also $D_{UV}$ has to be a
ring automorphism of ${\mathcal O}_{U\cap V}\oplus {\mathcal O}_{U\cap V}$ 
for the ring structure described above. This is easily seen to be equivalent 
to $\del_{UV} \in \Gamma (U\cap V, \op{Der}({\mathcal O}_{Y}))$ and hence 
the sheaf of rings ${\mathcal O}_{\widetilde{Y}}$ is described by the cocycle
$\{ \del_{UV} \} \in Z^{1}({\mathfrak U},T_{Y})$ which represents the 
Kodaira-Spencer class of the infinitesimal deformation $\widetilde{Y} \to B$
(see \cite[$0_{\op{IV}}$ 20]{egaiv} for more details). 

Let $h_{U}$, $g_{UV}$, $a_{U}$ be as in 
lemma~\ref{lem-vertical-tangent-spaces}. The triple $(\widetilde{Y},
\widetilde{F},\widetilde{\nabla})$ now is encoded in a pair of cochains
$\{ \tilde{g}_{UV} \} \in C^{1}({\mathfrak U}, \op{GL}(n,{\mathcal 
O}_{\widetilde{Y}}))$ and $\{ \tilde{a}_{U} \} \in C^{0}({\mathfrak U},
\op{End}({\mathcal O}_{\widetilde{Y}}^{\oplus n}))$ satisfying the same 
cocycle and connection conditions (\ref{eq-cocycle}) and 
(\ref{eq-connection}). Again we can write 
\begin{align*}
\tilde{g}_{UV} = & g_{UV} + \varepsilon x_{UV}, \\
\tilde{a}_{U} = & a_{U} + \varepsilon m_{U},
\end{align*}
$\{ m_{U} \} \in C^{0}({\mathfrak U}, \op{End}({\mathcal O}_{Y}^{\oplus n})
\otimes \Omega^{1}_{Y})$, but only this time $\{ x_{UV} \} \in C^{1}(
{\mathfrak U}, \op{End}_{\bbc}({\mathcal O}_{Y}^{\oplus n}))$. Since the 
decomposition $\tilde{g}_{UV} = g_{UV} + \varepsilon x_{UV}$ comes from a
trivialization $\widetilde{F}_{|U} \cong {\mathcal O}_{\widetilde{Y}|U}^{
\oplus n} \stackrel{c_{U}}{\to} ({\mathcal O}_{U} \oplus 
{\mathcal O}_{U})^{\oplus
n}$ we have that $x_{UV} : {\mathcal O}_{U\cap V}^{\oplus n} \to {\mathcal 
O}_{U\cap V}^{\oplus n}$ is a differential operator of order $\leq 1$ and that
$x_{UV} - \del_{UV}\otimes \op{id}_{\bbc^{n}}$ is a ${\mathcal O}_{U\cap 
V}$-linear endomorphism of ${\mathcal O}_{U\cap V}^{\oplus n}$. In other 
words, the Kodaira-Spencer class of the deformation $\widetilde{Y}$ can be 
recovered from $\tilde{g}_{UV}$ as the symbol of the differential operator
$x_{UV}$. Thus $\thick{x}_{UV} := h_{U}^{-1}\circ x_{UV}
\circ h_{V}$ is a section of $\mycal{E}(F)$ and since the identity $x_{UW} = 
g_{UV}x_{VW} + x_{UV}g_{VW}$ is equivalent to $\thick{x}_{UW} = 
\thick{x}_{UV} + \thick{x}_{VW}$ we have $\{ \thick{x}_{UV} \} \in 
Z^{1}({\mathfrak U}, \mycal{E}(F))$. Rewriting again the connection condition
(\ref{eq-connection}) in terms of $\thick{m}_{U} = (h_{U}\otimes\op{id})\circ 
m_{U} \circ h_{U}^{-1}$ we get 
\[
\thick{m}_{U} - \thick{m}_{V} = d^\nabla \circ \thick{x}_{UV} - 
\ell(\thick{x}_{UV}) \circ d^\nabla.
\]
Together with the integrability of $\nabla$ this gives $(\thick{x}_{UV},
\thick{m}_{UV}) \in Z^{1}({\mathfrak U}, \liea{Y}{F,\nabla})$ which 
completes the proof of part (a) of the lemma. 

For part (b) we invoke (similarly to the proof of
Lemma~\ref{lem-vertical-tangent-spaces}) the isomorphism
(\ref{tangent-stack-isomorphism}). Note that unlike part (b) of
Lemma~\ref{lem-vertical-tangent-spaces} no formality is claimed here
since already the deformation-obstruction complexes $\tau_{\leq 1}Rf_{*}T_{f}$ 
need
not be formal. However if the fibers of $f$ have no infinitesimal
automorphisms the complex $\tau_{\leq 1}Rf_{*}T_{f}$ is automatically
formal and we again have 
\begin{align*}
\tau_{\leq 1}\bbr f_{\op{DR}*}\liea{\fdr}{F_{\op{un}},
\nabla_{\op{un}}}) & = [\bbr^{0}f_{\op{DR}*}\liea{\fdr}{F_{\op{un}},
\nabla_{\op{un}}}) \stackrel{0}{\to} \bbr^{1}f_{\op{DR}*}\liea{\fdr}{F_{\op{un}},
\nabla_{\op{un}}})]\\
\tau_{\leq 1}\bbr f_{\op{Dol}*}\liea{\fdol}{F_{\op{un}},
\nabla_{\op{un}}})& = [\bbr^{0} f_{\op{Dol}*}\liea{\fdol}{E_{\op{un}},
\theta_{\op{un}}} \stackrel{0}{\to} \bbr^{1}f_{\op{Dol}*}\liea{\fdol}{F_{\op{un}},
\nabla_{\op{un}}})] 
\end{align*}
due to Simpson formality result \cite[Lemma~3.5]{simpson-higgs} 
\cite[Proposition~10.5]{simpson-moduli2}.
\hfill $\Box$

\bigskip

Next we use the information about the vertical and full tangent spaces of the
fibration  $\pdr : \stackdr{X/S}{n} \to S$, that we have just obtained, to give
a description of the non-abelian Gauss-Manin connection. From now on, until
the end of this section, we will work only with the part $\stackdrreg{X/S}{n}$ 
of $\stackdr{X/S}{n}$ over which the morphism $\pdr$ is smooth. Let 
$(Y;F,\nabla) \in \stackdrreg{X/S}{n}$ be a closed 
point lying over $s \in S$. The non
abelian Gauss-Manin connection is given by a map $GM_{(Y;F,\nabla)} : T_{S,s} 
\to T_{\stackdrreg{X/S}{n},(Y;F,\nabla)}$ splitting the differential 
$d\pdr : T_{\stackdrreg{X/S}{n},(Y;F,\nabla)} \to T_{S,s}$. On the other 
hand,  we know from the proof of Lemma~\ref{lem-full-tangent-spaces} that 
$d\pdr$ fits in the commutative diagram
\[
\xymatrix{T_{\stackdrreg{X/S}{n},(Y;F,\nabla)} \ar[r]^-{d\pdr} 
\ar[d]_-{\kappa_{(Y;F,\nabla)}} & T_{S,s} \ar[d]^-{\kappa_{Y}} \\
\hh(\tau_{\leq 1}\bbr\Gamma(\liea{Y}{F,\nabla})))
\ar[r]_-{R\Gamma(\sigma)} &
\hh(\tau_{\leq 1}R\Gamma(T_{Y}))
}
\]
where $\kappa_{Y}$ and $\kappa_{(Y;F,\nabla)}$ denote the Kodaira-Spencer 
maps for $Y$ and for the triple $(Y;F,\nabla)$ respectively and 
$R\Gamma(\sigma)$ is the map in cohomology induced by the symbol 
morphism of complexes $\sigma : \liea{Y}{F,\nabla} \to [T_{Y} \to 0 \to 
\ldots]$. Recall
next 
that since $\nabla : T_{Y} \to \liea{Y}{F}$ is a connection, it splits the 
symbol map for $F$ and that due to the integrability of $\nabla$ we have 
$\ad{\nabla}\circ \nabla = 0$. Thus $\nabla$ can be thought of as a morphism 
of complexes $[T_{Y} \to 0 \to \ldots] \to \liea{Y}{F,\nabla}$ that splits the 
short exact sequence of complexes
\[
\xymatrix@1{0 \ar[r] & {\lieg{Y}{F,\nabla}} \ar[r] & {\liea{Y}{F,\nabla}} 
\ar[r]^{\sigma} & [T_{Y} \to 0 \to \ldots] \ar[r] \ar@/^1pc/[l]^\nabla 
\ar[r] & 0.
}
\]
Now the intrinsic description of the non-abelian Gauss-Manin connection is 
given by the following lemma.

\begin{lem} \label{lem-GM-smooth} The map $GM_{(Y;F,\nabla)}$ fits in the
commutative diagram
\[
\xymatrix{T_{S,s} \ar[d]_-{\kappa_{Y}} \ar[r]^-{GM_{(Y;F,\nabla)}} &
T_{\stackdrreg{X/S}{n},(Y;F,\nabla)} \ar[d]^-{\kappa_{(Y;F,\nabla)}} \\
\hh(R\Gamma(T_{Y}))
\ar[r]_-{R\Gamma(\nabla)} & 
\hh(\bbr\Gamma(\liea{Y}{F,\nabla}))
}
\]
\end{lem}
{\bf Proof.} Recall the definition of the map $GM_{(Y;F,\nabla)}$ (cf. {\bf 
NA(vii)} and \cite[Section~8]{simpson-moduli2}). Start with a tangent vector 
$v \in T_{S,s}$. Geometrically $v$ corresponds to a morphism $i_{v} : B \to 
S$ mapping the closed point $o \in B$ to $s \in S$. In other words $v$ gives 
a commutative diagram 
\[
\xymatrix@1{Y \ar@{^{(}->}[r]^{j} \ar[d] & X^{v} \ar[d] \\
o \ar[r] & B
}
\]
with $X^{v} := X\times_{S} B$ and $j$-the natural inclusion. Furthermore, 
since $\stackdr{X/S}{n}$ represents the moduli functor of rank $n$ local 
systems on $X/S$, the vector $GM_{(Y;F,\nabla)}(v) \in 
T_{\stackdrreg{X/S}{n},(Y;F,\nabla)}$ can 
be interpreted as a relative rank $n$ local system on $X^{v}/B$ that 
restricts to $(F,\nabla)$ on $Y$. Now it remains only to note that
specifying $(F,\nabla)$ is the same as specifying a structure
$\varphi : p^{*}_{1}F \widetilde{\to} p^{*}_{2}F$ of a crystal of 
vector bundles over $Y$ and that $j^{*}$ is an equivalence between the 
categories of crystals over $X^{v}/B$ and crystals over $Y/B$ 
\cite[Appendix]{grothendieck-crystals}, 
\cite[Proposition~8.4]{simpson-moduli2}.

To make this explicit put $q^{v} : X^{v} \to (Y\times Y)^{\wedge}$ for
the natural map and let $q^{v*}\nabla : F \to F\otimes 
j^{*}\Omega^{1}_{X^{v}/B}$ denote the map $e \mapsto
j^{*}((p^{*}_{2}e - \varphi p^{*}_{2}e) \mod J^{2})$, where as usual 
$J$ stands for the ideal sheaf of $Y \subset (Y\times Y)^{\wedge}$. If 
$U \subset Y$ is an affine open and $c_{U} : {\mathcal O}_{X^{v}|U} \to 
{\mathcal O}_{U}\oplus {\mathcal O}_{U}$ is as in 
Lemma~\ref{lem-full-tangent-spaces}, then $c_{U}$ induces an isomorphism
$\tilde{c}_{U} : \Omega^{1}_{X^{v}/B|U} \to \Omega^{1}_{U} \oplus {\mathcal 
O}_{U}$ that is uniquely characterized by the property $\tilde{c}_{U}
\circ d \circ c_{U}^{-1} (f,a) = (df, a)$. Let ${\mathfrak U}$ be a 
\v{C}ech covering of $Y$ and let $\{\del_{UV} \}$ be the 
${\mathfrak U}$-cocycle representing $\kappa_{Y}(v)$. The sheaf 
$j^{*}\Omega^{1}_{X^{v}/B}$ is completely determined by the cocycle $\{ 
\widetilde{D}_{UV} \} \in Z^{1}({\mathfrak U}, \op{End}(\Omega^{1}_{Y}
\oplus {\mathcal O}_{Y}))$,
where $\widetilde{D}_{UV}((\alpha,a)) = (\alpha, a + i_{\sigma(\del_{UV})}
\alpha)$ and the map $(d^{q^{v*}\nabla})_{|U}$ is given by the morphism 
$F_{U} \to 
(F_{U} \oplus F_{U}\otimes \Omega^{1}_{U})$, $e \mapsto e \oplus \nabla e$. 

Now, by what we said above, $\kappa_{(Y;F,\nabla)}(GM_{(Y;F,\nabla)}(v))$ is 
the Kodaira-Spencer class of the unique (up to isomorphism) relative local 
system $(F^{v},\nabla^{v})$ on $X^{v}/B$ that restricts to $q^{v*}\nabla$.
Finally, the expression for $q^{v*}\nabla$ in terms of ${\mathfrak 
U}$ and Lemma~\ref{lem-full-tangent-spaces} show that $R\Gamma(\nabla)(v)$ 
restricts exactly to $q^{v*}\nabla$. 
\hfill $\Box$

\bigskip

Suppose now we have an algebraic section $a : S \to \stackdrreg{X/S}{n}$. We 
can use the previous lemma to express the property of $a$ being horizontal for 
the non-abelian Gauss-Manin connection in geometric terms. Specifying $a$ is
the same as specifying a relative rank $n$ local system on $X/S$, i.e. a 
vector bundle $F \to X$ of rank $n$ and a relative integrable connection
$d^{\nablaf} : F \to F\otimes \Omega^{1}_{f}$. Denote by 
$\lieg{f}{F,\nablaf}$
and $\liea{f}{F,\nabla_{\! f}}$ the relative versions of complexes in 
Lemma~\ref{lem-vertical-tangent-spaces} and 
Lemma~\ref{lem-full-tangent-spaces} 
respectively. As before, the relative connection $\nablaf$ induces a 
morphism 
of complexes $\nablaf : [T_{f} \to 0 \to \ldots] \to \liea{f}{F,
\nablaf}$ and 
we have the following proposition.

\begin{prop} \label{prop-integrable-section} The section $a$ is 
horizontal with 
respect to the non-abelian Gauss-Manin connection on $\stackdrreg{X/S}{n}$ if
and only if one of the following equivalent conditions holds
\begin{list}{{\em (\arabic{inner})}}{\usecounter{inner}}
\item There exists a global integrable connection $\nabla :
T_{X} \to \liea{X}{F}$ which induces 
\[
\nablaf : T_{f} \to \liea{f}{F}.
\]
\item The following diagram of objects in $D^{b}(S)$ commutes
\[
\xymatrix@1{T_{S} \ar[rr]^-{\kappa_{(X/S; F, \nablaf)}} 
\ar[dr]_-{\kappa_{X/S}} &
& \tau_{\leq 0}({\bbr f_{*}\liea{f}{F,\nablaf}}[1]) \\ & 
\tau_{\leq 0}(Rf_{*}T_{f}[1])
\ar[ur]_-{\tau_{\leq 0}(Rf_{*}\nablaf[1])}
}
\]
\end{list}
\end{prop}
{\bf Proof.} The fact that $(2$) is equivalent to the horizontality of $a$ is
almost immediate. Indeed, according to Remark~\ref{rem-D-schemes}(iii), $a$ is
horizontal if and only if $da = GM$ as maps from $T_{S}$ to $a^{*}
T_{\stackdrreg{X/S}{n}}$. In order to compare $da$ and $GM$ recall
first that due to \cite[Proposition~1.4.15 of Expos\'{e}~XVIII]{sga4} 
the construction $\hh$
induces a $1$-equivalence of the category of of quasi-isomorphism
classes of complexes of sheaves of amplitude one and the $2$-category
of Picard stacks. But by Lemma~\ref{lem-full-tangent-spaces}
we have $a^{*}T_{\stackdrreg{X/S}{n}} = 
\hh(\tau_{\leq 1}(\bbr f_{*}\liea{f}{F,\nablaf}))$ and so we can identify 
$da$ with $\kappa_{(X/S; F, 
\nablaf)}$. Similarly by Lemma~\ref{lem-GM-smooth}
we have $GM = Rf_{*}\nablaf\circ \kappa_{X/S}$.

Consequently, we need only to check the equivalence of $(1)$ and $(2)$. 
The short exact sequence of sheaves
\[
\xymatrix@1{
0 \ar[r] & T_{f} \ar[r] & T^{\sim}_{X} \ar[r]^-{df} & f^{-1}T_{S} \ar[r] & 0,
}
\]
pushes forward to a distinguished triangle in $D^{b}(S)$:
\[
\xymatrix@1{
 Rf_{*}T_{f} \ar[r] & Rf_{*}T^{\sim}_{X} \ar[r]^-{df} &
Rf_{*}f^{-1}T_{S} \ar[r]  & Rf_{*}T_{f}[1],
}
\]
and the Kodaira-Spencer map $\kappa_{X/S} : T_{S} \to \tau_{\leq 0}
(Rf_{*}T_{f}[1])$ is just  the $\tau_{\leq 0}$-truncation of the edge
homomorphism $Rf_{*}f^{-1}T_{S} \to Rf_{*}T_{f}[1]$. Let
$\kappa^{1}_{X/S} : T_{S} \to R^{1}f_{*}T_{f}$ denote the naive
Kodaira-Spencer map, i.e. the composition of $\kappa_{X/S}$ with the
natural morphism 
\[
\tau_{\leq 0}(Rf_{*}T_{f}[1]) \to R^{1}f_{*}T_{f} 
\]
in $D^{b}(S)$.

In order to express $\kappa_{(X/S; F, \nablaf)}$ as an edge homomorphism, we
will just have to recast the local calculation from  
Lemma~\ref{lem-full-tangent-spaces} into the global setting of the family
$f : X \to S$. Following section~\ref{ssec-crystals} put  $\lieat{X}{F} 
\subset \liea{X}{F}$ for all differential operators whose symbol is in 
$T^{\sim}_{X}$. Since $T^{\sim}_{X}$ centralizes $T_{f}$ we have
that $L_{v}(f^{*}\Omega^{1}_{f}) \subset f^{*}\Omega^{1}_{f}$ for all 
$v \in T^{\sim}_{X}$ and hence $L_{\sigma(\del)}$ descends to a 
map in $\op{End}_{\bbc}(\Omega^{1}_{f})$ for a $\del \in \lieat{X}{F}$. Now, 
the same line of reasoning as in Lemma-Definition~\ref{lem-def-full-tangent}
shows that the map $\ellt : \lieat{X}{F} \to \lieat{X}{F\otimes 
\Omega^{1}_{f}}$ is well defined and that $\ad{\nablaf}(\del) := 
\nablaf\circ \del - \ellt(\del)\circ \nablaf$ gives a complex of sheaves on 
$X$:
\[
\lieat{X}{F,\nablaf} := [ 
\xymatrix@1{{\lieat{X}{F}} \ar[r]^-{\ad{\nablaf}} & {\op{End}(F)}\otimes 
\Omega^{1}_{f} \ar[r]^-{\ad{\nablaf}} & {\op{End}(F)\otimes 
\Omega^{2}_{f}} \ar[r] & {\ldots} }].
\]
Finally, due to the calculation in Lemma~\ref{lem-full-tangent-spaces} we 
have that $\kappa_{(X/S; F, \nablaf)}$ is nothing but the 
$\tau_{\leq 0}$-trun\-ca\-tion of the edge 
homomorphism for the distinguished triangle in $D^{b}(S)$:
\[
\xymatrix@1{
{\bbr f_{*}\liea{f}{F,\nablaf}} \ar[r] & {\bbr f_{*}\lieat{X}{F,\nablaf}}
\ar[r]^-{df} & Rf_{*}f^{-1}T_{S} \ar[r] & {\bbr f_{*}\liea{f}{F,\nablaf}[1]}
}
\]
obtained as a push forward of the  short exact 
sequence of complexes
\[
\xymatrix@1{
0 \ar[r] & {\liea{f}{F,\nablaf}} \ar[r] & {\lieat{X}{F,\nablaf}}
\ar[r]^-{df} &
f^{-1}T_{S} \ar[r] & 0.
}
\]
Here, as usual, $f^{-1}T_{S}$ is thought of as a complex concentrated in 
degree zero. Again we put $\kappa^{1}_{(X/S; F, \nablaf)}$ for the naive
Kodaira-Spencer map given as the composition 
$T_{S} \to \tau_{\leq 0}(\bbr
f_{*}\liea{f}{F,\nablaf}[1]) \to \bbr^{1}
f_{*}\liea{f}{F,\nablaf}$.

Furthermore, the two complexes defining $\kappa_{X/S}$ and 
$\kappa_{(X/S;F,\nablaf)}$ are tied up in the following 
commutative diagram 
\begin{equation} \label{eq-lift}
\xymatrix{
& & 0 & 0 \\
& & f^{-1}T_{S} \ar@{=}[r] \ar[u] & f^{-1}T_{S} \ar[u] \\
0 \ar[r] & {\lieg{f}{F,\nablaf}} \ar[r] & {\lieat{X}{F,\nablaf}} \ar[r] 
\ar[u] & 
T^{\sim}_{X} \ar[r] \ar[u] \ar@/^1pc/@{-->}[l]^{\nabla^{\sim}} & 0 \\
0 \ar[r] & {\lieg{f}{F,\nablaf}} \ar[r] \ar@{=}[u] & {\liea{f}{F,\nablaf}} 
\ar[r] 
\ar[u] & 
T_{f} \ar[r] \ar@/^1pc/[l]^\nablaf \ar[u] & 0 \\
& 0 \ar[u] & 0 \ar[u] & 0 \ar[u]}
\end{equation}
Now, we are ready to show the implication $(1) \Rightarrow (2)$. The
condition $(1)$ guarantees the existence of a global 
$\nabla : F \to \liea{X}{F}$ inducing $\nablaf$. The restriction of 
$\nabla$ to $T^{\sim}_{X}$ gives a 
morphism of complexes $\nabla^{\sim} : T^{\sim}_{X}\to \lieat{X}{F,\nablaf}$
that lifts $\nablaf$ in the diagram~\ref{eq-lift}. Thus, the triple of maps 
$(\nablaf,\nabla^{\sim},\op{id}_{f^{-1}T_{S}})$ gives a morphism between the 
third and the second column of diagram (\ref{eq-lift}) and a posteriori a 
morphism between the distinguished triangles in $D^{b}(S)$ which one
obtains after applying $f_{*}$. In particular, this yields $(2)$.

In order to prove that $(2)$ implies $(1)$ it is enough to show that 
the connection $\nablaf$ lifts to a morphism 
$\nabla^{\sim}$ splitting the second row 
of diagram (\ref{eq-lift}). Indeed, if $\nabla^{\sim}$ exists we can extend 
it by ${\mathcal O}_{X}$-linearity to a connection $\nabla : T_{X} \to
\liea{X}{F}$. Furthermore $\nabla$ ought to be integrable for
$\nabla^{\sim}$ is a morphism of complexes. 

Let 
\[
e_{(X/S;F,\nablaf)} \in H^{1}(X, \op{Hom}_{f^{-1}{\mathcal O}_{S}}(
f^{-1}T_{S},\liea{f}{F,\nablaf})) \text{ and }
e_{X/S} \in    H^{1}(X, 
\op{Hom}_{f^{-1}{\mathcal O}_{S}}(f^{-1}T_{S},T_{f}))
\] 
denote the extension 
classes of the second and the third column of (\ref{eq-lift}) respectively.
The lifting $\nabla^{\sim}$ of $\nablaf$ will exist if and only if the 
push-forward of the extension class $e_{X/S}$ via the morphism $\nablaf :
T_{f} \to \liea{f}{F,\nablaf}$ coincides with $e_{(X/S;F,\nablaf)}$. 

Observe that the Lerray spectral sequence gives a diagram
\[
\xymatrix{
0 \ar[d] & 0 \ar[d] \\
H^{1}(S,\Omega^{1}_{S}\otimes f_{*}T_{f}) \ar[d] 
\ar[r]^-{h^{1}(f_{*}\nablaf)} & H^{1}(S,\Omega^{1}_{S}\otimes 
{\bbr^{0}f_{*}\liea{f}{F,\nablaf})} 
\ar[d] \ar@/^1pc/[l]^-{h^{1}(f_{*}\sigma_{f})} \\
H^{1}(X, \op{Hom}_{f^{-1}{\mathcal O}_{S}}(
f^{-1}T_{S},T_{f})) \ar[d]^-{q} \ar[r]^-{h^{1}(\nablaf)} &
 H^{1}(X, \op{Hom}_{f^{-1}{\mathcal O}_{S}}(f^{-1}T_{S},\liea{f}{F,\nablaf})) 
\ar[d]_-{Q} \ar@/^1pc/[l]^-{h^{1}(\sigma_{f})} \\
H^{0}(S,\Omega^{1}_{S}\otimes R^{1}f_{*}T_{f}) \ar[r]^-{R^{1}f_{*}\nablaf} &
H^{0}(S,\Omega^{1}_{S}\otimes \bbr^{1}f_{*}\liea{f}{F,\nablaf}) 
\ar@/^1pc/[l]^-{R^{1}f_{*}\sigma_{f}}
}
\]
which is commutative for both $\sigma_{f}$ and $\nablaf$.
Set $e = h^{1}(\nablaf)(e_{X/S}) - e_{(X/S;F,\nablaf)}$.  Since $q(e_{X/S}) 
= \kappa^{1}_{X/S}$ and $Q(e_{(X/S;F,\nablaf)}) =
\kappa^{1}_{(X/S;F,\nablaf)}$ we 
have that 
\[
Q(e) = R^{1}f_{*}\nablaf\circ h^{0}(\kappa_{X/S}) -
h^{0}(\kappa_{(X/S;F,\nablaf)})
= 0
\]
by condition $(2)$. Thus $e \in H^{1}(S,\Omega^{1}_{S}\otimes 
\bbr^{0}f_{*}\liea{f}{F,\nablaf})$. Due to the smoothness of the stack
$\stackdrreg{X/S}{n}$ we have an equality 
\[
\op{Ext}^{1}_{{\mathcal
O}_{S}}(T_{S}, \bbr^{0}f_{*}\liea{f}{F,\nablaf}) 
= \op{Ext}^{1}_{{\mathcal
O}_{S}}(\ker(\kappa^{1}_{(X/S;F,\nablaf)}),
\bbr^{0}f_{*}\liea{f}{F,\nablaf}).
\] 
On the other hand the condition $(2)$ implies that
the image of $e$ in the group 
\[
\op{Ext}^{1}_{{\mathcal
O}_{S}}(\ker(\kappa^{1}_{(X/S;F,\nablaf)}),
\bbr^{0}f_{*}\liea{f}{F,\nablaf})
\] 
is equal to zero which proves 
the proposition. 
\hfill $\Box$

\bigskip

\begin{rem} \label{rem-global-lifting} {\bf (i)} The cohomology
class $e$ which appears in the proof of 
Proposition~\ref{prop-integrable-section} is the obstruction for the 
relative connection $(F,\nablaf)$ to lift to a connection on the  whole $X$.
To make this somewhat more explicit observe first that $e$ lies even deeper 
inside $H^{1}(X, 
\op{Hom}_{f^{-1}{\mathcal O}_{S}}(f^{-1}T_{S},\liea{f}{F,\nablaf}))$.
Indeed, diagram (\ref{eq-lift}) implies that 
\[
e_{X/S} =
h^{1}(\sigma_{f})(e_{(X/S;F,
\nablaf)})
\] 
which combined with the identity $\sigma_{f}\circ 
\nablaf = \op{id}$ gives $h^{1}(\sigma_{f})(e) = 0$. Thus $e \in 
H^{0}(S,\Omega^{1}_{S}\otimes \bbr^{0}f_{*}\lieg{f}{F,\nablaf})$. It is 
easy to represent $e$ by a \v{C}ech cocycle. Let ${\mathfrak U}$ be a \v{C}ech
covering of $S$. Due to Proposition~\ref{prop-integrable-section} 
for any $U \in {\mathfrak U}$ we can choose a lifting $\nabla^{\sim}_{U} :
T^{\sim}_{X_{U}} \to \lieat{X_{U}}{F,\nablaf}$ of $\nablaf$. The difference
$e_{UV} := \nabla^{\sim}_{V} - \nabla^{\sim}_{U}$ for two $U, V \in {\mathfrak
 U}$ is trivial on $T_{f} \subset T^{\sim}_{X}$ and hence belongs to
$\Gamma(X_{U\cap V}, \op{Hom}_{f^{-1}{\mathcal O}_{S}}(f^{-1}T_{S},
\lieg{f}{F,\nablaf}) = \Gamma(U\cap V, \Omega^{1}_{S}\otimes 
\bbr^{0}f_{*}\lieg{f}{F,\nablaf})$. In particular, the cocycle 
$\{ e_{UV} \} \in Z^{1}({\mathfrak U}, \Omega^{1}_{S}\otimes 
\bbr^{0}f_{*}\lieg{f}{F,\nablaf})$ represents $e$.

\medskip

\noindent
{\bf (ii)} The proof of Proposition~\ref{prop-integrable-section} explains
also the geometric meaning of the natural weakening of condition $(2)$
which takes into account only the naive parts of the Kodaira-Spencer
maps. Namely we have the equivalence of the following two conditions
\begin{itemize}
\item[$(1')$] {\em Locally in $S$ there exists a global integrable connection 
$\nabla :
T_{X} \to \liea{X}{F}$ which induces }
\[
\nablaf : T_{f} \to \liea{f}{F}.
\]
\item[$(2')$] {\em The following diagram of sheaves on $S$
\[
\xymatrix@1{T_{S} \ar[rr]^-{\kappa^{1}_{(X/S; F, \nablaf)}} 
\ar[dr]_-{\kappa^{1}_{X/S}} &
& {\bbr^{1} f_{*}\liea{f}{F,\nablaf}} \\ & 
R^{1}f_{*}T_{f})
\ar[ur]_-{R^{1}f_{*}\nablaf}
}
\]
commutes.}
\end{itemize}
The implication $(1') \Rightarrow (2')$ follows in exactly the same way
as in the proof of Proposition~\ref{prop-integrable-section}. For the
proof $(2') \Rightarrow (1')$ one only has to notice that when $S$ is
replaced by an affine open n $U \subset S$ and when $X$ is replaced by
$X\times_{S} U$ the obstruction class $e$ vanishes since 
$H^{1}(U,\Omega_{U}\otimes \bbr^{0}f_{*}\liea{f}{F,\nablaf}) = 0$.

\end{rem}

\subsection{General variations - some speculations} 
\label{ssec-general-variations}

In this section we probe a possible general framework for non-abelian
Hodge theory in the "weight one" case.
The section is not used in the rest of the paper and may be skipped  by
the reader.

The results of Simpson discussed in Section~\ref{sec-nahodge}{\bf
(NA(i-vi))} suggest the  following general notion which is essentially
due to Simpson:

\begin{defi} \label{defi-abstract-nah}
A {\em polarized complex non-abelian Hodge structure of 
weight one} consists of the following data
\begin{description}
\item[$\bbc$NAH1](Space) A complex algebraic variety (or stack) $M$.
\item[$\bbc$NAH2](Hodge filtrations) A variety $Z$ equipped with:  
\begin{itemize}
\item an algebraic  $\cstar$-action $\gamma : \cstar\times Z \to Z$ 
such that for any closed point $z \in Z$ the limits $\lim_{t \to 0}
\gamma_{t}(z)$ exists in $Z$;
\item a $\cstar$-equivariant morphism $\zeta : Z \to {\mathbb P}^{1}$;
\item an isomorphism $Z_{1} \cong M$.
\end{itemize}
\item[$\bbc$NAH3](Opposedness of the Hodge filtrations) There exists a real
analytic trivialization $\phi : Z \to M\times {\mathbb P}^{1}$ such that
\begin{itemize}
\item The $\phi$-constant sections of $\zeta : Z \to {\mathbb P}^{1}$
are holomorphic;
\item If $D(\zeta)$ denotes the Douady space of sections of
$\zeta$ and if $M^{\phi} \subset D(\zeta)$ is the subset of all
$\phi$-constant sections, then there exists a neighborhood $M^{\phi}
\subset \mycal{S} \subset D(\zeta)$ such that for every $x, y \in
{\mathbb P}^{1}$ the natural evaluation map 
\[
\op{ev}_{x,y} : D(\zeta) \to Z_{x}\times Z_{y}
\]
induces an isomorphism of $\mycal{S}$ and a neighborhood of the
diagonal.
\end{itemize}
\item[$\bbc$NAH4](Polarization) An algebraic relative form
$\Omega \in H^{0}(Z,\Omega^{2}_{\zeta}\otimes \zeta^{*}{\mathcal O}_{\mathbb
P}^{1}(2))$ satisfying $\gamma_{t}^{*}\Omega = 
t\Omega$ for all $t \in \cstar$ and such that the restriction 
$\Omega_{t} := \Omega_{|Z_{t}}$ is a symplectic form when restricted 
to the smooth locus of any fiber.
\item[$\bbc$NAH5](Splitting of the Hodge filtrations)   
Two morphisms
$h_{0} : Z_{0} \to B_{0}$ and $h_{\infty} : Z_{\infty} \to B_{\infty}$
whose fibers are Lagrangian for $\Omega_{0}$ and $\Omega_{\infty}$
respectively and generically transversal to the closures of the
$\cstar$-orbits.
\end{description}

\medskip

A polarized complex non-abelian Hodge structure will be
called {\em real} if in addition $Z$ is equipped with an antiholomorphic
involution $\sigma : Z \to Z$ covering the antipodal involution on
${\mathbb P}^{1}$ so that $\zeta$ and $\Omega$ are real and
$\mycal{S}^{\sigma} = M^{\phi}$.
\end{defi}

\begin{rem} \label{rem-abstract-nahs} (i) In fact Simpson
\cite[Section~5]{simpson-santa-cruz} defines 
a filtration of a scheme (or a stack) $M$ as a ${\mathbb G}_{m}$-scheme 
(or stack) $Z$ equipped with a ${\mathbb G}_{m}$-equivariant map to
${\mathbb A}^{1}$ and an identification $Z_{1} \cong M$. In view of this
condition $\bbc${\bf NAH2} can be though of as giving two filtrations on
$M$ - one corresponding to $Z_{|{\mathbb P}^{1}\setminus\{0\}}$ and the
other to $Z_{|{\mathbb P}^{1}\setminus\{\infty\}}$

\medskip

\noindent
(ii) The notion of an abstract
non-abelian Hodge structure of weight one proposed above specializes to
some well known geometric objects. For example if 
$M$ is smooth, then a real non-abelian Hodge structure of weight one is nothing but the
but a hyperk\"{a}hler structure with $\cstar$-equivariant twistor
space\footnote{This shouldn't be confused with the much stronger notion
of a hyperk\"{a}hler cone used say in \cite{ranee-instantons}.} This follows
simply from the twistor interpretation of hyperk\"{a}hler structures from
\cite{hitchinetal}. In fact in this case one can replace condition
$\bbc${\bf NAH3} by the much simpler to check (but equivalent condition)
that $\zeta$ has at least one $\sigma$-invariant holomorphic section 
with normal bundle
isomorphic to ${\mathcal O}_{{\mathbb P}^{1}}(1)\otimes \bbc^{\dim
Z/{\mathbb P}^{1}}$. 

Similarly if we drop conditions $\bbc${\bf
NAH4} and  $\bbc${\bf NAH5} and the $\cstar$-equivariance condition 
from $\bbc${\bf
NAH2}  then we recover the notion of a (possibly singular) quaternionic
variety introduced by Deligne in \cite{deligne-twistor-letter}.

\medskip

\noindent
(iii) The twistor family $\zeta : Z \to {\mathbb P}^{1}$ corresponding
to a non-abelian Hodge structure has many peculiar geometric properties.

For example the assumption that $Z$ is algebraic immediately implies
that $Z$ and a posteriori $M$ cannot be proper. 
In addition the $\cstar$-equivariance required in $\bbc${\bf NAH2}
implies that $\zeta$ trivializes holomorphically 
over ${\mathbb P}^{1}\setminus\{ 0, \infty \}$. In particular $M = Z_{1}$,
$Z_{0}$ and $Z_{\infty}$ are the only non-isomorphic members of the
family $\zeta$. 

If in addition we assume that $M$ is smooth then the condition that
for every $z \in Z$ the limit $\lim_{t \to
0}tz$ exists, then by \cite[Section~4]{kaledin-cotangent} we see that 
$Z_{0}$ is
birational to the total space $T$ of the cotangent bundle of
$Z_{0}^{\cstar}$ and that moreover the symplectic form
$\Omega_{0}$  transforms into the standard symplectic form on
$T$ at least over a big open set. In particular
the universal categorical quotient of $Z_{0}$ by $\cstar$ exists and is
birational to $Z_{0}^{\cstar}$. By the
$\cstar$-equivariance of $\zeta$ then it is clear that the  universal 
categorical
quotient of $Z_{|{\mathbb P}^{1}\setminus \{0\}}$ by $\cstar$ will
exist and will be birational to $Z_{0}^{\cstar}$ as well. Thus we get a
rational map $Z_{|{\mathbb P}^{1}\setminus \{0\}} \to
Z_{0}^{\cstar}\times {\mathbb P}^{1}$ which furnishes the family of
algebraic symplectic family $Z_{|{\mathbb P}^{1}\setminus \{0\}} \to
{\mathbb P}^{1}$ with a Lagrangian foliation which on the fiber $Z_{0}$
is generically transversal to $h_{0}$. Similar analysis holds
over ${\mathbb P}^{1}\setminus \{\infty \}$.

This structure of $Z$ is consistent with the picture we got in
Section~\ref{ssec-abelian-non-abelian-hodge}{\bf (NA(iv))}.
\end{rem}

\bigskip

\noindent 
Definition~\ref{defi-abstract-nah} seems (and is in fact) quite
restrictive. Some obvious examples of non-abelian Hodge structures of
weight one are:

\begin{ex} \label{ex-quaternionic} Every polarized 
complex hodge structure of weight one (= a regular equivariant 
quaternionic vector space \cite{simpson-mixed-twistor}
\cite[Section~1.1.7]{kaledin-cotangent}) is
also a non-abelian polarized Hodge structure of weight one.
\end{ex}

\begin{ex} \label{ex-projective} Let $V$ be an $n+1$-dimensional
real vector space and let $V_{\bbc} : V\otimes_{\bbr}\bbc$ be its
complexification. Put $P :=
\op{Proj}(S^{\bullet}V_{\bbc}^{\vee})$ for the projectivization of
$V_{\bbc}$.
Take $M
\subset P\times P^{\vee}$ to be the affine subvariety
defined as $M = \{ (x,h) | x \not\in h \}$. The Euler sequence on
$P$ identifies $M$ with the twisted cotangent bundle of
$P$ corresponding to the hyperplane class $c_{1}({\mathcal
O}_{P}(1)) \in H^{1}(P, \Omega^{1}_{P})$. In particular there is a
tautological holomorphic family $Z^{+} \to {\mathbb A}^{1} = H^{1}(P,
\Omega^{1}_{P})$ of $T^{\vee}_{P}$-torsors on $P$. By construction we
have
\[
Z^{+}_{|{\mathbb A}^{1}\setminus\{0\}} \cong {\mathbb
A}^{1}\setminus\{0\}\times M \text{ and } Z^{+}_{0} = T^{\vee}_{P}.
\]
Moreover we can use the real structure on $V_{\mathbb C}$ and the
antipodal real structure on $\bbp^{1}$ to glue $Z^{+}$ with its
conjugate family $Z^{-}$ and obtain a $\cstar$-equivariant family
$\zeta : Z \to \bbp^{1}$. The family $\zeta$ is the twistor family of
the hyperk\"{a}hler manifold $T^{\vee}_{P}$ which is obtained from the
quaternionic space ${\mathbb H}^{n}$ by hyperk\"{a}hler reduction by
$U(n)$ \cite{calabi}. In particular $\bbc${\bf NAH1-4} will
automatically hold. To see that $\bbc${\bf NAH5} holds observe that
$P$ is a toric variety and so we have a Poisson action of 
$(\cstar)^{n}$ on $T^{\vee}_{P}$ for the standard algebraic symplectic
form. Define $h_{0}
: T^{\vee}_{P} \to \bbc^{n}$ to be just the moment map for the
$(\cstar)^{n}$-action. It is well known that the fibers of $h_{0}$ are
coisotropic (see e.g. \cite[Section~1.5]{ginzburg-book}) and it is not
hard to see that the fiber of $h_{0}$ over $0 \in \bbc^{n}$ is a
(reducible) subvariety of dimension $n$. Thus the generic fiber of
$h_{0}$ is Lagrangian and generically transversal to the fibers of
$T^{\vee}_{P}  \to P$.
\end{ex}

\begin{ex} \label{ex-abelian-nah} Similarly to
Example~\ref{ex-projective} one can start not with a complex
projective space but with a principally polarized abelian variety
$(A,\theta)$. In this case one takes  $M$ to be the twisted 
cotangent bundle of $A$ corresponding to the class $c_{1}(\theta) \in
H^{1}(A,\Omega^{1}_{A})$ and $Z^{+} \to \bbc\theta \subset
H^{1}(A,\Omega^{1}_{A})$ to be again the universal torsor. Again
$Z^{+}$ glues with its complex conjugate to produce the twistor space
of the hyperk\"{a}hler manifold $T^{\vee}_{A}$. The hyperk\"{a}hler structure
on $T^{\vee}_{A}$ exists  since the K\"{a}hler metric
corresponding to $\theta$ is flat (see \cite{calabi} or
\cite{biquard-gauduchon-cotangent}). Finally the splitting of the Hodge
filtrations follows from the fact that $T_{A}^{\vee}$ is a trivial
bundle and so we can take $h_{0}$ to be the projection onto the fiber
$T^{\vee}_{A,0}$ at $0$.
\end{ex}

\begin{ex} \label{ex-geometric-nah} Given a smooth projective
variety $X/\bbc$ and a complex reductive algebraic group $G$, then we
get a non-abelian Hodge structure by taking $M = \hdr{1}{X}{G}$, $Z = 
\hdel{1}{X}{G}$ and the twistor lines come from the choice of harmonic
metric (these are called the preferred sections in
\cite{simpson-santa-cruz}). As explained in detail in 
\cite{simpson-santa-cruz} (see also
Section~\ref{ssec-abelian-non-abelian-hodge}{\bf (NA(vi))}) the fibers
of $\hdel{1}{X}{G}$ over $0$ and $\infty$ are naturally identified
with $\hdol{1}{X}{G}$ and $\overline{\hdol{1}{X}{G}}$
respectively and the splitting of the Hodge filtrations is given by
the Hitchin map at $0$ and its complex conjugate at $\infty$.
\end{ex}

\

\bigskip

\noindent
Note that in fact all of the previous examples are to some extend
special cases of Example~\ref{ex-geometric-nah}. Indeed
Examples~\ref{ex-quaternionic} and \ref{ex-abelian-nah} correspond to
e.g. taking $X$  to be an abelian variety of dimension $g$ (or a
smooth curve of genus $g$) and $G = \bbc$ or $G = \cstar$
respectively. For Example~\ref{ex-projective} a slightly different
approach is necessary. Fix an elliptic curve $E$ and let $G =
\op{SL}_{n}(\bbc)$. In this case the moduli space
$H^{1}(E,\op{SL}_{n+1}({\mathcal O}_{E}))$ of semistable vector
bundles 
of rank $n+1$
and degree $0$ can be naturally identified (see
e.g. \cite{tu-elliptic}) with the symmetric product
$E^{(n)}$ trough the Fourier-Mukai transform on $E$. On the other
hand one can use the Abel-Jacobi map for $E$ to identify $E^{(n)}$
with an $n$-dimensional complex projective space $P$ and so we have an
$S_{n}$-Galois cover $\pi : E^{\times n} \to P$, where $S_{n}$ denotes
the symmetric group on $n$ letters. It is clear that the moduli space
$\hdol{1}{E}{\op{SL}_{n+1}(\bbc)}$ of Higgs bundles  
is birational\footnote{In fact for $n = 1$ the moduli space
$\hdol{1}{E}{\op{SL}_{2}(\bbc)}$ can be naturally identified with the
Hilbert quotient of $T^{\vee}E$ by $\{\pm 1\}$ i.e. with the blow up of
the ordinary quotient  $T^{\vee}E/\{\pm 1\}$ at its four ordinary
double points. It is reasonable to expect therefore that more generally
$\hdol{1}{E}{\op{SL}_{n+1}(\bbc)}$ will be just the Hilbert quotient
of $T^{\vee}E^{\times n}$ by $S_{n}$ but this is irrelevant for our
discussion.}  to the $S_{n}$-quotient of the cotangent bundle to the abelian
variety $E^{\times n}$ and hence birational to $T^{\vee}_{P}$. This
gives a hyperk\"{a}hler identification of big open sets of $T^{\vee}_{P}$
and $\hdol{1}{E}{\op{SL}_{n+1}(\bbc)}$. The corresponding splittings
of the Hodge filtrations seem to be different however. Indeed the
Lagrangian fibration for $\hdol{1}{E}{\op{SL}_{n+1}(\bbc)}$ is given by the
Hitchin map and hence has fibers abelian varieties. In contrast the
Lagrangian fibration for $T^{\vee}_{P}$ has toric varieties as fibers.

\bigskip

\noindent
In view of the above discussion it seems reasonable to try to relax the
conditions $\bbc${\bf NAH1-5}. For example instead of working with a
$\cstar$-scheme (or stack) $M$ one may work with a formal scheme (or a
formal stack) equipped with a $\cstar$ action. Similarly it seems
reasonable to only require the existence of rational
maps $h_{0} : Z_{0} \dashrightarrow B_{0}$ and $h_{\infty} :
Z_{\infty} \dashrightarrow B_{\infty}$ with Lagrangian fibers. With
these relaxations one can now speculate about the existence of more
exotic examples.

\begin{ex-speculation} \label{ex-cotangent-flags}  Let $G$ be a 
complex semi simple group and
let $P \subset G$ be a parabolic. It is known
\cite[Theorem~1]{biquard-gauduchon-cotangent} (see also 
\cite{kaledin-cotangent} and \cite{feix})
that for any homogeneous K\"{a}hler metric
$g$ on the partial flag variety $G/P$  there exists a unique
$G$-invariant hyperk\"{a}hler metric on $T^{\vee}_{G/P}$ which restricts
to $g$ on the zero section. As shown in
\cite[Section~3.3]{biquard-twistors} the twistor space $Z$ of such a
hyperk\"{a}hler metric again has the property that
$Z_{|\bbp^{1}\setminus\{\infty\}}$ is the tautological family of
$\Omega^{1}_{G/P}$-torsors with class $[\omega] \in
H^{1}(G/P,\Omega^{1}_{G/P})$ where $\omega$ is the K\"{a}hler class of
$g$. In particular we can take $M$ to be the twisted cotangent bundle
$T^{\vee}_{G/P}([\omega])$. In fact this description  can be used
\cite{mirkovic-hyperkahler} as a starting point for a completely
algebraic description of the hyperk\"{a}hler structure on $T^{\vee}_{G/P}$
which is especially well suited for inducing hyperk\"{a}hler structures on
associated spaces.

To interpret this data as a non-abelian
Hodge structure of weight one we also need to specify the splitting of
the Hodge filtrations. For this we may proceed in two ways. One
alternative is to replace $Z_{0} = T^{\vee}_{G/P}$ with a formal
neighborhood of the zero section in $T^{\vee}_{G/P}$. In that case $M$
should be replaced by the scheme parameterizing all pairs $(x,\gamma)$
where $x \in G/P$ and $\gamma$ is a section of
$T^{\vee}_{G/P}([\omega]) \to G/P$ over a formal neighborhood of $x
\in G/P$ (compare with Remark~\ref{rem-D-algebras}(ii)). 
In this case the splitting of the Hodge filtration on say
$Z_{0}$ will be just the restriction of the moment map $\mu :
T^{\vee}_{G/P} \to \op{Lie}(G)^{\vee}$ for the standard $G$-action and
symplectic form on $T^{\vee}_{G/P}$. Another possibility is to keep
$Z$ as it is and to look for a Lagrangian foliation transversal to the
fibers of the cotangent bundle. One choice might come from taking 
preimages of $P$-coadjoint orbits via the moment map $\mu_{P} :
T^{\vee}_{G/P} \to \op{Lie}(P)^{\vee}$.
\end{ex-speculation}

\begin{ex-speculation} \label{ex-coadjoint-orbits} A variation of the
previous example is to take $Z$ to be the twistor space of the
Kronheimer hyperk\"{a}hler metric \cite{kronheimer-hyperkahler},
\cite{biquard-twistors} of a coadjoint orbit ${\mathbb O} \subset
{\mathfrak g}^{\vee}$ of a complex reductive Lie algebra ${\mathfrak g}$.
In this case the twistor fiber $Z_{1}$ is naturally identified with
the moduli space of solutions of the twisted complex Nahm equations
 or equivalently with the moduli space of
rotationally symmetric logarithmic $\lambda$ connections on
$\cstar$ \cite{biquard-twistors}.

The latter interpretation suggests that it may be possible to define 
the splitting map $h_{0} :  {\mathbb O} \to B_{0}$ by working with
Jacobians of spectral covers of logarithmic spaces. This is a very
interesting question which merits serious consideration.
\end{ex-speculation}

\begin{ex-speculation} \label{ex-hilbert-scheme} Let $C$ be a smooth
curve of genus $g$ and let $Y = \op{tot}(T^{\vee}_{C})$ be the total
space of its cotangent bundle. The surface $Y$ is naturally
holomorphic symplectic and so the Hilbert scheme $Y^{[m]}$ of
$0$-dimensional schemes of length $m$ on $Y$ will be a smooth
holomorphic symplectic manifold. In fact due to
\cite[Theorem~1]{kaledin-cotangent} one knows that $Y$ carries an
incomplete $U(1)$-equivariant hyperk\"{a}hler metric and so one expects
that $Y^{[m]}$ is hyperk\"{a}hler as well. The analysis of
\cite{kaledin-cotangent} and the twistor space
construction in the forthcoming thesis \cite{feix}
show that at least when one works with the formal completion of the
zero section in $Y$ the corresponding twistor family behaves exactly as
required by $\bbc${\bf NAH1-4}. More precisely given any smooth
K\"{a}hler manifold $(V,\omega)$ whose K\"{a}hler
metric is real analytic \cite[Theorem~1]{kaledin-cotangent} and
\cite{feix} construct a unique $U(1)$-equivariant
hyperk\"{a}hler structure on a tubular neighborhood of the zero
section of $T^{\vee}_{V}$ which restricts to $\omega$. Furthermore it
is shown in \cite{feix} (and implicitly in \cite{kaledin-cotangent})
that the twistor space of this hyperk\"{a}hler 
structure is trivial over $\cstar$ with a fiber which is biholomorphic
to a tubular neigborhood $N$ of $V$ in the $\omega$-twisted cotangent
bundle $T^{\vee}_{\omega,V}$. Here $V \subset
\op{tot}(T^{\vee}_{\omega,V})$ 
is emebedded as a real-analytic submanifold via the section corresponding
to the K\"{a}hler form $\omega$. Moreover Kaledin's interpretation of
the points of $N$ as regular extended connections on $V$ shows that
$N$ can be identified with a tubular neighborhood of the diagonal of
$V\times V$ and so the twistor space of $N$ at infinity can be
interpretted as the complex conjugate of the deformation of $V \subset
V\times V$ to its normal cone. In these setup the twistor lines for
$N$ become just the intrinsic holomorphic exponential map for the
K\"{a}hler manifold $V$ as described by \cite[Section~2.9]{kapranov}
(see also \cite{bershadskyetal}, \cite{calabi}).

One might also speculate that there should
be a natural Lagrangian foliation of $Y^{[m]}$ which is
transversal to the fibers of the Lagrangian fibration $Y^{[m]}$ and can
serve as a splitting of the Hodge filtrations (at least in the formal case).

This is supported further by the beautiful analysis of \cite{hurtubise}
(see also \cite[Chapter~7]{nakajima-hsbook})
which shows that for infinitely many choices of $m$ the Hilbert scheme
$Y^{[m]}$ is birational to the moduli space
$\hdol{1}{C}{\op{GL}_{n}(\bbc)}$ and that the holomorphic symplectic
forms match on a big open set after this birational identification.
\end{ex-speculation}

\bigskip

\noindent
Following \cite[Section~8]{simpson-santa-cruz}
one can speak of variations of non-abelian Hodge
structures of weight one. By definition this is a morphism of varieties
$\zeta: Z \to S\to {\mathbb P}^{1}$, so that 

\begin{itemize}
\item $Z$ is equipped with a $\cstar$-action covering the standard
action on ${\mathbb P}^{1}$ and such that for every $s \in S$ the
restriction $Z_{|\{ s\}\times \bbp^{1}}$ satisfies $\bbc${\bf
NAH1-4}.
\item There exist schemes $B_{0} \to S$ and $B_{\infty} \to S$ and
morphisms $h_{0} : Z_{|S\times \{0\}} \to B_{0}$ and 
$h_{\infty} : Z_{|S\times \{\infty\}} \to B_{\infty}$ over $S$
specializing to splittings of the Hodge filtrations for every $s \in
S$.
\item The family $Z_{|S\times \{1\}} \to S$ is a local system of
schemes on $S$, i.e. $Z_{|S\times \{1\}}$ is equipped with an action of
the formal groupoid $S_{\op{DR}}$ over $S$ \ref{rem-deRham-groupoid}. 
\item Griffiths transversality holds, i.e. the natural action of the
formal groupoid $S_{\op{DR}}$ on $Z_{|S\times \{1\}}$ extends to an
action of the Hodge formal groupoid $S_{\op{Hod}}$ on $Z_{|S\times
(\bbp^{1}\setminus \infty)}$. Here $S_{\op{Hod}} := \xymatrix@1{
N \ar@<.5ex>[r]\ar@<-.5ex>[r]  & S\times {\mathbb A}^{1}
}$ with $N$ being the formal completion of $S$ inside its deformation
to  the normal cone.
\end{itemize}

\bigskip

\noindent
Moreover given a non-abelian Hodge structure $(M, Z \to
\bbp^{1},\phi,\Omega,h_{0},h_{\infty})$  of weight one we can talk of 
non-abelian $(p,p)$ classes in exactly the
same way as before. Namely these are the closed points in $M$ that
correspond to $\cstar$ fixed points in $Z_{0}$ under the
trivialization $\phi$. 

In particular it makes sense to ask if Theorem~\ref{thm-main} will hold in
this more general context. We formulate this as a question for further
study:

\begin{que} \label{que-pp-general} Let $Z \to S\times {\mathbb P}^{1}$
be a variation of non-abelian Hodge structures of weight one over a
smooth quasi projective $S$. Assume that we are given 
section $a$ of $Z_{|S\times \{ 1 \}} \to S$ which is horizontal for
the non-abelian connection, i.e. which comes from a section of 
$Z_{|S\times \{ 1 \}} \to S_{\op{DR}}$. Asume further that for some
point $o \in S$ the point $a(o)$ is a
non-abelian $(p,p)$ class in $Z_{(o,1)}$ is it true that $a(s)$ is a
non-abelian $(p,p)$-class for all $s \in S$?
\end{que}

It will be very interesting to investigate this question for some
example which is not of geometric origin. For this one will need to
build interesting variations of non-abelian Hodge structures. For
almost all examples considered above this can be done in the obvious
manner. In that respect Example-Speculation~\ref{ex-hilbert-scheme}
looks especially promising. Indeed one expects that for any smooth
family of curves $C \to S$ and any global K\"{a}hler metric on $C$ the
corresponding family of hyperk\"{a}hler structures produced by
Kaledin's construction will constitute a variation of non-abelian
Hodge structures of weight one. Question~\ref{que-pp-general} becomes
very intriguing in this setup since Nakajima had shown
\cite[Proposition~7.5]{nakajima-hsbook} that 
\[
(Y^{[n]})^{\cstar} = \coprod_{\nu} C^{(\nu)}
\]
where $\nu$ runs over all partitions of $n$ and 
\[
C^{(\nu)} :=
C^{(\alpha_{1})}\times \ldots \times C^{(\alpha_{n})} \text{ for a partition } 
\nu = (1^{\alpha_{1}},2^{\alpha^{2}}, \ldots, n^{\alpha_{n}}).
\]

Similarly one might attempt to construct examples of non-abelian
variations starting from
Example-Speculation~\ref{ex-coadjoint-orbits}. One possibility is to
work with a local system of complex reductive Lie algebras and to 
consider the
minimal nilpotent orbit in each fiber. Since the Lie bracket is
assumed to be horizontal for the corresponding connection it is clear
then that the connection will preserve the minimal nilpotent
orbit. Thus the family of minimal nilpotent orbits will give us a
local system of schemes. Moreover the analysis carried out in 
\cite{ranee-instantons} \cite{ranee-quantization} shows that in
several 
cases the zero fiber of the
twistor space is isomorphic (at least locally in the \'{e}tale
topology) to the total space of the holomorphic cotangent bundle of 
the real minimal nilpotent orbit
equipped with its natural K\"{a}hler structure. Another possibility
suggested from the geometric quantization approach
\cite{ranee-instantons} is to work with a fixed complex nilpotent
orbit but to vary the holomorphic polarization on it. Then from the
geometric quantization point of view one expects a flat connection on
on the space of bundle of quantization spaces over the parameter space
of all polarizations. Since the quantization spaces are typically
spaces of functions (or sections in a line bundle) on the variety it
is not unreasonable to look for a connection on the actual family of
polarized orbits. One such connection with finite monodromy is
explicitly described in \cite{ranee-sl3}.

\section{The main theorem} \label{sec-main}

As explained in Section~\ref{ssec-abelian-non-abelian-hodge}{\bf (NA(v))} 
Simpson
characterizes \cite{simpson-higgs} the
non-abelian $(p,p)$ classes as local systems whose Higgs bundles are
fixed under the natural $\cstar$-action. Thus
Theorem~\ref{thm-main} can be viewed as a statement  comparing the
non-abelian Gauss-Manin connection on the stack
$\stackdr{X/S}{n}$ of relative local systems with the $\cstar$-action of 
on the stack $\stackdol{X/S}{n}$  of relative Higgs bundles. It is hard
to work out such a comparison in practice because the homeomorphism
$\tau_{X}$ in the non-abelian Hodge theorem (see 
\ref{ssec-abelian-non-abelian-hodge}{\bf (NA(iii))}) is incompatible
with both the Gauss-Manin connection and the $\cstar$-action.

To circumvent such difficulties we use (see Section~\ref{ssec-euler})
the fact that the $\cstar$-fixed
points can be detected infinitesimally. This reduces the question to a
cohomological calculation for ordinary local systems of vector spaces
which can then be carried out by using Simpson's higher K\"{a}hler
identities \cite[Lemma~2.2]{simpson-higgs}. 

It is instructive to first analyze the case when we are dealing with
sections that pass trough smooth points of the moduli stack. This is the
subject of Section~\ref{ssec-smooth}. Finally in Section~\ref{ssec-general} we
complete the proof of Theorem~\ref{thm-main} by reducing the general
case to the smooth situation.

\subsection{The Euler vector field} \label{ssec-euler}

Let $M$ be a scheme over $\op{Spec}(\bbc)$ which is equipped with 
an algebraic $\cstar$-action $\mu : \cstar \times M \to M$. Since the
group $\cstar = \op{Spec}(\bbc[t,t^{-1}])$ is parallelizable the action
$\mu$ induces a canonical Euler vector field $\eta \in H^{0}(M,T_{M})$
defined as $\eta := i^{*}(d\mu)(d/dt,0)$, with $i : M \hookrightarrow 
\cstar \times M$, $i(x) := (1,x)$ being the natural inclusion. From the
definition of $\eta$ it is clear that a smooth closed point $x \in M$ is
$\cstar$-fixed if and only if the image $\eta_{x}$ of $\eta$ under the
natural evaluation map $\op{ev}_{x} : T_{M} \to T_{M,x} := 
T_{M}\otimes k(x)$ is zero.

More generally consider a smooth Artin stack $\mycal{X}$. 
Assume further that $\mycal{X}$ is equipped
with a $\cstar$-action $\mu : \cstar \times \mycal{X} \to \mycal{X}$.
Then we can again use the flow to construct 
the Euler vector field for $\mu$ which in
this case will be a 1-morphism of stacks $\eta : \mycal{X} \to
T_{\mycal{X}}$ which is a section of the structure morphism
$\pr_{\mycal{X}} : T_{\mycal{X}} \to \mycal{X}$. 
The $\cstar$-fixed sub stack $\mycal{X}^{\cstar}$ of $\mycal{X}$ is 
defined as the stacky zero
locus of the vector field $\eta$. More precisely $\mycal{X}^{\cstar}
:= \mycal{X}\times_{\ver_{\mycal{X}},T_{\mycal{X}},\eta} \mycal{X}$
where as usual $\ver_{\mycal{X}}$ is the vertex of the cone stack
$T_{\mycal{X}}$ (see \ref{appss-truncated-cotangent} for notation).

\begin{rem} \label{rem-quotient-fixed-points}
Note that if $\mycal{X}$ is of the form $[R/G]$ with $R$ a
quasi projective scheme over $\bbc$ and $G$ a complex reductive group and
if in addition the $\cstar$-action on $\mycal{X}$ descends from an
action of $G\times \cstar$ on $R$, then the sub stack
$\mycal{X}^{\cstar}$ is just the quotient stack $[R^{\,\cstar}/G]$.
\end{rem}

Let now $S$ be a smooth complex quasi-projective variety and let 
$f : X \to S$ be a smooth projective morphism. Consider the stack 
$\pdr : \stackdr{X/S}{n} \to S$ of relative rank $n$ local systems
and the stack $\phod : \stackhod{X/S}{n} \to S\times {\mathbb A}^{1}$ 
of $\lambda$-connections of rank $n$
(see \cite[Proposition~4.1]{simpson-santa-cruz}
for existence).  Let $a_{\op{DR}} : 
S \to \stackdrreg{X/S}{n}$ be an algebraic section which
is horizontal with respect to the non-abelian Gauss-Manin connection.

For future reference we formulate the following lemma which is well
known to the experts:

\begin{lem} \label{lem-dol-section} Assume $S$ is projective.
\begin{description}
\item[(i)] There exists a canonical extension of $a_{\op{DR}}$ to a section 
$a_{\op{Hod}} : S\times {\mathbb A}^{1} 
\to \stackhodreg{X/S}{n}$. In particular there is an
algebraic section $a_{\op{Dol}} : S \to \stackdolreg{X/S}{n}$ which
corresponds to $a_{\op{DR}}$ through the non-abelian Hodge theorem.
\item[(ii)] For a point $s \in S$ the local system 
$a_{\op{DR}}(s) \in \stackdr{X_{s}}{n}$ will underly a $\bbc$VHS iff 
$a_{\op{Dol}}(s) \in \stackdol{X_{s}}{n}^{\cstar}$, i.e. iff the Euler
vector field $\eta$ vanishes at $a_{\op{Dol}}(s)$.
\end{description}
\end{lem}
{\bf Proof.} Let $(F,\nabla_{f})$ be a bundle with a relative integrable
connection on $X$ representing the section $a_{DR}$. The horizontality
of $a_{DR}$ with respect to the non-abelian Gauss-Manin connection
implies (see Proposition~\ref{prop-integrable-section}(1)) the
existence of a global integrable connection $\nabla : F \to
F\otimes\Omega^{1}_{X}$ which induces $\nabla_{f}$. Since $S$ (and
therefore $X$) is projective we can invoke
\cite[Corollary~3.10]{simpson-higgs} and conclude that
$(F,\nabla)$ will correspond to a
global Higgs bundle $(E,\theta)$ on $X$ which is filtered by Higgs
subbundles and for which the associated graded is a direct sum of
stable Higgs bundles with vanishing rational Chern classes. Using
$(F,\nabla)$ and $(E,\theta)$ one can produce a twistor line
${\mathbb A}^{1} \to \stackhod{X}{n}$ in exactly the
same way as in \cite[Section~4]{simpson-santa-cruz}. Namely if
$\delbar_{F}$ and $\delbar_{E}$ denote the complex structure operators
for $F$ and $E$ respectively we define $\nabla_{\lambda} :=
\lambda\nabla + (1-\lambda)\theta$ and $\delbar_{\lambda} :=
\lambda\delbar_{F} + (1 - \lambda)\delbar_{E}$. The section 
${\mathbb A}^{1} \to \stackhod{X}{n}$, $\lambda \mapsto
(\delbar_{\lambda},\nabla_{\lambda})$ gives then the twistor line in
question. Consider the restriction $a : S \to \stackhod{X/S}{n}$ 
of the family of global $\lambda$-connections
$\delbar_{\lambda} + \nabla_{\lambda}$ to a family of relative
$\lambda$-connections along  the fibers of $f$. Therefore part (i) of
the lemma will be proven if we can show that $a$ will factor as
\[
S \stackrel{a_{\op{Hod}}}{\to} \stackhodreg{X/S}{n} \hookrightarrow 
\stackhod{X/S}{n}
\] 
where as usual $\stackhodreg{X/S}{n} \subset \stackhod{X/S}{n}$ denotes
the part of $\stackhod{X/S}{n}$ where $\phod$ is smooth. This follows
from the important fact (proven in the
subsection "Griffiths transversality revisited" of
\cite[Section~11]{simpson-santa-cruz}) that the map
$\pdr : \stackhodreg{X/S}{n} \to S$ is trivial 
\'{e}tale locally on $\stackhodreg{X/S}{n}$. Alternatively one can use
the fact that $a$ comes from a section  ${\mathbb A}^{1} \to
\stackhod{X}{n}$ and \cite[Theorem~9.1]{simpson-santa-cruz}.

\medskip

\noindent
The statement in part (ii) is just \cite[Lemma~4.1]{simpson-higgs} 
combined with
Remark~\ref{rem-quotient-fixed-points} and
\cite[Lemma~3.5]{simpson-higgs}. 

\ \hfill $\Box$

\bigskip

The previous lemma shows that in order to understand the relation
between the non-abelian $(p,p)$ classes and the non-abelian Gauss-Manin
connection we have to study the way $\eta$ interacts with the
non-abelian Gauss-Manin connection. For this we will need a description
of $\eta$ similar to the concrete description of $T_{\pdol}$ obtained in
Lemma~\ref{lem-vertical-tangent-spaces}.

Let $\bun{X/S}{n}$ denote the moduli stack of rank
$n$ vector bundles on $X$ with vanishing rational Chern classes along
the fibers of $f$. Let $\pi : \bun{X/S}{n} \to S$ be the structure
map and let $\bunreg{X/S}{n}$ denote the open sub stack in
$\bun{X/S}{n}$ over which the morphism $\pi$ is smooth.
\glosstex(moduli)[p]{bundles} 
\glosstex(moduli)[p]{bundles-regular}
The forgetting the Higgs field induces a
morphism of stacks
$p : \stackdol{X/S}{n} \to \bun{X/S}{n}$ which, as explained in
Section~\ref{ssec-abelian-non-abelian-hodge}{\bf (NA(iv))} should be 
though of as the projection to
the $(1,0)$ part of the non-abelian Dolbeault cohomology. By
definition the action of $\cstar$ on $\stackdol{X/S}{n}$ respects $p$
and so $\eta$ can be interpreted as a section of the relative tangent
stack $T_{p} := T_{\stackdol{X/S}{n}/\bun{X/S}{n}} \to
\stackdol{X/S}{n}$.

According to Lemma~\ref{lem-vertical-tangent-spaces}(ii) we have 
\[
T_{\pdol} = \hh(\tau_{\leq
1}\bbr{\fdol}_{*}\lieg{\fdol}{E_{\op{un}},\theta_{\op{un}}}),
\]
and so given an $S$-point $a_{\op{Dol}} : S \to \stackdolreg{X/S}{n}$
represented by a relative Higgs bundle $(E,\theta_{f}) \to X$ we
have $a_{\op{Dol}}^{*}T_{\pdol} = \hh(\tau_{\leq 1}\bbr
f_{*}\lieg{f}{E,\theta_{f}})$. The explicit description of $\eta$ we
need is now given by the following

\begin{lem} \label{lem-euler-field} Let $(E,\theta_{f})$ be a Higgs
bundle on $X/S$ representing a section $a_{\op{Dol}} : S \to
\stackdolreg{X/S}{n}$. Then 
\begin{description}
\item[(i)] There is an isomorphism of Picard stacks over $S$
\[
a_{\op{Dol}}^{*}T_{p} = \hh(\tau_{\leq 1}\op{Cone}[\bbr f_{*}\sigma_{\leq
0}\lieg{f}{E,\theta_{f}}[-1] \to \bbr f_{*}\sigma_{\geq
1}\lieg{f}{E,\theta_{f}}]),
\]
where $\sigma$ denotes the stupid truncation.
\item[(ii)] The section $a_{\op{Dol}}^{*}\eta : S \to
a_{\op{Dol}}^{*}T_{p} \hookrightarrow a_{\op{Dol}}^{*}T_{\pdol}$
corresponds to the natural morphism of complexes on $X$:
\[
{\mathcal O}_{X}[-1] \stackrel{\theta_{f}}{\to} \sigma_{\geq
1}\lieg{f}{E,\theta_{f}} 
\]
\end{description}
\end{lem}
{\bf Proof.} By definition $E_{\op{un}}$ is the pull-back from the
universal bundle on $\bun{X/S}{n}\times_{S} X$ via the projection
$p\times \op{id}_{X} : \stackdol{X/S}{n}\times_{S} X \to
\bun{X/S}{n}\times_{S} X$. Furthermore the \v{C}ech calculation in the
proof of Lemma~\ref{lem-vertical-tangent-spaces} shows that the
pullback of the 
differential morphism $dp : T_{\pdol} \to p^{*}T_{\pi}$ via
$a_{\op{Dol}}^{*}$ is induced
from the natural morphism of complexes 
\settowidth{\dgrmone}{$\op{End}(E)\otimes \Omega^{1}_{f}$ \ }
\settowidth{\dgrmtwo}{$\op{End}(E)$ \ }
\begin{equation} \label{differential}
\lieg{f}{E,\theta_{f}} = 
\left[
\begin{minipage}[c]{\dgrmone}
$\xymatrix{
\op{End}(E) \ar[d]_{\op{ad}_{\theta_{f}}} \\
\op{End}(E)\otimes \Omega^{1}_{f} \ar[d]_{\op{ad}_{\theta_{f}}} \\
\op{End}(E)\otimes \Omega^{2}_{f} \ar[d]_{\op{ad}_{\theta_{f}}} \\
\text{\vdots}
}$
\end{minipage}
\; \right] \longrightarrow
\left[
\begin{minipage}[c]{\dgrmtwo}
$\xymatrix{
\op{End}(E) \ar[d] \\
0 \ar[d] \\
0 \ar[d] \\
\text{\vdots}
}$
\end{minipage}
\; \right] = 
\sigma_{\leq 0}\lieg{f}{E,\theta_{f}}
\end{equation}
On the other hand the morphism (\ref{differential}) fits in a short
exact sequence of complexes
\[
0 \to \sigma_{\geq 1}\lieg{f}{E,\theta_{f}} \to \lieg{f}{E,\theta_{f}}
\to \sigma_{\leq 0}\lieg{f}{E,\theta_{f}} \to 0
\]
which pushes forward to a distinguished triangle in $D^{b}(S)$:
\begin{equation} \label{differential-downstairs}
\bbr f_{*}\sigma_{\leq 0}\lieg{f}{E,\theta_{f}}[-1] \to \bbr
f_{*}\sigma_{\geq 1}\lieg{f}{E,\theta_{f}} \to \bbr
f_{*}\lieg{f}{E,\theta_{f}} \to \bbr f_{*}\sigma_{\leq
0}\lieg{f}{E,\theta_{f}}.
\end{equation}
The pullback of $T_{p}$ via the section $a_{\op{Dol}}$ is naturally
identified with the fiber product 
\[
a_{\op{Dol}}^{*}T_{\pdol}\times_{a_{\op{Dol}}^{*}(dp), (p\circ
a_{\op{Dol}})^{*}T_{\pi}, \ver_{(p\circ
a_{\op{Dol}})^{*}T_{\pi}}} S,
\]
which due to (\ref{differential-downstairs}) and the functoriality
of $\hh$-construction coincides with 
\[
\hh(\tau_{\leq 1}\op{Cone}[\bbr
f_{*}\sigma_{\leq 
0}\lieg{f}{E,\theta_{f}}[-1] \to \bbr f_{*}\sigma_{\geq
1}\lieg{f}{E,\theta_{f}}]).
\] 
This proves part (i) of the lemma. 

\medskip

The proof of part (ii) is almost tautological. Indeed the group
$\cstar$ acts on $\stackdol{X/S}{n}$ by
rescaling the Higgs fields and so the flow vector field
$a_{\op{Dol}}^{*}\eta$  
corresponds to the global section $\theta_{f} \in
H^{0}(X,\op{End}(E)\otimes \Omega^{1}_{f})$ which is naturally
interpreted as a morphism $\theta_{f} : {\mathcal O}_{X}[-1] \to
\sigma_{\geq 1}\lieg{X/S}{E,\theta_{f}}$. Therefore in terms of the
cohomological description (i) of $a_{\op{Dol}}^{*}T_{p}$ the section
$a_{\op{Dol}}^{*}\eta : S \to a_{\op{Dol}}^{*}T_{p}$ corresponds to
the $\tau_{\leq 1}$ truncation of the morphism 
\[
Rf_{*}{\mathcal
O}_{X}[-1] \to \bbr f_{*}\sigma_{\geq
1}\lieg{f}{E,\theta_{f}} \to 
\op{Cone}[\bbr
f_{*}\sigma_{\leq 
0}\lieg{f}{E,\theta_{f}}[-1] \to \bbr f_{*}\sigma_{\geq
1}\lieg{f}{E,\theta_{f}}].
\]
This finishes the proof of the lemma. \hfill $\Box$

\bigskip
\bigskip

\noindent
The $\tau_{\leq 1}$ truncation of the push forward of the
morphism $\theta_{f} : {\mathcal O}_{X}[-1] \to \sigma_{\geq
1}\lieg{f}{E,\theta_{f}}$ induces naturally a global section
$R^{1}f_{*}\theta_{f}$ of the
coherent sheaf $\bbr^{1}f_{*}\lieg{f}{E,\theta_{f}}$ on $S$. Let
\[
\epsilon_{f} \in H^{0}(S,
\bbr^{1}f_{*}\sigma_{\geq 1}\lieg{f}{E,\theta_{f}}/f_{*}\op{End}(E))
\]
denote the image of $R^{1}f_{*}\theta_{f}$. The 
previous lemma has the following simple but important corollary.

\begin{cor} \label{cor-fixed-point}
Let $f : X \to S$ be a smooth projective family over a smooth
quasi-projective scheme $S$. Let $(E,\theta_{f})$ be a relative Higgs
bundle on $X$ representing a section $a_{\op{Dol}} : S \to
\stackdolreg{X/S}{n}$. Then for point $s \in S$ we have
$a_{\op{Dol}}(s) \in \stackdol{X/S}{n}^{\cstar}$ iff 
$\epsilon_{f}(s) = 0$. 
\end{cor}
{\bf Proof.} By construction $\mycal{H}^{0}\bbr f_{*}\sigma_{\geq
1}\lieg{f}{E,\theta_{f}}] = 0$ and hence we have a quasi-isomorphism
\[
\begin{split}
\tau_{\leq 1} \op{Cone}[\bbr
f_{*}\sigma_{\leq 
0}\lieg{f}{E,\theta_{f}}[-1] \to \bbr f_{*}\sigma_{\geq
1}\lieg{f}{E,\theta_{f}}] & \cong (\bbr^{1} f_{*}\sigma_{\geq
1}\lieg{f}{E,\theta_{f}}/\bbr^{0}f_{*}\lieg{f}{E,\theta_{f}})[-1] \\
& = (\bbr^{1} f_{*}\sigma_{\geq
1}\lieg{f}{E,\theta_{f}}/f_{*}\op{End}(E))[-1].
\end{split}
\]
Combined with Lemma~\ref{lem-euler-field}(ii) this proves the
corollary. \hfill $\Box$

\
\bigskip
\bigskip

\begin{rem} \label{rem-euler-gerbe} One can recast the statement of
Corollary~\ref{cor-fixed-point} in more geometric terms. Name\-ly, the
cohomological description of Lemma~\ref{lem-euler-field}(ii) leads to
the following concrete description of 
the residual gerbe of the $S$-point 
$a_{\op{Dol}}^{*}\eta$ of $a_{\op{Dol}}^{*}T_{\pdol}$.

Since $a_{\op{Dol}}^{*}\eta$ is a
section of the structure morphism $a_{\op{Dol}}^{*}T_{\pdol} \to S$
it follows that the coarse
moduli space of the $S$-stack $a_{\op{Dol}}^{*}\eta(S)$ (taken with the
reduced sub stack structure) is isomorphic to $S$. In particular to
specify $a_{\op{Dol}}^{*}\eta$ one only needs to specify its residual
gerbe $\mycal{G} \to S$ \cite[Section~5]{laumon-stacks}. The proof of
Corollary~\ref{cor-fixed-point} implies that $a_{\op{Dol}}^{*}T_{p}$
is isomorphic as a Picard stack to $\hh(f_{*}\op{End}(E) \to
\bbr^{1}f_{*}\sigma_{\geq 1}\lieg{f}{E,\theta_{f}})$. Furthermore
we have  
\[
\begin{split}
\bbr^{1}f_{*}\sigma_{\geq 1}\lieg{f}{E,\theta_{f}}) & = 
\bbr^{0}f_{*}(\sigma_{\geq 1}\lieg{f}{E,\theta_{f}})[-1]) \\
& = 
\bbr^{0}f_{*}[\xymatrix@1{
\op{End}(E)\otimes \Omega^{1}_{f} \ar[r]^-{\op{ad}_{\theta_{f}}} & 
\op{End}(E)\otimes \Omega^{2}_{f} \ar[r]^-{\op{ad}_{\theta_{f}}} &
\ldots
}
].
\end{split}
\]
By examining the long exact sequence of sheaves corresponding to the
distinguished triangle (\ref{differential-downstairs}) one sees that 
the morphism $f_{*}\op{End}(E) \to \bbr^{0}f_{*}(\sigma_{\geq
1}\lieg{f}{E,\theta_{f}})[-1])$ is just the push forward of 
the morphism of complexes 
\[
\xymatrix{
\op{End}(E) \ar[r]^-{\op{ad}_{\theta_{f}}} \ar[d]
& \op{End}(E)\otimes \Omega^{1}_{f} \ar[d]^{\op{ad}_{\theta_{f}}} \\
0 \ar[r] \ar[d] & \op{End}(E)\otimes \Omega^{2}_{f} 
\ar[d]^{\op{ad}_{\theta_{f}}} \\ 
\text{\vdots} & \text{\vdots}
}
\]
Note moreover that since $\theta_{f}
\in H^{0}(X, \op{End}(E)\otimes \Omega^{1}_{f})$ is in the kernel of
$\op{ad}_{\theta_{f}}$, the push forward of $\theta_{f}$ can be viewed
as a global section of the sheaf $\bbr^{0}f_{*}(\sigma_{\geq
1}\lieg{f}{E,\theta_{f}})[-1])$ on $S$.

We will abuse notation and write $\op{ad}_{\theta_{f}} :
f_{*}\op{End}(E) \to \bbr^{0}f_{*}(\sigma_{\geq
1}\lieg{f}{E,\theta_{f}})[-1])$ and $\theta_{f} \in H^{0}(S,
\bbr^{0}f_{*}(\sigma_{\geq 1}\lieg{f}{E,\theta_{f}})[-1]))$.  

The sections of $a_{\op{Dol}}^{*}T_{p}$ over an $S$-scheme $T \to S$
are pairs $(A,\alpha)$ where $A \to T$ is a $f_{T*}\op{End}(E_{T})$ torsor
and $\alpha  : A \to \bbr^{0}f_{T*}(\sigma_{\geq
1}\lieg{f_{T}}{E_{T},\theta_{f_{T}}})[-1]))$ is
$\op{ad}_{\theta_{f_{T}}}$-equivariant map. By
Lemma~\ref{lem-euler-field}(ii) it follows that $a_{\op{Dol}}^{*}\eta$
is represented by the pair $(f_{*}\op{End}(E), \theta_{f} +
\op{ad}_{\theta_{f}}(\bullet))$. In particular the residual gerbe
$\mycal{G}$ 
of $a_{\op{Dol}}^{*}\eta$ is just the neutral gerbe
$B(\ker(\op{ad}_{\theta_{f}}))$ on $S$.

It is instructive to compare the statement of Corollary~\ref{cor-fixed-point}
with this geometric description. Given a point $s \in S$ put $E_{s} :=
E_{|X_{s}}$ and $\theta_{s} := \theta_{f|X_{s}}$. Let as before $B$
denote the nilpotent scheme
$\op{Spec}(\bbc[t]/(t^2))$ and let $p_{X_{s}} :
X_{s}\times B \to X_{s}$ and $p_{B} : X_{s}\times B \to B$ denote the
natural projections. Identify  $B$ with the first infinitesimal
neighborhood of $1 \in \cstar$, i.e. put $t = \lambda -1$ where
$\lambda$ is a standard coordinate on $\cstar$, i.e. 
$\cstar = \op{Spec}(\bbc[\lambda,\lambda^{-1}])$.
The point $\epsilon_{f}(s) =
(E_{s},\theta_{s}) \in \stackdol{X_{s}}{n}$ will be fixed under the
$\cstar$-action iff there is an automorphism $\psi : p_{X_{s}}^{*}E
\widetilde{\to} p_{X_{s}}^{*}E$ satisfying: (a) $\psi_{|X_{s}\times
\op{Spec}(\bbc)} = \op{id}_{E_{s}}$; and (b) $\psi\circ
p_{X_{s}}^{*}\theta_{s} = (\lambda p_{X_{s}}^{*}\theta_{s})\circ
\psi$, where by abuse of notation $\lambda$ denotes the 1-jet of the
function $\lambda$ at $1 \in \cstar$.

Now due to (a) and (b) we can write $\psi = \op{id} + t\varphi$ with
$\varphi \in \op{End}(E_{s})$ and $\lambda
p_{X_{s}}^{*}\theta_{s} = \theta_{s} + t\theta_{s}$ and hence
$\epsilon_{f}(s) \in \stackdol{X_{s}}{n}^{\cstar}$ iff $(\op{id} +
t\varphi)\circ \theta_{s} = (\theta_{s} + t\theta_{s})\circ (\op{id} +
t\varphi)$, i.e. if and only if  $\theta_{s} +
\op{ad}_{\theta_{s}}(\phi) = 0$. This of course is exactly the
statement of Corollary~\ref{cor-fixed-point}
\end{rem}

\

\bigskip

We are now ready to prove Theorem~\ref{thm-main}. With the hope of
making the exposition more accessible we first treat the case when $f
: X \to S$ is of relative dimension one and the section $a_{\op{DR}} :
S \to \stackdr{X/S}{n}$ passes only trough smooth points of moduli,
i.e. when $a_{\op{DR}} : S \to \stackdrreg{X/S}{n} \subset
\stackdr{X/S}{n}$.

\subsection{A warmup - the smooth case} \label{ssec-smooth}

For the duration of this section assume that $S$ is a smooth
quasi-projective variety and that $f : X \to S$ is a smooth family of
integral curves over $S$. Let $a_{\op{DR}} : S \to
\stackdrreg{X/S}{n}$ be an algebraic section which is horizontal
with respect to the non-abelian Gauss-Manin connection. Let
$(F,\nablaf)$ be a relative connection on $X$ representing
$a_{\op{DR}}$. Proposition~\ref{prop-integrable-section} implies that
there exists a global integrable connection $(F,\nabla )$ which
induces $\nablaf$ along the fibers of $f$. By replacing if necessary
$(F,\nabla)$ by its semi simplification
(cf. \cite[Lemma~3.5]{simpson-higgs}) we may assume without a loss of
generality that $(F,\nabla)$ is a reductive local system on $X$. 
Therefore the statement of
Theorem~\ref{thm-main} becomes equivalent to the following

\begin{prop} \label{prop-smooth-curves} Let $f : X\to S$ be a smooth
family of integral curves over a projective base $S$. Let $(F,\nabla)$
be a reductive local system on $X$ and assume that there exists a
point $o \in S$ so that the induced local system $(F_{o},\nabla_{o})$
on $X_{o}$ underlies a $\bbc$VHS and such that $(F_{o},\nabla_{o})$ is
a smooth point of $\stackdr{X_{o}}{n}$. Then for any $s \in S$ the local
system $(F_{s},\theta_{s})$ underlies a $\bbc$VHS on $X_{s}$.
\end{prop}
{\bf Proof.} Let $s \in S$ be any point. By the Lefschetz hyperplane
section theorem we can always find a smooth curve $C \subset S$ cut
out by hyperplanes of sufficiently high degree so that $C$ contains
the two points $o, s \in S$. Since the horizontality with respect to
the non-abelian Gauss-Manin connection as well as the $\cstar$-action
on the Dolbeault moduli spaces are stable under base change we may
assume without losing generality that $S$ is a smooth projective curve.

Let $(E,\theta)$ be the Higgs bundle on $X$ corresponding
to $(F,\nabla)$ and let $(E,\theta_{f})$ be the induced relative Higgs
bundle. Consider the section 
\[
\epsilon_{f} \in H^{0}(S,\bbr^{1}f_{*}\sigma_{\geq
1}\lieg{f}{E,\theta_{f}}/f_{*}\op{End}(E))
\]
defined after the proof of Lemma~\ref{lem-euler-field}. From the long
exact sequence of sheaves corresponding to the distinguished triangle
(\ref{differential-downstairs}) it follows that $\bbr^{1}f_{*}\sigma_{\geq
1}\lieg{f}{E,\theta_{f}}/f_{*}\op{End}(E) \subset
\bbr^{1}f_{*}\lieg{f}{E,\theta_{f}}$ and so $\epsilon_{f}$ can be thought of
as an element in
$H^{0}(S,\bbr^{1}f_{*}\lieg{f}{E,\theta_{f}})$. 

Since $(F,\nabla)$ represents a horizontal section of the local system
of stacks $\stackdr{X/S}{n} \to S$ which passes trough a smooth point 
of $\stackdr{X_{o}}{n}$ it follows that $(F_{s},\nabla_{s})$ will be a
smooth point of $\stackdr{X_{s}}{n}$ for all $s$.
Hence, according to Corollary~\ref{cor-fixed-point} the proposition will be
proven if we 
show that $\epsilon_{f} = 0$ in 
$H^{0}(S,\bbr^{1}f_{*}\lieg{f}{E,\theta_{f}})$. 
By hypothesis $(E_{o},\theta_{o})$ is $\cstar$-fixed and so by
Corollary~\ref{cor-fixed-point} we have
$\epsilon_{f}(o) = 0$. On the other hand the sheaf
$\bbr^{1}f_{*}\lieg{f}{E,\theta_{f}}$ can be naturally extended to a
Higgs sheaf on $S$. Indeed consider first direct image
\cite[Section~4]{simpson-some-families}  of 
the Higgs bundle $(\op{End}(E),\op{ad}_{\theta})$ 
under the projective map $f : X \to
S$. To describe this direct image note that since $f$ is a fibration
of curves over a curve we have 
\[
\begin{split}
\lieg{f}{E,\theta_{f}} & = 
[\xymatrix@1{
\op{End}(E) \ar[r]^-{\op{ad}_{\theta_{f}}} & \op{End}(E)\otimes
\Omega^{1}_{f}}] \\
\lieg{X}{E,\theta} & = 
[\xymatrix@1{
\op{End}(E) \ar[r]^-{\op{ad}_{\theta}} & \op{End}(E)\otimes
\Omega^{1}_{X} \ar[r]^-{\op{ad}_{\theta}} & \op{End}(E)\otimes
\Omega^{2}_{X}}] 
\end{split}
\]
and hence we have a short exact sequence of complexes on $X$
\begin{equation} \label{higgs-pushforward}
0 \to \lieg{f}{E,\theta_{f}}[-1]\otimes f^{*}\Omega^{1}_{S} \to 
\lieg{X}{E,\theta} \to \lieg{f}{E,\theta_{f}} \to 0
\end{equation}
As explained in \cite[Section~4]{simpson-some-families} the first
direct image of $(\op{End}(E),\op{ad}_{\theta})$ under $f$ is just the
sheaf $\bbr^{1}f_{*}\lieg{f}{E,\theta_{f}}$ together with the first
edge homomorphism $\delta : \bbr^{1}f_{*}\lieg{f}{E,\theta_{f}} \to
\bbr^{1}f_{*}\lieg{f}{E,\theta_{f}}\otimes \Omega^{1}_{S}$ for the
long exact sequence of hyper-derived images corresponding to
(\ref{higgs-pushforward}). 

Recall next that $\epsilon_{f}$ was defined as the image of
$R^{1}f_{*}\theta_{f}$ in $\bbr^{1}f_{*}\lieg{f}{E,\theta_{f}}$. On
the other hand we have a natural element $\epsilon \in
H^{0}(S,\bbr^{1}f_{*}\lieg{X}{E,\theta})$ defined analogously as
the image of $R^{1}f_{*}\theta$. Since $\theta$ restricts to the
relative Higgs field $\theta_{f}$ it follows that $\epsilon$ maps to
$\epsilon_{f}$ under the natural map $\bbr^{1}f_{*}\lieg{X}{E,\theta}
\to \bbr^{1}f_{*}\lieg{f}{E,\theta_{f}}$. Therefore from the long 
exact sequence of hyper-derived images corresponding to
(\ref{higgs-pushforward}) we conclude that $\delta(\epsilon_{f}) = 0$,
i.e. $\epsilon_{f}$ is in fact an element in the $0$-th cohomology
group $H^{0}(S,(\bbr^{1}f_{*}\lieg{X}{E,\theta},\delta))$ of the Higgs
sheaf $(\bbr^{1}f_{*}\lieg{X}{E,\theta},\delta)$. 

To finish the proof of the proposition it remains only to note that by
\cite[Corollary~5.2]{simpson-some-families} the Higgs bundle
$(\bbr^{1}f_{*}\lieg{X}{E,\theta},\delta)$ on $S$ corresponds to the
local system $(\bbr^{1}f_{*}\liea{f}{F,\nablaf},D)$ with $D :
\bbr^{1}f_{*}\liea{f}{F,\nablaf} \to 
\bbr^{1}f_{*}\liea{f}{F,\nablaf}\otimes \Omega^{1}_{S}$ being just the
first edge homomorphism in the long 
exact sequence of hyper-derived images corresponding to
\begin{equation} \label{localsystem-pushforward}
0 \to \liea{f}{F,\nabla_{f}}[-1]\otimes \Omega^{1}_{S} \to
\liea{X}{F,\nabla} \to \liea{f}{F,\nabla_{f}} \to 0.
\end{equation}
Combined with Simpson formality result \cite[Lemma~1.2]{simpson-higgs}
this gives 
\[
\epsilon \in
H^{0}(S,(\bbr^{1}f_{*}\lieg{X}{E,\theta},\delta)) = 
H^{0}(S,(\bbr^{1}f_{*}\liea{f}{F,\nablaf},D)).
\] 
In particular
$\epsilon_{f}$ can be interpreted as an algebraic section of the sheaf 
$\bbr^{1}f_{*}\liea{f}{F,\nablaf}$ which is horizontal with respect to
$D$. Hence if $\epsilon$ vanishes at a point it vanishes everywhere
which proves the proposition. \hfill $\Box$

\

\bigskip

\begin{rem} \label{rem-no-curves} (i) The proof of the proposition works
in essentially the same way without reducing to the case when $S$ is a
curve. Indeed to make the above argument work in general one only
needs to replace (\ref{higgs-pushforward})  (respectively
(\ref{localsystem-pushforward})) by the short exact sequence
of complexes 
\[
0 \to \lieg{f}{E,\theta_{f}}[-1]\otimes f^{*}\Omega^{1}_{S} \to 
\lieg{X}{E,\theta}/I^{2}\lieg{X}{E,\theta} \to \lieg{f}{E,\theta_{f}}
\to 0,
\]
where as usual 
\[
\begin{split}
I^{1}\lieg{X}{E,\theta} & :=
\op{im}[\lieg{X}{E,\theta}\otimes f^{*}\Omega_{S} \to
\lieg{X}{E,\theta}] \\
I^{2}\lieg{X}{E,\theta} & :=
\op{im}[I^{1}\lieg{X}{E,\theta}\otimes f^{*}\Omega_{S} \to
\lieg{X}{E,\theta}]
\end{split}
\]
(see \cite[Section~4]{simpson-some-families}). 

\medskip

\noindent
(ii) The hypothesis that $f$ is of relative dimension one is
superfluous in the statement of
Proposition~\ref{prop-smooth-curves} and was only used to simplify the
exposition. The proof for the case of
general fiber dimension works in exactly the same way. The reason we
chose to state Proposition~\ref{prop-smooth-curves} for for fiber
dimension one only will become clear in the next section where we
show that the smoothness of the morphisms of $\pdr$ and $\pdol$ is
easier to control when one is dealing with curves.
\end{rem}

\

\bigskip

\noindent
The proof of Proposition~\ref{prop-smooth-curves} 
has the following immediate corollary

\begin{cor} \label{cor-quasi-projective} Let $f : X \to S$ be
 a smooth
projective morphism with connected fibers to a smooth quasi-projective
$S$. Assume that $f$ extends to a projective morphism $\bar{f} :
\overline{X} \to \overline{S}$ where $\overline{X}$ and $\overline{S}$
are smooth and projective and $\overline{S}\setminus S$ and 
$\overline{X}\setminus X$ are divisors with strict normal
crossings. Let $(F,\nabla)$ be a reductive local system on $\overline{X}$ 
and assume that there exists a
point $o \in S$ so that the induced local system $(F_{o},\nabla_{o})$
on $X_{o}$ underlies a $\bbc$VHS. Then for any $s \in S$ the local
system $(F_{s},\theta_{s})$ underlies a $\bbc$VHS on $X_{s}$.
\end{cor}
{\bf Proof.}  The proof is essentially the same but at the last stage
we need to invoke a much stronger theorem of Simpson
\cite[Corollary~5.12]{simpson-some-families} which asserts that the
filtered regular Higgs bundle (= parabolic Higgs bundle)
$(\bbr^{1}f_{*}\lieg{X}{E,\theta},\delta)$ on $S$ 
corresponds to the local system
$(\bbr^{1}f_{*}\liea{f}{F,\nablaf},D)$.
\hfill $\Box$

\

\bigskip

\begin{rem} (i) Note that the reduction to the situation where $S$ is
a curve is really essential here since the theory of regular filtered
Higgs bundles \cite{simpson-non-compact-curves} is developed only in
the curve case.

\medskip

\noindent
(ii) It is tempting to try and extend the previous corollary to a
situation where $(F,\nabla)$ is a reductive local system on 
$\overline{X}$ with logarithmic poles along $\overline{X}\setminus
X$. 

There are several obstructions to carrying out the arguments in this
case. First of all it is unclear whether such a local system will
have a harmonic metric. There are however some existence results
under mild additional assumptions. For example
\cite[Theorem~11.4]{biquard} guarantees the existence of a harmonic
metric for $(F,\nabla)$ over $X$ as long as $(F,\nabla)$ is stable
as a parabolic local system and the residual primitive homogeneous
local systems on the the total space of the normal bundle of
$\overline{X}\setminus X$ are semi simple. Similarly in
\cite{corlette-open} the existence of a harmonic metric on
$(F,\nabla)$ is proven under the assumptions that $(F,\nabla)$ is
reductive and has a quasi-unipotent monodromy at infinity.
Ideally one would like to apply these results to a
logarithmic $(F,\nabla)$ representing a section $a_{\op{DR}}$ which is
horizontal for the non-abelian Gauss-Manin connection. 
For this one will have to analyze the behavior of $a_{\op{DR}}$ at
infinity which is an interesting question in its own right.

The second problem is in the lack of a push forward result similar to
\cite[Corollary~5.12]{simpson-some-families} which works for parabolic
local systems with poles along $\overline{X}\setminus X$. Specifically
one needs to find the right growth conditions for sections in a
harmonic bundle near a point in a normal crossings divisor.
\end{rem}

\subsection{A reduction - the general case} \label{ssec-general}

In this section we show how to reduce the statement
Theorem~\ref{thm-main} to the situation in
Proposition~\ref{prop-smooth-curves}. For this we have to explain how
to: (a) deal with
morphisms $f : X \to S$ of fiber dimension bigger than one; and (b)
how to deal  with
sections $a_{\op{DR}} : S \to \stackdr{X/S}{n}$ that do not
necessarily land in $\stackdrreg{X/S}{n}$.

The problem (a) is quite mild. In fact Remark~\ref{rem-no-curves}(ii) 
shows how (a) can be tackled directly within the
method of proof of Proposition~\ref{prop-smooth-curves}. Alternatively
we may use the Lefschetz hyperplane section theorem and its
non-abelian version - the Mehta-Ramanathan type restriction result
\cite[Proposition~3.6]{simpson-higgs} (see also
\cite[Theorem~7.2.1]{huybrechts-lehn}).  This second route is
preferable since it will put us into a favorable setup for dealing
with (b).

Concretely let $f : X \to S$ be a smooth
projective morphism with connected fibers. Assume that $S$ is
quasi-projective and let  $a_{\op{DR}} : S
\to \stackdr{X/S}{n}$ be an algebraic section which is horizontal with
respect to the non-abelian Gauss-Manin connection. Since the statement
of Theorem~\ref{thm-main} is stable under base changes $T \to S$ with
$T$ smooth and projective we may assume without a loss of generality
that $f : X \to S$ has a section $\xi : S to X$.

Let $i : C
\hookrightarrow X$ be a general enough intersection of relative
hyperplanes so that $f\circ i : C \to S$ is smooth with connected
fibers and $\dim_{\bbc}(C/S) = 1$. By the Lefschetz hyperplane section
theorem we have a surjection $\pi_{1}(C_{s}) \twoheadrightarrow
\pi_{1}(X_{s})$ for every $s \in S$ and therefore $i$ induces a closed
immersion of algebraic stacks $\stackb{X/S}{n} \subset
\stackb{C/S}{n}$. The strong form of the Riemann-Hilbert
correspondence \cite[Proposition~7.8]{simpson-moduli2} (we need the
existence of $\xi$ for that!)  yields then
a closed immersion of analytic stacks $\stackdr{X/S}{n}^{\op{an}}
\subset \stackdr{C/S}{n}^{\op{an}}$ which in turn implies that the
natural morphism $i^{*}_{\op{DR}}
\stackdr{X/S}{n} \to \stackdr{C/S}{n}$ of Artin
algebraic stacks is also a closed immersion. Furthermore due to the
functoriality of the algebraic construction of the non-abelian
Gauss-Manin connection the morphism $i^{*}_{\op{DR}}$ must be
horizontal and so $i^{*}_{\op{DR}}\circ a_{\op{DR}} : S \to
\stackdr{C/S}{n}$ is an algebraic section which is horizontal
with respect to the non-abelian Gauss-Manin connection as well. Next
observe that the surjectivity of $\pi_{1}(C_{s}) \to
\pi_{1}(X_{s})$ puts us in a position to apply
\cite[Corollary~4.3]{simpson-higgs} and conclude that $a_{\op{DR}}(s)$
underlies a $\bbc$VHS if and only if $i^{*}_{\op{DR}}\circ
a_{\op{DR}}(s)$ underlies a $\bbc$VHS.

The next fact we need is that due to \cite[Proposition~3.6]{simpson-higgs} 
the semistability of Higgs bundles  is preserved by restrictions to
hyperplane sections. Hence in the same way as in the proof of
Lemma~\ref{lem-dol-section} we can use
\cite[Corollary~3.10]{simpson-higgs} to reduce Theorem~\ref{thm-main}
to the following 

\begin{lem} \label{lem-not-smooth} Let $f : X \to S$ be a smooth
morphism with connected fibers of relative dimension one. Assume that
$S$ is quasi-projective and that $f$ admits a section $\xi$. Let
$a_{\op{DR}} : S \to \stackdr{X/S}{n}$ be a section represented by a
relative connection $(F,\nablaf)$. Then 
\begin{itemize}
\item[(a)] $a_{\op{DR}}$ is horizontal
iff there exists a global connection $\nabla : F \to
F\otimes\Omega^{1}_{X}$ which induces $\nablaf$.
\item[(b)] if $S$ is projective, then $(F_{s},\nabla_{s})$ underlies a
$\bbc$VHS iff $\epsilon_{f}(s) = 0$
\end{itemize}
\end{lem}
{\bf Proof.} If $a_{\op{DR}}$ happens to map $S$ into
$\stackdrreg{X/S}{n}$, then this is the content of
Proposition~\ref{prop-integrable-section}. 
For a general $a_{\op{DR}}$ we need to analyze the singularities of
the morphism $\pdr$. 

Let $f : X \to S$ be a smooth fibration of curves of genus $g > 1$.
Let $G$ be a complex
reductive group. It is well known (see e.g. \cite[Corollaries~4.5.2
and 8.1.9]{behrend-thesis}) that the stack $\bun{X/S}{G}$ is a
smooth stack over $S$ of relative dimension $(g-1)\dim G$. In
particular $\bunreg{X/S}{n} = \bun{X/S}{n}$. Moreover from the
definition of $\stackdol{X/S}{n}$ and the fact that $\dim (X/S) = 1$
it is clear that $\stackdol{X/S}{n}$ can be identified with the
relative cotangent stack $T^{\vee}_{\pi} \to \bun{X/S}{n}$ of $\pi :
\bun{X/S}{n} \to S$. Here by $T^{\vee}_{\pi}$ one means the vector
bundle stack
$\op{Spec}(S^{\bullet}R^{1}\op{pr}_{*}\op{End}(E_{\op{un}}))$ 
where as usual $E_{\op{un}} \to \bun{X/S}{n}\times_{S} X$ is the
universal bundle and $\op{pr} : \bun{X/S}{n}\times_{S} X \to
\bun{X/S}{n}$ is the natural projection.
 
This indicates that it is not unreasonable to expect that
$\stackdol{X/S}{n}$ (and hence $\stackdr{X/S}{n}$) will be close
enough to being smooth. In fact it is easy to see that the stack
$\stackdr{X/S}{n}$ (respectively $\stackdol{X/S}{n}$) embeds in a
stack which is smooth over $S$. 
To construct such an embedding one uses the following well known (see
for example \cite[Section~2.11]{beilinson-drinfeld-langlands}) 
rigidification
trick. 

Let $\stackdr{X/S(\log\xi)}{n}$ be the stack parameterizing
relative local systems on $f : X \to S$ with logarithmic poles along
$\xi(S)$. Let $(F_{\op{un}},\nabla_{\op{un}}) \to
\stackdr{X/S(\log\xi)}{n}\times_{S} X$ be the universal relative
local system and let $\mycal{F} := (\op{id}\times
\xi)^{*}F_{\op{un}}$. The residue of $\nabla_{\op{un}}$ along $\xi(S)$
is a section $\op{Res}(\nabla_{\op{un}}) \in H^{0}(\stackdr{X/S(\log\xi)}{n},
\op{End}(\mycal{F}))$ and $\stackdr{X/S}{n}$ is just the closed
substack of $\stackdr{X/S(\log\xi)}{n}$ cut out by the equation
$\op{Res}(\nabla_{\op{un}}) = 0$. Furthermore the maximal substack
$\stackdrreg{X/S(\log \xi)}{n} \subset \stackdr{X/S(\log\xi)}{n}$
which is smooth over $S$ can be described explicitly in this case.
Indeed, note first that the same argument as in the proof of
Lemma~\ref{lem-vertical-tangent-spaces} shows that the
deformation-obstruction complex for $\stackdr{X/S(\log\xi)}{n}$ is
just the complex 
\begin{equation} \label{logarithmic-complex}
\xymatrix@1{\op{End}(F_{\op{un}}) \ar[r]^-{\op{ad}_{\nabla_{\op{un}}}} 
& \op{End}(F_{\op{un}})\otimes \op{pr}_{X}^{*}\Omega^{1}_{f}(\log
\xi(S))
}.
\end{equation}
Consequently if we put $\fdr : \stackdr{X/S(\log\xi)}{n}\times_{S} X \to
\stackdr{X/S(\log\xi)}{n}$ for the natural projection we can
characterize $\stackdrreg{X/S(\log\xi)}{n}$ as the open sub stack over
which the morphism of coherent sheaves 
\[
\xymatrix@1{R^{1}{\fdr}_{*}\op{End}(F_{\op{un}})
\ar[rrr]^-{R^{1}{\fdr}_{*}\op{ad}_{\nabla_{\op{un}}}} & & &
R^{1}{\fdr}_{*}(\op{End}(F_{\op{un}})\otimes 
\op{pr}_{X}^{*}\Omega^{1}_{f}(\log
\xi(S)))
}
\]
is surjective. In particular $\stackdrreg{X/S(\log\xi)}{n}$ is smooth
of dimension $(2g - 1)n^{2}$ over $S$. Moreover if $T
\to S$ is an  $S$-scheme and if $(F,\nabla_{f_{T}})$ is in
$\stackdr{X/S(\log\xi)}{n}(T)$ observe that by relative duality
$(F,\nabla_{f_{T}})$ will belong to $\stackdrreg{X/S(\log\xi)}{n}(T)$
iff the morphism 
\[
\op{ad}_{\nabla_{f_{T}}}:
f_{T*}(\op{End}(F)(-\xi_{T}(T)) \to f_{T*}(\op{End}(F)\otimes
\Omega^{1}_{f_{T}})
\]
is injective.

On the other hand if $(F,\nabla_{f_{T}})$ is in $\stackdr{X/S}{n}(T)$
to begin with, then a section $s \in
\Gamma(U,f_{T*}(\op{End}(F)(-\xi_{T}(T)))$ for some open $U \subset S$
will be in the kernel of $\op{ad}_{\nabla_{f_{T}}}$ if and only if $s$
is a $\nabla_{T}$-horizontal section of $\op{End}(F)$ which vanishes
along $\xi_{T}(T)$, i.e. if and only if $s = 0$. This shows that 
$\stackdr{X/S}{n}$ is in fact a sub stack of  $\stackdrreg{X/S(\log \xi)}{n}$.

Observe next  that the family $\stackdr{X/S(\log \xi)}{n} \to S$ also
has an algebraic integrable connection which can be defined in the
same way as the non-abelian Gauss-Manin connection in terms of formal
groupoids \cite[Section~8]{simpson-santa-cruz}. Also the fact that we
have an inclusion of
formal groupoids $X_{\op{DR}} \subset X_{\op{DR}}(\log \xi(S))$ 
(see e.g. Remark~\ref{rem-deRham-groupoid} for notation) implies that
$\stackdr{X/S}{n} \subset \stackdrreg{X/S(\log \xi)}{n}$ is an
inclusion of crystals of stacks. Thus $a_{\op{DR}} : S \to
\stackdrreg{X/S(\log \xi)}{n}$ is a horizontal section.

Finally the fact that (\ref{logarithmic-complex}) is the deformation
obstruction complex for the stack 
$\stackdr{X/S(\log \xi)}{n}$ combined with
the smoothness of $\stackdrreg{X/S(\log \xi)}{n}$ puts us in a
situation where the arguments we used to prove
Proposition~\ref{prop-integrable-section} work verbatim. This shows
that $\nablaf$ comes from a global logarithmic connection $\nabla
: F \to F\otimes \Omega^{1}_{X}(\log \xi(S))$. But $\xi(S)$ is
transversal to the fibers of $f$ and so $\op{Res}(\nabla) =
\op{Res}(\nablaf) = 0$. This proves part (a) of the lemma. 

\medskip

\noindent
Similarly to prove part (b) we have to view the global Higgs field
$(E,\theta)$ corresponding to $(F,\nabla)$ as a section in
$\stackdolreg{X/S(\log \xi)}{n}$. Again we can identify the
deformation obstruction complex of $\stackdolreg{X/S(\log \xi)}{n}$ as
the complex
\[
\xymatrix@1{\op{End}(E_{\op{un}})
\ar[r]^-{\op{ad}_{\theta_{\op{un}}}} & \op{End}(E_{\op{un}})\otimes
\op{pr}_{X}^{*}\Omega^{1}_{f}(\log\xi(S))},
\]
and the same reasoning as in the proof of Lemma~\ref{lem-euler-field}(ii)
and Corollary~\ref{cor-fixed-point} shows that $(E_{s},\theta_{s})$
will be $\cstar$-fixed only when the section
\[
\epsilon_{f} \in H^{0}(S, f_{*}(\op{End}(E)\otimes
\Omega^{1}_{f})/f_{*}\op{End}(E)) \subset  H^{0}(S,
f_{*}(\op{End}(E)\otimes
\Omega^{1}_{f}(\log\xi(S))/f_{*}\op{End}(E))
\]
vanishes at $s$. The lemma is proven. \hfill $\Box$

\

\bigskip

\begin{rem} \label{rem-2-stacks} (i) It is shown in
\cite[Proposition~2.11.2]{beilinson-drinfeld-langlands} that in the
assumptions of the lemma the
morphism $\stackdr{X/S}{n} \to S$ is a l.c.i. morphism of dimension
$(2g -2)n^{2} + n$. 

\medskip

\noindent
(ii) The rigidification trick used in the proof of the previous lemma
is not really necessary and in fact by using crystals of 2-stacks one
should be able to prove Proposition~\ref{prop-integrable-section}(1)
for arbitrary sections $a_{\op{DR}} : S \to \stackdr{X/S}{n}$.
\end{rem}

\

\bigskip
\bigskip

\appendix

\Appendix{Tangent stacks} \label{app-tangent-stacks}

\medskip

\noindent
The language of algebraic stacks is the natural framework for describing 
moduli problems in algebraic geometry. It grew out of M. Artin's
approach to moduli \cite{artin-algebraization1, artin-algebraization2} and 
is by now a standard tool in deformation theory. Since the only comprehensive
treatment of the theory of Artin algebraic stacks is the Orsay preprint
\cite{laumon-stacks} we will review briefly the definition and the basic 
properties of the tangent stack of an algebraic stack. Our main references are 
\cite{artin-stacks}, \cite{laumon-stacks} and \cite[Appendix]{vistoli-stacks}.

\subsection{Algebraic stacks} \label{app-astacks}

The main problem one encounters in constructing a moduli space parameterizing
a given family of geometric objects is the problem of represenability. Very
often the set of equivalence classes of our objects is too wild and does not 
carry any natural geometric structure. Typically the main obstacle for finding
such a structure is the different size of the equivalence classes. To  
remedy that one tries to retain somehow the information about the many 
representatives of a given equivalence class. Thus one is naturally lead to 
replace the set of equivalence classes by the category of all of their 
representatives. This category has a rather special nature since the morphisms
between any two objects come from an equivalence relation and are therefore all
invertible. Categories for which all morphisms are invertible are called
{\em groupoids}. In Artin's approach to moduli the first step is to try and 
represent a given moduli problem not by a scheme but by a category comprised
of groupoids endowed with extra geometric structure. When the formal properties
of such a geometric structure are written down one gets the notion of an 
algebraic (or Artin) stack. This is very similar to the process of putting a
scheme structure on set with the key difference that the points in the 
stack are objects in a category (and hence have intrinsically defined 
automorphisms) rather than elements in a set. The properties 
one needs in order to do geometry on a groupoid had crystallized in the
seminal works of Mumford \cite{mumford-picard}, Deligne and Mumford 
\cite{deligne-mumford} and Artin \cite{artin-stacks}. Even though the actual 
definition given below is rather formal, it is concrete enough to allow us to
operate with an algebraic stack in the same way as with any other object in
algebraic geometry.

\medskip

We will need the notion of a stack over a base scheme $S$. Denote by 
$(\op{Sch}/S)$ the category of schemes over $S$. As explained above intuitively
one should think of a stack as a collection of groupoids endowed with
geometric
\glosstex(categories)[p]{schemes}
structure. All groupoids form a $2$-category $(\op{Grp})$.
\glosstex(categories)[p]{groupoids}
The objects of 
$(\op{Grp})$ are the groupoids, the $1$-morphisms are the functors between 
groupoids and the $2$-morphisms are the  isomorphisms of functors
between groupoids. 

\begin{defi} A groupoid over $S$ (or a pre-sheaf of groupoids over $S$) is a
lax functor 
\[
\mycal{X} : (\op{Sch}/S)^{\op{op}} \to (\op{Grp}).
\]
In other 
words $\mycal{X}$ assigns a groupoid $\mycal{X}_{U}$ to any $S$-scheme $U$
(the groupoid of sections of $\mycal{X}$ over $U$),
a change of base functor $\varphi^{*} : \mycal{X}_{V} \to \mycal{X}_{U}$
to any morphism $\varphi : U \to V$ of $S$-schemes and a canonical isomorphism
of functors 
\[
\psi^{*}\circ \varphi^{*} \cong (\varphi\circ \psi)^{*}
\] 
for any $W \stackrel{\psi}{\to} V$ and $V \stackrel{\psi}{\to} U$ - composable
arrows in $(\op{Sch}/S)$. Furthermore these canonical isomorphisms have to 
satisfy the standard cocycle condition.
\end{defi}

Some good examples to keep in mind are:

\begin{ex} \label{ex-groupoids}
(i) Every algebraic space (cf. \cite{knutson})  $X$ over $S$ has an 
associated $S$-groupoid
which we will again denote by $X$. For any $U \in \op{Ob}(\op{Sch}/S)$
the groupoid $X_{U}$ is just the discrete category corresponding to the 
set of $U$ points $\op{Hom}(U,X)$ of $X$. For any morphism $V 
\stackrel{\varphi}{\to} U$ the change of base morphism $\varphi^{*}$
is just the restriction $\op{Hom}(U,X) \to \op{Hom}(V,X)$.

\medskip

\noindent
(ii) Let $G$ be an affine group scheme over $S$ Let $X$ be an algebraic space
over $S$ equipped with a $G$-action. We have the quotient groupoid 
$[X/G] : (\op{Sch}/S)^{\op{op}} \to (\op{Grp})$ for which $[X/G]_{U}$ is 
the groupoid of all pairs $(P,a)$ where $P$ is a $G$-torsor on $U$ and 
$a : P \to X\times_{S} U$ is a $G$-equivariant morphism. 

As special case is to take $X = S$ equipped with the trivial $G$ action. The
quotient $[S/G]$ is just the groupoid of all $G$-torsors over $S$. It  is 
called the {\em classifying groupoid} of $G$ and is usually denoted by 
$BG$. 
\glosstex(categories)[p]{classifying-groupoid}
\medskip

\noindent
(iii) Let $X$ be an algebraic space over $S$. The $S$-groupoid $\op{Qcoh}_{X/S}$
of quasi-coherent ${\mathcal O}_{X}$-modules is defined as follows. For any
$U \to S$ the groupoid of sections of $\op{Qcoh}_{X/S}$
\glosstex(categories)[p]{qcsheaves-groupoid}
over $U$ is the category whose objects are all quasi-coherent sheaves on 
$X\times_{S} U$ which are flat over $U$ and whose morphisms are the 
isomorphisms of quasi-coherent sheaves. For any morphism $V 
\stackrel{\varphi}{\to} U$ in $(\op{Sch}/S)$ the change of base functor 
$\varphi^{*}$ is just the pull-back via $\op{id}_{X}\times_{S} \varphi$.

\medskip

\noindent
(iv) Fix a reductive group $G$. Let $f : X \to S$ be smooth 
projective and let $x : S \to X$ be a section. Denote by $\rb{X/S}{x}{G}$
the family of representation spaces of the fundamental groups of the 
fibers of $f$. In other words $\rb{X/S}{x}{G}_{s} := 
\op{Hom}(\pi_{1}(X_{s}, x(s)), G)$. Let $\rdr{X/S}{x}{G}$ denote the
fine moduli
scheme of principal $G$-bundles on $X$ with  a relative integrable connection 
and a frame over $x$ constructed in \cite{simpson-moduli1}. Finally put 
$\rdol{X/S}{x}{G}$ for the fine moduli scheme of relative semistable principal
Higgs bundles with vanishing rational Chern classes and a frame over $x$
constructed in \cite{simpson-moduli1}. The group $G$ acts on all of these
spaces and by taking quotients we get pre-sheaves of groupoids
$\stackb{X/S}{G} = [\rb{X/S}{x}{G}/G]$, $\stackdr{X/S}{G} = 
[\rdr{X/S}{x}{G}/G]$ and $\stackdol{X/S}{G} = [\rdol{X/S}{x}{G}/G]$ 
corresponding to the moduli problems for representations, local systems and
Higgs bundles respectively.

\end{ex}

\medskip

In general, when dealing with moduli, one starts with some class ${\mathfrak
X}$ of geometric
objects over $S$ (e.g. schemes, sheaves, maps, etc.) and an  equivalence 
relation ``$\sim$'' on  ${\mathfrak X}$. Next one tries to represent 
the functor to sets  
\glosstex(categories)[p]{sets}
\[
\xymatrix@R=1pt{
(\op{Sch}/S)^{\op{op}} \ar[r]^-{X^{\natural}} & (\op{Set}) \\
 (U \to S) \ar[r] & {\left\{ \begin{minipage}[c]{2.5in} the set of equivalence
classes of families in ${\mathfrak X}$ parameterized by $U$
\end{minipage}\right\}}
}
\]
by a $S$-scheme. This usually fails since $X^{\natural}$ is rarely
a sheaf in any reasonable topology and so cannot be representable. On the
other hand, very often the pre-sheaf of groupoids
\[
\xymatrix@R=1pt{
(\op{Sch}/S)^{\op{op}} \ar[r]^-{\mycal{X}} & (\op{Grp}) \\
 (U \to S) \ar[r] & {\left\{ \begin{minipage}[c]{2.5in} the groupoid with
objects - the families in  ${\mathfrak X}$ parameterized by $U$ and
with 
morphisms - the equivalences of families
\end{minipage}\right\}}
}
\]
is a sheaf. Moreover $X^{\natural}$ is easily recovered from 
$\mycal{X}$ since for every $U \in \op{Ob} (\op{Sch}/S)$ the set 
$X^{\natural}(U)$ is just the set of connected components of the 
groupoid $\mycal{X}(U)$. For example the moduli functors 
$M^{\natural}_{\op{B}}(X/S,n)$, $M^{\natural}_{
\op{DR}}(X/S,n)$, $M^{\natural}_{\op{Dol}}(X/S,n)$ from 
\cite{simpson-moduli1} are obtained in this way from the pre-sheaves of 
groupoids $\stackb{X/S}{G}$, $\stackdr{X/S}{G}$ and $\stackdol{X/S}{G}$ 
respectively. The general principle is that instead of 
trying to represent $X^{\natural}$ by a scheme we may try to 
put enough geometric structure on $\mycal{X}$ so that it can be treated
as a scheme.

Roughly speaking the stacks are pre-sheaves of groupoids which become sheaves 
when considered in an appropriate topology.

\begin{defi} Let $\mycal{X}$ be a groupoid over $S$. The pre-sheaf $\mycal{X}$ 
is called a {\em stack}  in the fppf/\-smooth\-/\'{e}ta\-le topology if 
\begin{list}{{\em (\roman{inner})}}{\usecounter{inner}}
\item For any $U$ in $(\op{Sch}/S)$ and any two objects $x, y$ in 
$\mycal{X}(U)$
the pre-sheaf
\[
\xymatrix@R=1pt{
(\op{Sch}/U)^{\op{op}} \ar[r]^-{\op{Isom}(x,y)} & (\op{Set}) \\
(V \to U) \ar[r] & {\op{Hom}_{\mycal{X}_{V}}(x_{V},y_{V})}}
\]
is a sheaf in the fppf/smooth/\'{e}tale topology.

\item If $\{ V_{i} \stackrel{\varphi_{i}}{\to} U \}$ is a covering of 
$U \in \op{Ob}(\op{Sch}/S)$ in the fppf/smooth/\'{e}tale topology and 
if $(x_{i}, f_{ij})$ is a descend datum (that is - $x_{i} \in 
\op{Ob} \mycal{X}_{V_{i}}$ and $f_{ji} : x_{i|V_{ji}} \stackrel{\cong}{\to} 
x_{j|V_{ji}}$ are morphisms in $\mycal{X}_{V_{ji}}$ satisfying the 
cocycle condition) relative to $\{ V_{i} 
\stackrel{\varphi_{i}}{\to} U \}$, then $(x_{i}, f_{ij})$ is 
effective, i.e. - there exists an object $x$ in $\mycal{X}_{U}$ 
and isomorphisms $f_{i} : x_{|V_{i}} \stackrel{\cong}{\to} x_{i}$ in 
$\mycal{X}_{V_{i}}$ so that for every $i,j$ one has $f_{j|V_{ji}} = f_{ji}
\circ (f_{i|V_{ji}})$. 

\end{list} \label{def-stack}
\end{defi}

It is not hard to verify that all of the $S$-groupoids in 
example~\ref{ex-groupoids} are actually stacks. Checking that the pre-sheaves 
from Example~\ref{ex-groupoids} (i) (ii) and (iv) are stacks is 
straightforward. The proof that Example~\ref{ex-groupoids} (iii) is a stack 
can be found in \cite[Sections 1.1 and 1.2 of 
Expos\'{e} VIII]{sga1}). 

\medskip

\begin{defi} A {\em morphism} between two stacks over $S$ is just a 
$1$-morphisms of pre-sheaves of groupoids. A morphism $f : \mycal{X} \to 
\mycal{Y}$ is {\em injective} if for every $U \in (\op{Sch}/S)$ the 
functor $f_{U} : \mycal{X}_{U} \to \mycal{Y}_{U}$ is faithful. A
morphism $f : \mycal{X} \to \mycal{Y}$ is {\em surjective} if for every
$U \in (\op{Sch}/S)$ and every $y \in \op{Ob}\mycal{Y}_{U}$ there exists 
a covering 
family $\{V \to U\}$ in $(\op{Sch}/S)_{\op{fppf/smooth/\acute{e}tale}}$ and 
a $x \in 
\mycal{X}_{V}$ for which $f_{V}(x)$ is isomorphic to $y$ in $\mycal{Y}_{V}$.
\end{defi}

\medskip

In order to do geometry on a stack $\mycal{X}$ it is essential to be able 
to patch geometric data that is defined ``locally'' on $\mycal{X}$. For 
this one needs a notion of a fiber product of stacks.
Let $f : \mycal{X} \to \mycal{S}$ and $g : \mycal{Y} \to \mycal{S}$ be two 
morphisms of stacks. The fiber product $\mycal{X}\times_{\mycal{S}} 
\mycal{Y}$ is
the stack defined as follows. For any $U \to S$ in $(\op{Sch}/S)$ the 
objects of $(\mycal{X}\times_{\mycal{S}} \mycal{Y} )(U)$ are triples 
$(x, y, \alpha)$ where $x \in \op{Ob} \mycal{X}(U)$, $y \in \op{Ob} 
\mycal{Y}(U)$ and $\alpha : f_{U}(x) \to g_{U}(y)$ is a morphism in 
$\mycal{S}(U)$. A morphism between two $(x',y',\alpha')$ and 
$(x'', y'', \alpha'')$ in $(\mycal{X}\times_{\mycal{S}} \mycal{Y} )(U)$ is 
a pair $(a,b) \in \op{Hom}_{\mycal{X}(U)}(x',x'')\times 
\op{Hom}_{\mycal{Y}(U)}(y',y'')$ for which $\alpha''\circ f(a) = f(b)\circ
\alpha'$. Finally for every $\varphi : V \to U$ in $(\op{Sch}/S)$ the 
pull-back functor $\varphi^{*}$ is defined component wise on every 
$(x, y, \alpha)$ and every $(a,b)$. 

Once we have the notion of a fiber product in the 2-category of stacks we
can study the local behavior of a morphism. Especially useful are morphisms 
between stacks which on affine scheme patches behave as morphisms of schemes.

\begin{defi}
A morphism of stacks $f : \mycal{X} \to \mycal{Y}$ is called {\em 
representable} if for any scheme $U$ in $(\op{Sch}/S)$ and any morphism 
$U \to \mycal{Y}$ over $S$ the fiber product $\mycal{X}\times_{\mycal{Y}} U$
is equivalent to an algebraic space. 
\end{defi}

\begin{rem} It is not hard to characterize the representable morphisms of 
stacks in purely categorical terms. Since the $S$-groupoids corresponding 
to algebraic spaces are exactly the locally discrete one it is clear that
$f : \mycal{X} \to \mycal{Y}$ is representable if and only if the functor
$f_{U}: \mycal{X}_{U} \to \mycal{Y}_{U}$ is faithful for all $U$ in 
$(\op{Sch}/S)$. Informally $\mycal{X} \to \mycal{Y}$ is representable if 
$\mycal{X}$ is less ``stacky'' than $\mycal{Y}$.
\end{rem}

Due to this definition any fppf local property of morphisms of schemes which
is stable under base change makes sense for representable morphisms of
stacks as well. More precisely if $P$ is such a property we will say that
a representable $f : \mycal{X} \to \mycal{Y}$ has the property $P$ if for every
 $S$-scheme $U$ and every morphism $U \to \mycal{Y}$ the morphism of algebraic
spaces $\mycal{X}\times_{\mycal{Y}} U \to U$ has the property $P$. In 
particular we can speak of $f$ being surjective, universally bijective, 
universally open or closed, separated, quasi-compact, of finite type, flat,
smooth, \'{e}tale, etc.

\begin{ex} From the above remark it is clear that if $H \to G$ is a 
homomorphism of algebraic groups over a field $\thick{k}$, then the induced 
morphism of stacks $BH \to BG$ is representable iff $H \to G$ is a 
monomorphism. In particular for $\ast := \op{Spec}(\thick{k})$ and 
$BG = [\ast/G]$ we have that $\ast \to BG$ is representable and that 
$BG \to \ast$ is not representable.
\end{ex}

Now we are ready to introduce the Artin algebraic stacks. Heuristically these
are stacks that look like schemes if one looks at them from the view point
of the category of schemes.

\begin{lemma-defi} An {\em algebraic (geometric) stack} is a $S$-groupoid 
$\mycal{Z}$ 
such that 
\begin{list}{{\em (\arabic{inner})}}{\usecounter{inner}}
\item $\mycal{Z}$ is a stack in the fppf/smooth/\'{e}tale topology.
\item One of the following equivalent conditions holds
\begin{list}{{\em \alph{rom}.}}{\usecounter{rom}}
\item The diagonal $\Delta_{\mycal{Z}} : \mycal{Z} \to \mycal{Z}\times 
\mycal{Z}$ is representable, separated and quasi-compact.
\item For all $S$-algebraic spaces $X$, $Y$ and all morphisms $X \to 
\mycal{Z}$ and $Y \to \mycal{Z}$ the fiber product $X\times_{\mycal{Z}}Y$ is 
equivalent to an algebraic space over $S$.
\item For all $S$-affine schemes $X$, $Y$ and all morphisms $X \to 
\mycal{Z}$ and $Y \to \mycal{Z}$ the fiber product $X\times_{\mycal{Z}}Y$ is 
equivalent to an algebraic space over $S$.
\end{list}
\item There exists a $S$-algebraic space $Z$ and a smooth surjective
morphism $p: Z \to \mycal{Z}$. The pair $(Z,p)$ is called an {\em atlas} of 
$\mycal{Z}$.
\end{list}
\end{lemma-defi}
{\bf Proof.} \cite[Corollary~2.12]{laumon-stacks} \hfill $\Box$

\begin{rem} (i) An algebraic stack $\mycal{Z}$ is called a Deligne-Mumford 
stack  if it has an atlas $p: Z \to \mycal{Z}$ with $p$ - \'{e}tale and 
surjective.

\medskip

\noindent
(ii) An important theorem of M. Artin weakens considerably part (3) of the 
definition of an algebraic stack. Artin's criterion 
\cite[Theorem~6.1]{artin-stacks} asserts that a pre-sheaf of groupoids 
$\mycal{Z}$ is an algebraic stack if and only if $\mycal{Z}$ satisfies 
(1) and (2) and if there exists an algebraic space $Z$ and a surjective
fppf morphism $p : Z \to \mycal{Z}$.

\medskip

\noindent
(iii) Another remarkable feature of Artin algebraic stacks is that they
admit a rather concrete geometric description as a quotient of an algebraic
space by a smooth equivalence relation 
over $S$. Recall \cite[II, 1.1]{knutson} that an equivalence relation on 
$X$ is given by an algebraic $S$-space $R$ together with a 
monomorphism $\delta : R \to X\times_{S} X$ so that for 
every $S$ scheme $U$ the subset $R(U) \subset X(U)\times X(U)$ is 
the graph of an equivalence relation of sets. The quotient of $X$ by 
the equivalence relation $R$ is by definition the quotient sheaf (in the
fppf/smooth/\'{e}tale topology) of sets on $(\op{Sch}/S)$ for the diagram 
\[\xymatrix@1{
R \ar@<.5ex>[r]^{\op{pr}_{1}\circ \delta}\ar@<-.5ex>[r]_{\op{pr}_{2}\circ 
\delta}  & X.
}
\]
We say that $\xymatrix@1{
R \ar@<.5ex>[r] \ar@<-.5ex>[r]  & X
}$ is a {\em smooth} equivalence relation if the structure morphisms
$\op{pr}_{i}\circ \delta$ are smooth and of finite type. 

Given a smooth equivalence relation  $\xymatrix@1{
R \ar@<.5ex>[r] \ar@<-.5ex>[r]  & X
}$ we can construct not only the quotient sheaf of sets on $(\op{Sch}/S)$ but
a quotient algebraic stack $[X/R]$ as well. For any $U$ in 
$(\op{Sch}/S)$ consider the category $[X/R]'(U)$ whose objects are 
pairs $(V \to U, \alpha)$ where $V \to U$ is a smooth covering in 
$(\op{Sch}/S)$ and $\alpha : (\xymatrix@1{V\times_{U} V \ar@<.5ex>[r] 
\ar@<-.5ex>[r]  & V}) \to (\xymatrix@1{ R \ar@<.5ex>[r] \ar@<-.5ex>[r]  & 
X
})$ is a morphism of equivalence relations. For any two pairs $(V' \to U, 
\alpha')$ and $(V'' \to U, \alpha'')$ define 
\[
\op{Hom}_{[X/R]'(U)}(
(V' \to U, \alpha'), (V'' \to U, \alpha'')) = 
\left\{ 
\begin{minipage}[c]{3in} the set of all isomorphisms
of $f : V' \stackrel{\cong}{\to} V''$ over $U$ for which $\alpha'' = 
(f\times_{U} f, f)\circ \alpha'$. 
\end{minipage} 
\right\}
\]
Finally for any $\varphi : U\to V$ denote by
$\varphi^{*} : [X/R]'(V) \to  [X/R]'(U)$  the natural 
restriction functor. Clearly $[X/R]'$ is a presheaf 
of groupoids on $(\op{Sch}/S)$. In general $[X/R]'$ is not a stack
but only a pre-stack (i.e. satisfies only condition (i) in 
Definition~\ref{def-stack} but not condition (ii)). A straightforward analogue
\cite[Lemma~2.2]{laumon-stacks} of the usual plus construction which 
associates a canonical sheaf to any pre-sheaf allows us to stackify the 
pre-stack $[X/R]'$. Denote the resulting stack by $[X/R]$. It
has an obvious smooth atlas $X \to [X/R]$ and by construction
the two equivalence relations $\xymatrix@1{ X\times_{[X/R]} X 
\ar@<.5ex>[r] \ar@<-.5ex>[r]  & X}$ and $\xymatrix@1{R 
\ar@<.5ex>[r] \ar@<-.5ex>[r] & X}$ are 
canonically isomorphic. In particular $[X/R]$ is an Artin stack.

\medskip

\noindent
(iv) The language of smooth equivalence relations is very convenient for 
expressing how far a given stack is from being a scheme. For example if 
$\xymatrix@1{ R \ar@<-.5ex>[r] \ar@<.5ex>[r] & X}$ is a smooth equivalence 
relation and if $R \to X\times_{[X/R]} X$ is  
unramified, then $[X/R]$ is a Deligne-Mumford stack. If 
$R \to X\times_{[X/R]} X$ is an unramified monomorphism,
then $[X/R]$ is an algebraic space.
\end{rem}

\bigskip

All the standard geometric attributes and properties of schemes carry over 
to the realm of Artin algebraic stacks. In particular we can talk of a stack 
$\mycal{X}$ being locally noetherian, reduced, geometrically unibranched, 
regular, separated, quasi-compact, connected, irreducible etc. For example 
for any algebraic $S$-stack $\mycal{X}$ we have

\medskip

$\lozenge$ {\bf Points of ${\mycal X}$:} A point of ${\mycal X}$ is an 
equivalence class of objects 
\[
\xi \in \left. \left( \coprod_{\begin{minipage}[c]{1in} $K$ - field over $S$ 
\end{minipage}} \op{Ob} \mycal{X}(\op{Spec}(K)) \right)\right/\sim,
\]
where $x_{1} : \op{Spec}(K_{1}) \to \mycal{X}$ and $x_{2} : 
\op{Spec}(K_{2}) \to \mycal{X}$ are considered equivalent if there exists a
common field extension $K_{1} \subset K \supset K_{2}$ so that 
$x_{1|\op{Spec}(K)}$ and $x_{2|\op{Spec}(K)}$ are isomorphic in 
$\mycal{X}(\op{Spec}(K))$. 

\medskip

$\lozenge$ {\bf Dimension of $\mycal{X}$:} Let $\mycal{X}$ be locally 
noetherian and irreducible. Choose an irreducible atlas $X  \to \mycal{X}$ 
and let $R := X\times_{\mycal{X}} X$ be the corresponding smooth equivalence
relation. Define $\dim \mycal{X} := \dim X - \dim (R/X)$. Here the relative 
dimension of $R$ over $X$ makes sense since the structure morphisms 
$\op{pr}_{i}\circ \delta : R \to X$, $i = 1,2$ are both smooth and surjective
and hence $\dim (R/X) = \dim (\op{pr}_{1}\circ \delta) = \dim 
(\op{pr}_{2}\circ \delta) = \dim R - \dim X$. Alternatively we have 
$\dim \mycal{X} = 2 \dim X - \dim R$. It can be checked 
\cite[Lemma~5.18]{laumon-stacks} that this definition is correct and does not 
depend on the choice of the atlas $X \to \mycal{X}$. Observe that according
to this definition the dimension of an Artin stack can be negative. For example
for a group $G$ over a field one has $\dim BG = - \dim G$.

\medskip

$\lozenge$ {\bf Sheaves on $\mycal{X}$:} In order to talk about sheaves 
we will have to define an appropriate site first. The naive approach will
be to try and define a Grothendieck topology on $\mycal{X}$ by taking open
sub stacks as neighborhoods. However since the collection of all open sub stacks
in $\mycal{X}$ forms a 2-category instead of a category it is clear that 
this naive approach cannot work directly. As usual the problem can be resolved 
by taking algebraic spaces as neighborhoods.

With any Artin algebraic stack one associates a site $\mycal{X}_{\op{sm}}$
as follows.

\begin{itemize}
\item The objects of $\mycal{X}_{\op{sm}}$ are all smooth maps 
$U \to \mycal{X}$ where $U$ is an algebraic space over $S$.
\item The morphisms between $U \to \mycal{X}$ and $V \to \mycal{X}$ are 
diagrams of the form 
\[
\xymatrix@C=2pt{ U \ar[rrrr]^-{f} \ar@/_/[ddrr] &&&& V \ar@/^/[ddll]\\
& \ar@{=>}[rr]^-{\alpha} &&& \\
&& {\mycal{X}} &&
}
\]
which are commutative up to a natural transformation $\alpha$.
\item The covering families of a smooth $U \to \mycal{X}$ are families 
of morphisms
\[
\xymatrix@C=2pt{ U_{i} \ar[rrrr] \ar@/_/[ddrr] &&&& U \ar@/^/[ddll]\\
& \ar@{=>}[rr]^-{\alpha_{i}} &&& \\
&& {\mycal{X}} &&
}
\]
such that $U_{i} \to \mycal{X}$ and $U_{i}\to U$ are smooth maps and 
$U = \cup_{i} \op{im}(U_{i})$.
\end{itemize}

\medskip

If $\mycal{X}$ is a Deligne-Mumford stack we can define a site 
$\mycal{X}_{\op{\acute{e}t}}$ in exactly the same way. 

\medskip

\begin{defi} \label{def-sheaf-on-stack} A presheaf of sets (abelian groups, 
etc.) 
on an Artin algebraic stack $\mycal{X}$ is a functor 
${\mathfrak F} : \mycal{X}_{\op{sm}}^{\op{op}}\to (\op{Set})$ 
($(\op{Ab})$, etc.). 
\glosstex(categories)[p]{abelian-groups}
A pre-sheaf  ${\mathfrak F}$ on $\mycal{X}$ is a sheaf
in the smooth topology if for any smooth map $U \to \mycal{X}$ and any 
covering family $\{ U_{i} \to U \}$ in $\mycal{X}_{ \op{sm}}$ the diagram
\[
\xymatrix@1{ {{\mathfrak F}}(U) \ar[r]^-{a} & {\prod_{i} {\mathfrak F}(U_{i}) }
\ar@<.5ex>[r]^-{b} \ar@<-.5ex>[r]_-{c} & {\prod_{i,j} 
{\mathfrak F}(U_{i}\times_{U} U_{j})}
}
\]
is exact in the sense that $a$ is the difference kernel of $b$ and $c$.
\end{defi}

Alternatively if $\mycal{X}$ is presented as $[X/R]$ for some smooth 
equivalence relation $\xymatrix@1{R \ar@<.5ex>[r]^-{s} \ar@<-.5ex>[r]_-{t}
& X}$, then a sheaf on $\mycal{X}$ is the same as a sheaf $F$ on $X$ together 
with an isomorphism $s^{*}F \cong t^{*}F$ satisfying the obvious cocycle 
condition on $X\times_{\mycal{X}} X \times_{\mycal{X}} X$.
 
\begin{ex} (i) To every algebraic space $Y$ one can associate a pre-sheaf
$\op{Hom}(\bullet, Y) : \mycal{X}_{ \op{sm}}^{\op{op}}\to (\op{Set})$,
$(U \to \mycal{X}) \mapsto \op{Hom}(U,Y)$. By faithfully flat descend for
algebraic spaces \cite[II.3]{knutson} this is a sheaf.  For $Y = {\mathbb 
A}^{1}$ this sheaf is denoted by ${\mathcal O}_{\mycal{X}_{  \op{sm}}}$.
Notice that by definition ${\mathcal O}_{\mycal{X}_{  \op{sm}}}(U \to 
\mycal{X}) = {\mathcal O}_{U}$.

\medskip

\noindent
(ii) A sheaf of ${\mathcal O}_{\mycal{X}_{  \op{sm}}}$-modules is a 
sheaf ${\mathfrak F}$ on $\mycal{X}_{  \op{sm}}$ such that 
${\mathfrak F}(U \to \mycal{X})$ is a ${\mathcal O}_{\mycal{X}_{  
\op{sm}}}$-module compatibly with pullbacks. A sheaf of ${\mathcal O}_{
\mycal{X}_{  \op{sm}}}$-modules is {\em quasi-coherent (coherent)} if 
${\mathfrak F}(U \to \mycal{X})$ is quasi-coherent (coherent) for all 
$U \to \mycal{X}$ in $\mycal{X}_{  \op{sm}}$.

\end{ex} 

\subsection{The truncated cotangent complex} \label{appss-truncated-cotangent}

In this section we recall, following \cite[Chapter~9]{laumon-stacks}, the 
definition of a tangent stack of an Artin algebraic stack and its relation 
with the cotangent complex.

\medskip

Consider the functor $(\bullet)[\varepsilon] : (\op{Sch}/S) \to (\op{Sch}/S)$, 
defined by $U[\varepsilon] := U \times \op{Spec}(
\bbc[\varepsilon]/(\varepsilon^{2}))$. Denote by $i : U \hookrightarrow 
U[\varepsilon]$ and $r : U[\varepsilon] \to U$ the canonical closed immersion 
and retraction respectively. 

\begin{defi} \label{def-tangent-stack} The tangent groupoid of an $S$-groupoid
$\mycal{X}$ is the pre-sheaf of groupoids $T_{\mycal{X}/S}$ for which 
$T_{\mycal{X}/S}(U) := \mycal{X}(U[\varepsilon])$ and for any morphism 
$\varphi : V \to U$ in $(\op{Sch}/S)$ the base-change functor $\varphi^{*} :
T_{\mycal{X}/S}(U) \to T_{\mycal{X}/S}(V)$ is just the functor 
$(\varphi[\varepsilon])^{*} : \mycal{X}(U[\varepsilon]) \to 
\mycal{X}(V[\varepsilon])$.
\end{defi}

It is not hard to check \cite[Lemma~9.13]{laumon-stacks} that for 
an Artin stack $\mycal{X}$ the tangent groupoid  will also be an Artin stack.
It is equipped with a canonical structure morphism 
$\pr_{\mycal{X}} : 
T_{\mycal{X}/S} \to \mycal{X}$ defined by $\pr_{\mycal{X}}(\xi) = i^{*}\xi$
and with a vertex morphism $\ver_{\mycal{X}} : \mycal{X} \to 
T_{\mycal{X}/S}$ defined by $\ver_{\mycal{X}}(x) = r^{*}x$. Furthermore
there is a natural morphism of $\mycal{X}$-stacks $\gamma : 
{\mathbb A}^{1}\times T_{\mycal{X}/S} \to T_{\mycal{X}/S}$ which can be 
described as follows. For any $\lambda\in {\mathbb A}^{1}$ consider 
the translation
action $t_{\lambda} : \op{Spec}(\bbc[\varepsilon]/(\varepsilon^{2})) \to 
\op{Spec}(\bbc[\varepsilon]/(\varepsilon^{2}))$ induced by the morphism
of rings $a + b\varepsilon \mapsto a + \lambda b\varepsilon$. For any $U$ in 
$(\op{Sch}/S)$ then define the functor $\gamma(\lambda,\bullet) : 
\mycal{X}(U[\varepsilon]) \to \mycal{X}(U[\varepsilon])$ to be the base-change 
functor $(\op{id}_{U}\times t_{\lambda})^{*}$. We leave it to the reader to 
check that $\gamma$ preserves $\ver_{\mycal{X}}$ and is multiplicative up 
to canonical 2-morphisms and thus makes $(T_{\mycal{X}/S},\ver_{\mycal{X}})$
into a cone stack in the sense of \cite[Definition~1.5]{behrend-fantechi}.

\begin{rem} \label{rem-tangent-presentations} 
(i) If $\mycal{X}$ is a smooth algebraic stack presented as 
$[X/R]$ with $X$, $R$ being smooth algebraic spaces, then $\xymatrix@1{T_{R} 
\ar@<.5ex>[r] \ar@<-.5ex>[r] & T_{X}}$ is a presentation of $T_{\mycal{X}}$.  

\medskip 

\noindent
(ii) For a Deligne-Mumford stack $\mycal{X}$ we can define a sheaf 
of K\"{a}hler differentials $\Omega^{1}_{\mycal{X}/S}$ on $\mycal{X}_{  
\acute{e}t}$ by setting $\Omega^{1}_{\mycal{X}/S}(U \to \mycal{X}) = 
\Omega^{1}_{U/S}$. It is clear from the 
definition that we have $T_{\mycal{X}/S} = \op{Spec}(\op{Sym}^{\bullet}
\Omega^{1}_{\mycal{X}/S})$ where $\op{Sym}^{\bullet}
\Omega^{1}_{\mycal{X}/S}$ is the symmetric algebra of $\Omega^{1}_{
\mycal{X}/S}$. In other words in this case $T_{\mycal{X}/S}$ is an abelian 
cone stack in the sense of \cite[Definition~1.9]{behrend-fantechi}. If in
addition $\mycal{X}$ is smooth $T_{\mycal{X}/S}$ will be a vector bundle
stack.

\medskip

\noindent
(iii) Unfortunately, the construction in (ii) cannot be applied
directly to Artin stacks since in that case the naive definition of 
K\"{a}hler differentials used for Deligne-Mumford stacks does not work.
It turns out that for a general Artin stack $\mycal{X}$ the tangent stack
$T_{\mycal{X}/S}$ is again an abelian cone stack. However it is very
rare for $T_{\mycal{X}/S}$ to be a vector bundle stack even when
$\mycal{X}$ is smooth over $S$. Nevertheless $T_{\mycal{X}/S}$ admits
an interpretation in terms of sheaves of differentials similar to the
one in (ii). 

Before we briefly explain this interpretation
 (see \cite[Theorem~9.20]{laumon-stacks}
for more details), recall the following
construction. Given  any 
algebraic $S$-stack $\mycal{X}$ and any complex of sheaves of abelian groups 
$E^{0} \to E^{1}$ on $\mycal{X}$ one may consider the stack theoretic 
quotient of the translation action of $E^{0}$ on $E^{1}$. In this way we get
an $S$-stack $\hh(E^{\bullet}) := [E^{1}/E^{0}]$ having also a structure 
of a strictly commutative group stack over $\mycal{X}$ (see \cite[Section 
1.4 of Expos\'{e} XVIII]{sga4} and \cite[Section 2]{behrend-fantechi} for
details). 

Consider now an atlas $p : X \to \mycal{X}$ for $\mycal{X}$. Due to
the base change property of the K\"{a}hler differentials we get a
complex $\Omega^{1}_{X/S} \to \Omega^{1}_{X/\mycal{X}}$ of
quasi-coherent \'{e}tale sheaves on $X$. It is straightforward to
check that $\hh((\Omega^{1}_{X/S} \to
\Omega^{1}_{X/\mycal{X}})^{\vee})$ is just the stack quotient of
$T_{X/S}$ by the equivalence relation 
\[
\xymatrix@1{T_{X/S}\times_{X} T_{X/\mycal{X}} 
\ar@<.5ex>[r] \ar@<-.5ex>[r]
& T_{X/S}}.
\]
Combined with (i) this now gives a canonical 1-isomorphism
of the pullback $T_{\mycal{X}/S}\times_{\mycal{X}} X$ of
$T_{\mycal{X}/S}$ to $X$ and  
$\hh((\Omega^{1}_{X/S} \to \Omega^{1}_{X/\mycal{X}})^{\vee})$.

\end{rem}

\medskip

If one wants to go one step further and obtain an intrinsic description
of $T_{\mycal{X}/S}$ in terms of differentials on $\mycal{X}$
only, then one is naturally lead to using the cotangent complex of
$\mycal{X}$. 

\subsubsection*{Review of the cotangent complex} \label{sssec-cotangent}

Recall that to each morphism of schemes $f : X \to Y$ Illusie
\cite{illusie-cotangent} associates a canonical chain complex 
$L_{X/Y}  \in \op{Ob}(C^{[-\infty,0]}({\mathcal
O}_{X_{fppf}}))$ called the cotangent complex of $f$. The complex
$L_{X/Y} $ has quasi-coherent cohomology sheaves and is
augmented to $\Omega^{1}_{X/Y}$. The construction of
$L_{X/Y} $ is technical and requires a somewhat
advanced simplicial machinery.  Rather than recalling this elaborate
construction we just list those characteristic properties of the
cotangent complex that are relevant to our discussion.

\bigskip

\centerline{\sc Characteristic properties of the cotangent complex -
the case of schemes}

\medskip

$\bullet$ {\bf Functoriality:} $L_{X/Y} $ exhibits the
same functorial behavior as $\Omega^{1}_{X/Y}$. More precisely any
commutative diagram
\[
\xymatrix{ X' \ar[r]^-{g} \ar[d] & X \ar[d]^-{f} \\
Y' \ar[r] & Y
}
\]
gives rise to a map of complexes $g^{*}L_{X/Y} \to L_{X'/Y'}$.
If in addition $X$ or $Y'$ is flat over $Y$ then this canonical map
is a quasi-isomorphism.

Furthermore, for any morphisms of schemes $X \stackrel{f}{\to} Y \to
S$ the natural short exact sequence of complexes 
$f^{*}L_{Y/S} \to L_{X/S} \to L_{X/Y}$ 
extends to a distinguished triangle in $D({\mathcal
O}_{X_{fppf}})$:
\begin{equation} \label{transitivity-triangle}
f^{*}L_{Y/S} \to L_{X/S} \to L_{X/Y} \to
f^{*}L_{Y/S}[1]
\end{equation}
which depends functorially on $X \to Y \to S$. We will denote degree one map 
$L_{X/Y} \to f^{*}L_{Y/S}[1]$ in
(\ref{transitivity-triangle}) by $e_{S}(X/Y)$. The element
\[
e_{S}(X/Y) \in \op{Hom}_{D({\mathcal
O}_{X_{fppf}}))}(L_{X/Y},f^{*}L_{Y/S}[1]) =:
\op{Ext}^{1}_{{\mathcal
O}_{X}}(L_{X/Y},f^{*}L_{Y/S})
\] 
is called {\em the
Kodaira-Spencer class} of the morphism $f$. In the case when $f$ is
smooth $e_{S}(X/Y)$ coincides with the Atiyah class we introduced after 
Definition~\ref{def-connection-schemes}.

\medskip

$\bullet$ {\bf Relation to the K\"{a}hler differentials:} The
cotangent complex captures the local properties of $f$. For example if
$Y$ is noetherian and $f$ is locally of finite type then $f$ 
is smooth iff the augmentation map $L_{X/Y}  \to
\Omega^{1}_{X/Y}$ is a quasi-isomorphism and $f$ is a l.c.i. morphism
iff $L_{X/Y} $ is quasi-isomorphic to a complex of perfect
amplitude one. 

More generally let $f : X \to S$ be a morphism of schemes which admits a
factorization 
\[
\xymatrix{ X \ar[r]^-{i} \ar[rd]_-{f} & Y \ar[d]^-{p} \\ & S
}
\]
with $i$ a closed immersion defined by an ideal $I$ and $p$ smooth.
Then, in $D(X)$ there is a canonical isomorphism 
\[
\tau_{\geq -1}
L_{X/S} \simeq [ 0 \to I/I^{2} \stackrel{d}{\to} i^{*}
\Omega^{1}_{Y/S} \to 0]
\]
where $i^{*}\Omega^{1}_{Y/S}$ is placed in degree zero.

\medskip

$\bullet$ {\bf Relation to deformation theory:} If $X \to S$ is a
morphism of schemes, the cotangent complex $L_{X/S}$ 
governs the deformation theory of $X$ over $S$. 

Recall that a closed
immersion $i: S \to \overline{S}$ of schemes is called a square zero
extension of $S$ by a quasi-coherent sheaf $M$ if the ideal sheaf $I$
of $S$ in $\overline{S}$ is isomorphic to $M$ as a ${\mathcal O}_{S}$-module
and $I^{2} = 0$ in ${\mathcal O}_{\overline{S}}$.

For any square zero extension $i : S \to \overline{S}$ we have a
functoriality morphism $L_{S} \to 
L_{S/\overline{S}}$ which can be composed with the
truncation $L_{S/\overline{S}} \to \tau_{\geq -1}
L_{S/\overline{S}}$ to give a map $e(i) : L_{S}
\to \tau_{\geq -1}L_{S/\overline{S}} = M[1]$. The role of
the cotangent complex in deformation theory is explained by 
the following theorem

\begin{theo} 
Let $f : X \to S$ be a morphism of schemes and let $i : S \to
\overline{S}$ be a square zero extension of $S$ by a sheaf $M$. Then:
\begin{list}{{\em (\roman{inner})}}{\usecounter{inner}}
\item There exists an obstruction $\omega(X/S,i) \in
\op{Ext}^{2}(L_{X/S}, f^{*}M)$ whose vanishing is necessary
and sufficient  for the existence of a deformation $\bar{f} :
\overline{X} \to \overline{S}$ of $f$ over $\overline{S}$. Furthermore
the obstruction class $\omega(X/S,i)$ is the cup product 
\[
\omega(X/S,i) = f^{*}e(i)\cup e(X/S)
\]
of the class $e(i) \in \op{Ext}^{1}(L_{S},M)$ corresponding
to $i : S \to \overline{S}$ and the Kodaira-Spencer class $e(X/S) \in 
\op{Ext}^{1}(L_{X/S},f^{*}L_{S})$.
\item When $\omega(X/S,i) = 0$, the set of isomorphism classes of
deformations $\bar{f}$ is an affine space under
$\op{Ext}^{1}(L_{X/S}, f^{*}M)$ and the automorphism group
of a fixed deformation is canonically isomorphic to 
$\op{Ext}^{0}(L_{X/S}, f^{*}M)$.
\end{list}
\end{theo}
{\bf Proof.} \cite[Chapter~III]{illusie-cotangent} \hfill $\Box$

\bigskip

\bigskip

Laumon and Moret-Bailly extended Illusie theory of the cotangent
complex to the case of algebraic stacks
\cite[Section~9]{laumon-stacks}. To each 1-morphism $f : \mycal{X}
\to \mycal{Y}$ of $S$-stacks they associate a projective system of ind-objects 
$L^{\geq -n}_{\mycal{X}/\mycal{Y}} \in \op{Ob}(D^{[-n,1]}({\mathcal
O}_{\mycal{X}})) \subset \op{Ob}(D^{[-\infty,1]}({\mathcal
O}_{\mycal{X}}))$ such that for every $n$ the morphism $L^{\geq
-n-1}_{\mycal{X}/\mycal{Y}} \to L^{\geq -n}_{\mycal{X}/\mycal{Y}}$
induces an isomorphisms of the truncation $\tau_{\geq -n} L^{\geq
-n-1}_{\mycal{X}/\mycal{Y}}$ with $L^{\geq -n}_{\mycal{X}/\mycal{Y}}$.
The cotangent complex $L_{\mycal{X}/\mycal{Y}}$ of $f$ is defined as
the projective limit of $\{  L^{\geq
-n-1}_{\mycal{X}/\mycal{Y}} \}_{n \geq 0}$ and enjoys properties
paralleling the ones the cotangent complex for schemes has.

\bigskip

\centerline{{\sc Characteristic properties of the cotangent complex -
the case of stacks}}

\bigskip

\noindent
$\lozenge$ {\bf Normalization:} Whenever $\mycal{X}$ and $\mycal{Y}$
are representable by
algebraic spaces $X$ and $Y$ the object $L_{\mycal{X}/\mycal{Y}}$ can
be realized as the projective limit of the system $\{ \tau_{\geq -n}
L_{X/Y} \}_{n \geq 0}$.

\bigskip

\noindent
$\lozenge$ {\bf Functoriality:} Any 2-commutative diagram of 
algebraic stacks
\[
\begin{minipage}[c]{1in}
\[
{\xymatrix{
{\mycal{X}\,}' \drtwocell
\ar[r]^-{g} \ar[d]  & {\mycal{X}} \ar[d]^-{f}
\ar@{}[dl]|<<{\circlearrowleft} \\
{\mycal{Y}\,}' \ar[r] \ar@{}[ur]|<<{\circlearrowright}  & {\mycal{Y}}
}}
\]
\end{minipage}
\qquad \left( \text{or} \quad
\begin{minipage}[c]{1in}
\[
{\xymatrix{
{\mycal{X}\,}' \drtwocell~^{\dir{}}~_{\dir{}}~'{\dir{}}~`{\dir{}}
\ar[r]^-{g} \ar[d]  & {\mycal{X}} \ar[d]^-{f} \\
{\mycal{Y}\,}' \ar[r]   & {\mycal{Y}}
}}
\]
\end{minipage}
\quad \text{for short} \right)
\]
gives rise to a morphism ${\mathbb L}g^{*}L_{\mycal{X}/\mycal{Y}} \to
L_{\mycal{X}'/\mycal{Y}'}$. 
If in addition $\mycal{X}$ or $\mycal{Y}'$ is flat over $\mycal{Y}$
then this canonical map is an isomorphism.

Furthermore, for any morphisms of algebraic $S$-stacks $\mycal{X} 
\stackrel{f}{\to} \mycal{Y} \to \mycal{Z}$ there exists a morphism 
 $L_{\mycal{X}/\mycal{Y}} \to {\mathbb
L}f^{*}L_{\mycal{Y}/\mycal{Z}}[1]$ which fits into a  
distinguished triangle in $D({\mathcal
O}_{\mycal{X}_{fppf}}))$:
\begin{equation} \label{transitivity-triangle-stacks}
{\mathbb L}f^{*}L_{\mycal{Y}/\mycal{Z}} \to L_{\mycal{X}/\mycal{Z}} \to
L_{\mycal{X}/\mycal{Y}}  \to
{\mathbb L}f^{*}L_{\mycal{Y}/\mycal{Z}}[1]
\end{equation}
which depends functorially on $ \mycal{X} \to \mycal{Y} \to \mycal{Z}$.
Again we get a Kodaira-Spencer class $e_{S}(\mycal{X}/\mycal{Y}) \in 
\op{Hom}_{D({\mathcal
O}_{\mycal{X}_{fppf}}))}(L_{\mycal{X}/\mycal{Y}},({\mathbb
L}f^{*}L_{\mycal{Y}/S})[1])$ called {\em the
Kodaira-Spencer class} of the morphism $f$. 

\bigskip

\noindent
$\lozenge$ {\bf Relation to the K\"{a}hler differentials:} The
ind-object $L^{\geq - 0}_{\mycal{X}/\mycal{Y}}$ in
$D^{[-0,0]}({\mathcal O}_{\mycal X})$, that is - in the category of quasi
coherent sheaves on ${\mycal X}$, has an inductive limit which we will
denote by $\Omega^{1}_{\mycal{X}\mycal{Y}}$. 

If the morphism $f : \mycal{X} \to \mycal{Y}$ is smooth, then the
projective system 
\[
L_{\mycal{X}\mycal{Y}} = (\ldots \to L^{\geq - n
-1}_{\mycal{X}/\mycal{Y}}
\to L^{\geq - n}_{\mycal{X}/\mycal{Y}} \to \ldots \to L^{\geq
0}_{\mycal{X}/\mycal{Y}})
\]
is essentially constant and for any atlas $p : X \to \mycal{X}$ of
$\mycal{X}$ there is a canonical isomorphism of ind-objects in 
$D^{[0,1]}({\mathcal O}_{\mycal X})$:
\[
p^{*}L^{\geq 0}_{\mycal{X}/\mycal{Y}} \widetilde{\to}
[\Omega^{1}_{X/\mycal{Y}} \to \Omega^{1}_{X/\mycal{X}}].
\]

If $f : \mycal{X} \to \mycal{Y}$ is a locally finitely presentable
1-morphism of algebraic stacks, then $f$ is smooth if and only for
each $n \geq 0$ the complex $L^{\geq - n}_{\mycal{X}/\mycal{Y}}$ is of
perfect amplitude contained in $[0,1]$.

Finally if the morphism $f : \mycal{X} \to \mycal{Y}$ is
representable and smooth the projective system
$L_{\mycal{X}/\mycal{Y}}$ is essentially constant and is 
represented by the quasi-coherent sheaf
$\Omega^{1}_{\mycal{X}/\mycal{Y}}$ placed at degree 0. In other words
the natural augmentation morphism $L_{\mycal{X}/\mycal{Y}} \to
\Omega^{1}_{\mycal{X}/\mycal{Y}}$ is an isomorphism.

\bigskip

\noindent
The relation between the stacky cotangent complex and the deformation
theory of the corresponding stacks carries over verbatim from the
scheme case.

\bigskip

\noindent 
Finally let us remark that the tangent stack $T_{\mycal{X}/S}$ of an
algebraic stack can be reconstructed from the sheaf
$\Omega^{1}_{\mycal{X}/S}$ of K\"{a}hler differentials or more
precisely from its background ind-object $L^{\geq
0}_{\mycal{X}/S}$. In other words,
there is (cf. \cite[9.22.1]{laumon-stacks} a canonical 1-isomorphism of stacks
\begin{equation}
\label{tangent-stack-isomorphism}
\hh((L^{\geq 0}_{\mycal{X}/S})^{\vee}) \to T_{\mycal{X}/S}
\end{equation}
which when pulled back to an atlas $p : X \to \mycal{X}$ \quad induces the
canonical isomorphism 
\[
\hh((\Omega^{1}_{X/S} \to
\Omega^{1}_{\mycal{X}/S})^{\vee}) \to
T_{\mycal{X}/S}\times_{\mycal{X}} X. 
\]

\

{\sc

\bigskip

\noindent
L. Katzarkov, UC Irvine. \\ lkatzark@math.uci.edu

\bigskip

\noindent
T. Pantev, UPenn, \\
tpantev@math.upenn.edu 
}

\end{document}